\documentclass[11pt, a4paper, twoside]{article}
\usepackage[a4paper, margin=3.75cm]{geometry}
\usepackage{amsmath}
\usepackage{amssymb}
\usepackage{amsthm}
\usepackage{mathrsfs}
\usepackage{graphicx}

\usepackage{titlesec}
\titleformat{\section}{\normalfont\Large\bfseries}{\thesection.}{0.5em}{}
\titleformat{\subsection}[runin]{\normalfont\normalsize\bfseries}{\thesubsection.}{0.5em}{}
\titleformat{\subsubsection}[runin]{\normalfont\normalsize\bfseries}{\thesubsubsection.}{0.5em}{}

\usepackage{fancyhdr}
\usepackage[font=small, width=1.0\textwidth]{caption}

\usepackage[breaklinks]{hyperref}
\hypersetup{
  unicode=true,
  pdfborder={0 0 0},
  pdffitwindow=false,
  pdfstartview={FitH},
  pdfpagelayout={OneColumn},
  pdftitle={Rational Hermite-Pade Approximants vs Pade Approximants. Numerical Results},
  pdfauthor={Egor O. Dobrolyubov, Nikolay R. Ikonomov, Leonid A. Knizhnerman, Sergey P. Suetin},
  pdfsubject={arXiv},
  pdfkeywords={expansion in small parameter, analytic extension, Riemann surface, Nuttall's partition, Pade approximants, Hermite-Pade polynomials, Van der Pol equation, Katz's points},
  allcolors=blue,
  colorlinks=true
}

\def\({\left(}
\def\){\right)}
\def\[{\left[}
\def\]{\right]}
\def\Re{\operatorname{Re}}
\def\Im{\operatorname{Im}}

\def\mdeg{\operatorname{deg}}
\def\const{\operatorname{const}}

\def\supp{\operatorname{supp}}
\def\mcap{\operatorname{cap}}
\def\NN{\mathbb N}

\def\RR{\mathbb R}
\def\CC{\mathbb C}

\def\sM{\mathscr M}
\def\HH{\mathscr H}
\def\RS{\mathfrak R}
\def\mcap{\mathrm{cap}}
\def\nn{\mathbf{n}}

\let\pfi\varphi
\let\eps\varepsilon
\let\leq\leqslant

\let\myh\widehat
\let\myt\widetilde
\let\myo\overline
\def\myj{\mathbf j}
\def\myk{\mathbf k}

\theoremstyle{plain}
\newtheorem{theorem}{Theorem}
\newtheorem*{theoremS}{Stahl Theorem (1985--1986)}

\theoremstyle{definition}

\begin{document}

\title{\Large\bfseries Rational Hermite--Pad\'e Approximants\\vs Pad\'e Approximants. Numerical Results}

\author{Egor O. Dobrolyubov${}^1$, Nikolay R. Ikonomov${}^2$,\\
Leonid A. Knizhnerman${}^3$, Sergey P. Suetin${}^4$\\[2mm]
\normalsize ${}^1$Lomonosov Moscow State University,\\
\normalsize Chemistry Department, Russia,\\
\normalsize Emanuel Institute of Biochemical Physics,\\
\normalsize Russian Academy of Sciences, Russia\\
\normalsize \href{mailto:dobroljubov@phys.chem.msu.ru}{dobroljubov@phys.chem.msu.ru}\\[2mm]
\normalsize ${}^2$Institute of Mathematics and Informatics,\\
\normalsize Bulgarian Academy of Sciences, Bulgaria\\
\normalsize \href{mailto:nikonomov@math.bas.bg}{nikonomov@math.bas.bg}\\[2mm]
\normalsize ${}^3$Marchuk Institute of Numerical Mathematics,\\
\normalsize Russian Academy of Sciences, Russia\\
\normalsize \href{mailto:lknizhnerman@gmail.com}{lknizhnerman@gmail.com}\\[2mm]
\normalsize ${}^4$Steklov Mathematical Institute,\\
\normalsize Russian Academy of Sciences, Russia\\
\normalsize \href{mailto:suetin@mi-ras.ru}{suetin@mi-ras.ru}}

\date{ }
\maketitle

\begin{abstract}
The main purpose of the paper is to present some powerful data on the advantage of the rational approximation procedure based on Hermite--Pad\'e polynomials over the Pad\'e approximation procedure. The first part of the paper is devoted to some numerical examples in this direction. The second part will be devoted to some theoretical results.

In particular, we demonstrate our ideas about the advantage of rational Hermite--Pad\'e approximants over Pad\'e approximants analyzing the analytical structure of the frequency function $\nu$ of the free Van der Pol equation.

Bibliography:~\cite{VDy92}~titles.

\end{abstract}

% \renewcommand{\baselinestretch}{0.75}\normalsize
% \setcounter{tocdepth}{1}
% \tableofcontents
% \renewcommand{\baselinestretch}{1.0}\normalsize

\thispagestyle{empty}

\newpage
\fancyhf{}
\fancyhead[LE]{\small\thepage}
\fancyhead[CE]{\small E. O. Dobrolyubov, N. R. Ikonomov, L. A. Knizhnerman, S. P. Suetin}
\fancyhead[RO]{\small\thepage}
\fancyhead[CO]{\small Rational Hermite--Pad\'e Approximants vs Pad\'e Approximants}

\section{Introduction}\label{s1}

\subsection{Expansions in small parameter and constructive rational approximations.}\label{s1s1}

The diagonal Pad\'e approximant $[n/n]_f$ is the {\it best rational approximant} to a given power series $f$. It is unique and for a multivalued analytic function $f$ with a finite number of branch points in the extended complex plane $\myh{\CC}$ the corresponding sequence $\{[n/n]_f\}$ gives a possibility to recover the analytic properties of $f$ and to reconstruct its values $f(z)$. This is the main content of the seminal Stahl's Theory~\cite{Sta97b}.

In contrast to that, it is little to be known on the convergence properties of the rational functions based of Hermite--Pad\'e polynomials for the pair $f,f^2$ when $f$ is just from the same class which is treated over Stahl's Theory (of course, $f$ should not be a hyperelliptic function).

Nevertheless, it should be noted that during the last ten years there were obtained quite a lot of results on the convergence properties of the rational Hermite--Pad\'e approximants for algebraic functions which form a very natural subclass of the functions considered in Stahl's Theory (see~\cite{KoPaSuCh17}, \cite{Sue18d}, \cite{Sue20}, \cite{IkSu21}, \cite{Kom21}, \cite{Kom21b}, \cite{Sue21}, \cite{Sue22b} and the bibliography therein).
In particular, it was proved in~\cite{KoPaSuCh17} that the so-called ``Nuttall's domain'' is really a domain on the
$m$-sheeted Riemann surface of an algebraic function of degree~$m$. This property was conjectured in 1984 by Nuttall~\cite{Nut84} in connection with the special ``Nuttall's partition'' of a compact Riemann surface into sheets (see also~\cite{IkSu21}).
In 2021 Aleksandr Komlov~\cite{Kom21}, \cite{Kom21b} introduced a very profound generalization of the classical construction of Hermite--Pad\'e polynomials and proved that the corresponding rational functions recover the values of an algebraic function on all Nuttall's sheets but on the ``highest'' one. This result is in a good accordance with Stahl's Theory since for an hyperelliptic function the diagonal Pad\'e approximants recover the values $f(z)$ of the function on the initial sheet only (it is an old tradition to refer to the initial sheet as ``physical'' or ``zero'' sheet, so the sheets of an $m$-sheeted Riemann surface of $f$ are numerated usually as $0$-sheet, $1$-sheet, \dots, $(m-1)$-sheet and the recovered values of $f$ are $f(z^{(0)}),f(z^{(1)}),\dots,f(z^{(m-1)})$).

Thus up to date a lot of useful information has been accumulated about the limit properties of Hermite--Pad\'e polynomials as wells as the corresponding rational functions. Here we demonstrate some of these properties on several significant numerical examples including the frequency function for the free Van der Pol equation. We announce also some results which will be proved in the second part of the paper.

Interest in Pad\'e approximation has increased dramatically
in 1970--1980s due to the requirements of the perturbation theory in physics and computational mechanics. 
The resulting expansion in a small parameter then had to be analyzed in some way and used to evaluate the expanded function outside of the disk of convergence of the expansion.
In other words, the problem of the extension of a power series outside of its disk of convergence should be solved in a way.
The method of Pad\'e approximations turned out to be very popular from this point of view; see~\cite{BaGr96} and also \cite{VDy83} and \cite{VDy92} and the bibliography therein. It should be noted also that in 2012 Antonio Trias~\cite{Tri12} (see also~\cite{Tri18}) has created
the powerful Holomorphic Embedding Load-flow Method (HELM) to solve the power-flow equations of electrical power systems. This method is completely based on Pad\'e approximations and Stahl's Theory.

\subsection{Definition of diagonal Pad\'e approximant.}\label{s1s2}

Let a function $f$ be holomorphic at the point $z=0$, $f\in\HH(0)$, and be given by the expansion in $z$, i.e.,
\begin{equation}
f(z)=\sum_{k=0}^\infty c_kz^k,\quad |z|<R,\quad c_k=c_k(f).
\label{expan}
\end{equation}
For each $n\in\NN_0$ let polynomials $Q_{n}\not\equiv0$ and $P_{n}$ be such that $\mdeg{Q_{n}}\leq n$, $\mdeg{P_{n}}\leq n$ and as $z\to0$ the following relation holds
\begin{equation}
(Q_{n}f-P_{n})(z)=O(z^{2n+1}).
\label{padedef}
\end{equation}
The relation~\eqref{padedef} is equivalent to the next one
\begin{equation}
(Q_{n}f)(z)=P_{n}(z)+0\cdot z^{n+1}+\dots+0\cdot z^{2n}+ O(z^{2n+1}).
\label{padedenom}
\end{equation}
Thus from~\eqref{padedenom} we obtain the following system of $n$ homogenous equations to the $n+1$ unknowns, the coefficients of the polynomials $Q_{n}(z)=\sum\limits_{j=0}^n q_j z^j$, $q_j=q_j(n)$:
\begin{equation}
\sum_{j=0}^n q_j\cdot c_k(z^j f)=0,\quad k=n+1,\dots,2n.
\label{padesystem}
\end{equation}
Since a nontrivial solution of the system~\eqref{padesystem} always exists, we obtain the existence of a polynomial $Q_{n}\not\equiv0$, $\mdeg{Q_{n}}\leq n$, which satisfies the relation~\eqref{padedenom}. Now the corresponding polynomial $P_{n}$, $\mdeg{P_{n}}\leq n$, is determined directly from the polynomial $Q_{n}$:
\begin{equation}
P_{n}(z)=\sum_{k=0}^n c_k(Q_{n}f)z^k.
\label{padenum}
\end{equation}
It follows directly from the above that namely the polynomial $Q_n$ is the main object of the convergence theory for Pad\'e approximants. Indeed, the polynomial $Q_n$ has to be found from the relation~\eqref{padedenom} and when $Q_n$ is known, the polynomial $P_n$ is determined over $Q_n$ by the explicit representation~\eqref{padenum}. This is in a good agreement with Stahl's Theory since under the conditions of Stahl's results the limit zero distributions for both of this polynomials are just the same.

The rational function $[n/n]_f:=P_{n}/Q_{n}$ is uniquely determined by the relation~\eqref{padedef} and is called the {\it diagonal Pad\'e approximant} (PA) of order $n$ to the power series~\eqref{expan}.
Since in ``generic case'' (when $Q_{n}(0)\neq0$)~\eqref{padedef} is equivalent to the relation
\begin{equation}
f(z)-\frac{P_{n}(z)}{Q_{n}(z)}=O(z^{2n+1}),
\label{padegen}
\end{equation}
then the PA $[n/n]_f=P_{n}/Q_{n}$ is the {\it best local rational approximation} to the power series~\eqref{expan} over the class $\mathscr R_{n}$ of the rational functions of degree $\leq{n}$,
\begin{equation*}
\mathscr R_{n}:=\biggl\{r(z): r(z)=\frac{a_nz^n+\dots+a_0}{b_nz^n+\dots+b_0},\quad\sum_{j=0}^n |b_j|\neq0\biggr\}.
\end{equation*}
The main idea of the $n$\,th diagonal PA is to find a rational function with $n$ {\it free poles} which is the best rational approximation from the class $\mathscr R_n$ to a given power series in $z$.

Notice that it follows directly from~\cite{VDy83}, \cite{VDy92}, \cite{ByDu16}, \cite{Tri18}, \cite{KrDoSyPa20}, \cite{ChDoKr22} that the most interesting object in PA theory is the class of multivalued analytic functions with a finite number of singular points, in particular, the class of algebraic functions. This class is just the main object\footnote{In fact Stahl's Theory is much more powerful and is valid for each multivalued analytic function with a set of singular points of zero (logarithmic) capacity.}
 of Stahl's Theory.

The main Stahl's theorem is originally stated for expansions at the infinity point $\zeta=\infty$, The order of the rational function $[n/n]_{f}(z)$ is invariant under the reflection $z=1/\zeta$. So we have in generic case that
\begin{equation}
f_{*}(\zeta)-[n/n]_f(1/\zeta)=O\(\frac1{\zeta^{2n+1}}\),\quad \zeta\to\infty,
\label{padeinf}
\end{equation}
where $f_{*}(\zeta):=f(1/\zeta)$, $[n/n]_f(1/\zeta)=\zeta^nP_n(1/\zeta)/\bigl(\zeta^nQ_n(1/\zeta)\bigr)=p_n(\zeta)/q_n(\zeta)$.

\begin{theoremS}
Let $f_{*}\in\HH(\infty)$ be a germ of a multivalued analytic function with a finite number of singular points. Then there exists a unique (up to a set of zero capacity) domain $D=D(f_{*})$, $\infty\in D\subset\myh{\CC}$, with the following properties:

1) the germ $f_{*}$ extends to $D$ as a meromorphic (i.e., single-valued analytic) function, $f\in\sM(D)$;

2) the boundary $S:=\partial{D}$ of $D$ does not separate the extended complex plane $\myh{\CC}$ and is a union of a finite number of analytic arcs;

3) the sequence $\{[n/n]_f(1/\zeta)\}$ converges in capacity to $f_{*}(z)$ on the compact subsets of $D$ and the rate of convergence is geometrical, i.e., as $n\to\infty$
\begin{equation}
\bigl|f_{*}(\zeta)-[n/n]_f(1/\zeta)\bigr|^{1/n}\overset
\mcap\longrightarrow e^{-2g_D(\zeta,\infty)}<1,
\label{rateconv}
\end{equation}
where $g_D(\zeta,\infty)$ is Green's function for $D$ with the logarithmic singularity at the infinity point;

4) each pole of $f_{*}$ in $D$ of multiplicity $k$ attracts at least $k$ poles of the rational function $[n/n]_f(1/\zeta)$ as $n\to\infty$.
\end{theoremS}

The set $S=S(f_{*})$ is the so-called Stahl's compact set of minimal capacity (see~\cite{Sta97b} and \cite{Sta12}). It solves in the only way (up to a set of zero capacity) the problem of minimal capacity among all admissible (with respect to the given $f_{*}$) compact sets, i.e.,
\begin{equation}
\mcap(S)=\min\{\mcap(K):\infty\notin K, K\,\text{does not separate}\,\myh{\CC}, f_{*}\in\sM(\myh{\CC}\setminus{K})\}.
\label{mincap}
\end{equation}

Also from Stahl's Theory it follows that all but a ``small'' (in fact $o(n)$ as $n\to\infty$) number of zeros of numerator and denominator of $[n/n]_f(1/\zeta)$ attract the compact set $S$ and the corresponding normalized zero-counting measures converge in weak-$*$ topology to the equilibrium (i.e., Robin) measure of~$S$. The compact $S$ itself consists of a finite number of analytic arcs whose endpoints are some of the branch points of $f_{*}$ and the so-called Chebotar\"ev's points (i.e., the points of zero density of Robin's measure of $S$). The total number of the endpoints (the branch points of $f_{*}$ and Chebotar\"ev's points) is an even number, say $2p+2$. Let $P(\zeta)=\prod\limits_{j=1}^{2p+2}(\zeta-\zeta_j)$ be the corresponding monic polynomial with zeros at the endpoints of those analytic arcs. In generic case all the points $\zeta_j$ are pairwise distinct and thus we can introduce a double-sheeted hyperelliptic Riemann surface $\RS_2(w)$ defined by the equation $w^2=P(\zeta)$ which genus equals $p$. When $f_{*}$ is a germ of a hyperelliptic function, it extends from the infinity point $\zeta=\infty$ to the whole Riemann surface $R_2(w)$ as a single-valued meromorphic function. In the general situation, this is no longer the case. The Riemann surface $R_2(w)$ is called the {\it hyperelliptic Riemann surface associated by Stahl} with the given germ $f_{*}$. This surface plays a crucial role in the description of the {\it strong} asymptotics of polynomials $q_n$ and $p_n$, as well as PA; see~\cite{Nut86}, \cite{ApBuMaSu11}, \cite{ApYa15}.The lifting of Stahl's compact set $S$ onto $\RS_2(w)$ gives the boundary between
the\footnote{Recall that it is an old tradition to numerate the sheets starting from $0$-sheet.}
zero sheet $\RS^{(0)}_2(w)$ and the first sheet $\RS^{(1)}_2(w)$. Since there is a one-to-one correspondence between Stahl's domain $D=\myh{\CC}\setminus{S}$ and zero sheet $\RS^{(0)}_2(w)$, usually the identification is made of $D$ and $\RS^{(0)}_2(w)$.

The branch points of $f_{*}$ which belong to Stahl's compact set $S$ and form the endpoints of the corresponding analytic arcs, are usually called as ``{\it active}'' branch points. Since zeros and poles of diagonal PA are attracted to $S$ with density corresponding to Robin's measure of $S$, they mark the active branch points of $f_{*}$ and Chebotar\"ev's points in a very good way.
Indeed, from Stahl's Theory it follows that the structure of Robin's measure is such that at the active branch points it has infinite density and at Chebotar\"ev's points it has zero density.

%%% Fig 1
\begin{figure}[!t]
\begin{center}
\includegraphics[width=0.75\textwidth]{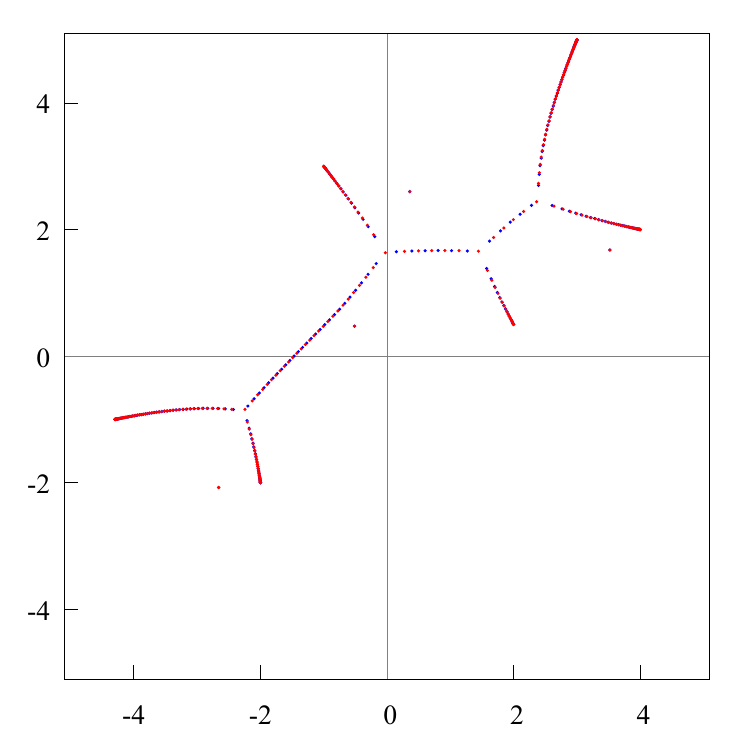}
\vskip-5mm
\caption{The zeros (blue points) and the poles (red points) of PA $[260/260]_f$ for the function $f$ given by~\eqref{fun_six} are plotted after the transformation $z\mapsto\zeta=1/z$. The numerical distributions of these zeros and poles are in full accordance with Stahl's Theory. All six branch points are active branch points and have infinite density of Robin's measure. The corresponding Stahl's compact set $S$ has also four Chebotar\"ev's points with zero density of Robin's measure. There are also four spurious zero-pole pairs which do not correspond to any singularity of $f$.}
\label{fig_six_pade_260}
\end{center}
\end{figure}

For example, let the function $f(z)$ be given by the explicit representation (see~\cite{IkKoSu15})
\begin{equation}
f(z)=\bigl((1-a_1z)(1-a_2z)(1-a_3z)(1-a_4z)(1-a_5z)
(1-a_6z)\bigr)^{1/6}
\label{fun_six}
\end{equation}
with some values $a_1,\dots,a_6$, $f(0)=1$.
In Fig.~\ref{fig_six_pade_260}
the zeros (blue points) and the poles (red points) of Pad\'e approximant $[n/n]_f$ of order
$n=260$ for the function $f$ given by the representation~\eqref{fun_six} are plotted after the transformation $z\mapsto\zeta=1/z$. The numerical distributions of these zeros and poles are in full accordance with Stahl's Theory~\cite{Sta97b}. All six branch points are active branch points and have infinite density of Robin's measure. There are also four Chebotar\"ev's points with zero density of Robin's measure. Notice that there are also four spurious\footnote{Sometimes they are called ``wondering'' pairs or Froissart doublets; see~\cite{ApBuMaSu11} and the bibliography therein.}
zero-pole pairs which do not correspond to any singularity of $f$.

\begin{figure}
\begin{center}
\includegraphics[width=0.7\textwidth]{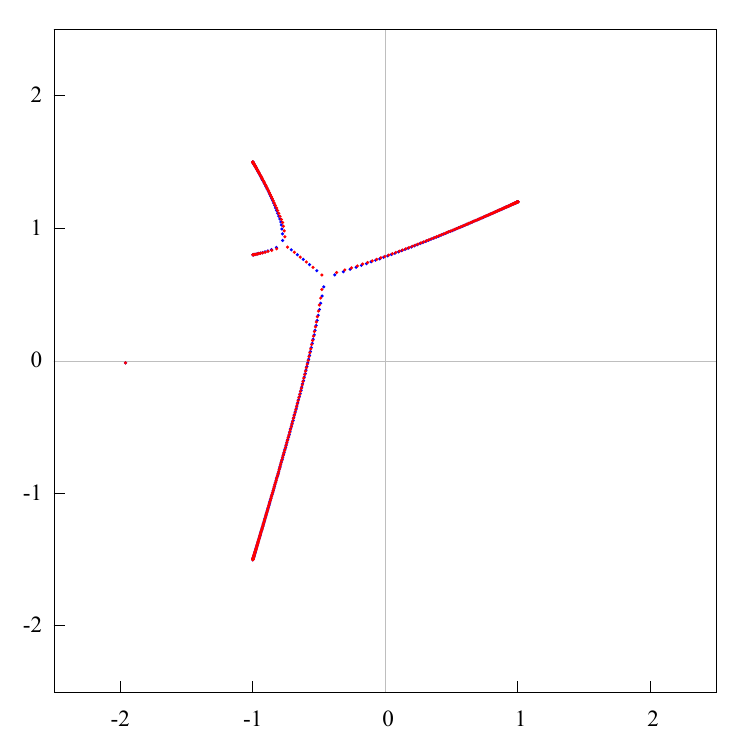}
\vskip-5mm
\caption{(Cf. Fig.~\ref{fig_log_pade_300}) The zeros (blue points) and the poles (red points) of PA $[300/300]_f$ for the function $f$ given by~\eqref{sq} are plotted after the transformation $z\mapsto\zeta=1/z$. The numerical distributions of these zeros and poles are in full accordance with Stahl's Theory. All four branch points are active branch points and have infinite density of Robin's measure. There are also two Chebotar\"ev's points with zero density of Robin's measure. }
\label{fig_sq_pade_300}
% \end{center}
% \end{figure}

% \begin{figure}[!t]
% \begin{center}
\includegraphics[width=0.7\textwidth]{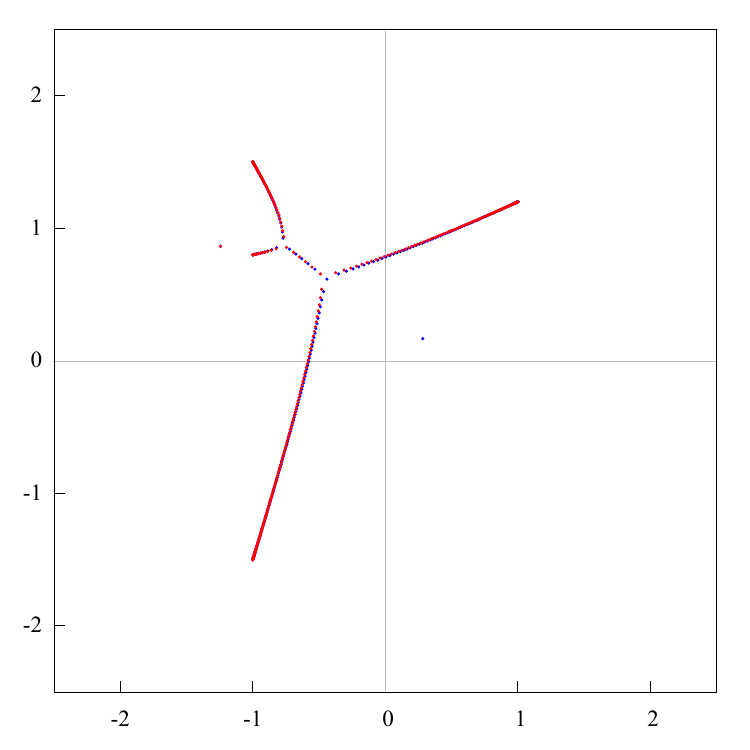}
\vskip-5mm
\caption{(Cf. Fig.~\ref{fig_sq_pade_300}) The zeros (blue points) and the poles (red points) of PA $[300/300]_f$ for the function $f$ given by~\eqref{log} are plotted after the transformation $z\mapsto\zeta=1/z$. The numerical distributions of these zeros and poles are in full accordance with Stahl's Theory. All four branch points are active branch points and have infinite density of Robin's measure. There are also two Chebotar\"ev's points with zero density of Robin's measure.}
\label{fig_log_pade_300}
\end{center}
\end{figure}

Thus over the limit behaviour (as $n\to\infty$) of zeros and poles of PA it is possible to recognize the active (in Stahl's sense) branch points of a function given by a power series~\eqref{expan}. However it is impossible to recognize the type of such branch points based of zeros and poles behaviour. To confirm this observation, let us consider two functions, given by the explicit representations:
\begin{gather}
f_1(z)=\biggl(\frac{1-a_1z}{1-a_2z}\biggr)^{1/2}
+\biggl(\frac{1-a_3z}{1-a_4z}\biggr)^{1/2},
\label{sq}\\
f_2(z)=\log\frac{1-a_1z}{1-a_2z}+\log\frac{1-a_3z}{1-a_4z}
\label{log}
\end{gather}
with some values $a_1,\dots,a_4$.
The types of singularities of these two functions are very different. Nevertheless, from Fig.~\ref{fig_sq_pade_300} and Fig.~\ref{fig_log_pade_300} it follows that Stahl's compact sets $S(f_1)$ and $S(f_2)$ coincide. The joint compact set\footnote{In fact it is a Chebotar\"ev's continuum.}
 $S$ contains all four branch points of each function and two Chebotar\"ev's points as well.

%%% Fig 2-3

Fortunately there is a generalization of Pad\'e construction for the pair of functions $f,f^2$ which gives a possibility to recognize at least square-root singularities via the behaviour of zeros of the corresponding Hermite--Pad\'e polynomials~\cite{Sue22b}. To be more precise, to solve the problem of recognizing the square-root singularities, we should use jointly type II HP polynomials for the pair $f,f^2$ and type I HP polynomials for the tuple $[1,f,f^2]$. Here we will demonstrate this approach numerically on some examples, including the frequency function for the free Van der Pol equation. Note that we also use here for numerical analysis a new construction of HP-polynomials for the tuple $[1,f,f^2,f^3]$ introduced recently by Komlov in~\cite{Kom21}, \cite{Kom21b} (see also~\cite{Sue18b}).

\subsection{Main features of PA and HP-polynomials.}\label{s1s3}

Let us summarize beforehand the main features of PA and HP-polynomials which will be discussed in detail below.

We start with an assumption that we are given $N=2n+1=3m+2$ coefficients of the expansion~\eqref{expan} and we have to use them in an optimal way to recover the properties of the corresponding analytic function. In this context, the paper~\cite{AmBoFe18} of John Boyd and coauthors should be noted. In that paper the authors calculated the type I Hermite--Pad\'e polynomials of degree $m=142$, which corresponds to using $6m+2=854$ terms of the series of the function $u(x;t)$ (see~\eqref{vanexp}) and the corresponding number of terms of the series of the frequency function for the free Van der Pol equation. The coefficients of the series have been calculated in~\cite{AmBoFe18} with $1000$ digits of accuracy. Clearly, the substantial resources were involved to calculate the necessary coefficients. Therefore, the coefficients found must be used optimally. In our opinion, when a fixed number (large enough) of the coefficients is found and someone intends to use those coefficients to analyze the properties of the function, the method of HP-polynomials is much more efficient then the method of PA.

\subsubsection{Pad\'e approximation.}\label{s1s1s1}{\ }

$\bullet$
based on a given $N=2n+1$ coefficients we can find a Pad\'e polynomials $Q_n$ and $P_n$ of degree $\leq{n}$;

$\bullet$
as $N\to\infty$, rational function $[n/n]_f=P_n/Q_n$ extends a given germ $f\in\HH(0)$ as a single valued analytic function into Stahl's domain (with respect to the point $z=0$) $D=\myh{\CC}\setminus{S}$;

$\bullet$
each pole of $f$ in $D$ of multiplicity $k$ attracts at least $k$ poles of the rational function $[n/n]$, i.e. the zeros of $Q_n$ mark the poles of $f$ in $D$ as $n\to\infty$;

$\bullet$
Pad\'e polynomials $Q_n$ and $P_n$ have just the same limit zero distribution which coincide with the equilibrium (Robin) measure of Stahl's compact set $S$;

$\bullet$
for the model function $f_{*}(\zeta)\in\HH(\infty)$ given by the representation~\eqref{funzhu} (see Sec.~\ref{s6}) we have that as $N\to\infty$
(cf.~\eqref{zhuhp2})
\begin{equation}
\bigl|f_*(\zeta)-[n/n]_f(1/\zeta)\bigr|^{1/N}\to
e^{-g_D(\zeta,\infty)}<1,\quad
\zeta\in \myh{\CC}\setminus{S}, \quad S=[-1,1].
\label{zhupade1}
\end{equation}

\subsubsection{Rational Hermite--Pad\'e approximation.}\label{s1s1s2}{\ }

$\bullet$
based on a given $N=3m+2$ coefficients we can find (see Section~\ref{s2}) type I Hermite--Pad\'e polynomials $Q_{m,0},Q_{m,1}$ and $Q_{m,2}$ of degree $\leq{m}$, type II HP-polynomials $P_{2m,0},P_{2m,1},P_{2m,2}$ of degree $\leq{2m}$ as well as the discriminant $D_m:=Q_{m,1}^2-4Q_{m,2}Q_{m,0}$ of degree $\leq{2m}$;

$\bullet$
as $N\to\infty$, the rational function $P_{2m,1}/P_{2m,0}$ extends a given germ $f\in\HH(\infty)$ as a single valued analytic function into the open set $\myh{\CC}\setminus{E}$, $E=[-1,1]$;

$\bullet$
each pole of $f$ in $\myh{\CC}\setminus{E}$ of multiplicity $k$ attracts at least $k$ poles of the rational function $P_{2m,1}/P_{2m,0}$ , i.e. the zeros of $P_{2m,0}$ mark the poles of $f$ in $\myh{\CC}\setminus{E}$;

$\bullet$
as $N\to\infty$,
type I HP-polynomials $Q_{m,0},Q_{m,1}$ and $Q_{m,2}$ have the limit zero distribution which corresponds to the compact set $F$ , $F=[a,b]$, $1<a<b$ (see Sec.~\ref{s6});
type II HP-polynomials $P_{2m,0},P_{2m,1},P_{2m,2}$ have the limit zero distribution which corresponds to the compact set $E$;
the discriminant $D_m$ has the limit zero distribution which corresponds to the compact set $F$; the ordered pair $(E,F)$ forms a so-called Nuttall's condenser~\cite{RaSu13};

$\bullet$
each square-root singularity of $f$ on $E$ attracts at least one zero of the discriminant $D_m$, i.e. the zeros of $D_m$ mark the square-root singularities of $f$;

$\bullet$
for the model function $f_{*}(\zeta)\in\HH(\infty)$ given by the representation~\eqref{funzhu} (see Sec.~\ref{s6}) we have that (cf.~\eqref{zhupade1}) as $N\to\infty$
\begin{equation}
\biggl|f(\zeta)-\frac{P_{2m,1}}{P_{2m,0}}(1/\zeta)\biggr|^{1/N}
\to e^{-\frac13G_F^{\lambda_E}(\zeta)-g_D(\zeta,\infty)}<
e^{-g_D(\zeta,\infty)}, \ \zeta\in \myh{\CC}\setminus{E}, \ E=[-1,1].
\label{zhuhp2}
\end{equation}

Thus, it follows from the above that when we are given $N = 2n + 1 = 3m +2$ coefficients of the power series \eqref{expan}, it is much more efficient to use type I and type II HP-polynomials corresponding to $m$ than to use PA corresponding to~$n$.

Emphasize that we discuss here the {\it rational} approximants based on HP-polynomials, but not the {\it quadratic} Shafer's approximants which are not constructive (see Section~\ref{s2}). In this interpretation we rely on the approach of Peter Henrici~\cite[Sec.~2]{Hen66}:
``\dots a procedure may be
called constructive if it yields the desired mathematical object \dots as the limit of a single sequence of rational functions of the data of the problem~\dots''.

\section{Hermite--Pad\'e polynomials}\label{s2}

\subsection{}\label{s2s1}
Let $f\in\HH(0)$ be given by a representation~\eqref{expan} and
 $f$ is not a hyperelliptic function.

For a fixed $n\in\NN$ let $Q_{n,0},Q_{n,1},Q_{n,2}\not\equiv0$,
be type I HP-polynomials of degree~$\leq{n}$ for the
tuple\footnote{Note
that Pad\'e polynomials $Q_n,-P_n$ could be considered as
Hermite--Pad\'e
polynomials for the tuple of functions $[1,f]$; see~\cite{Nut84}, \cite{Sta88}.}
of functions $[1,f,f^2]$, i.e.
\begin{equation}
(Q_{n,0}+Q_{n,1}f+Q_{n,2}f^2)(z)=O\(z^{3n+2}\),\quad z\to0.
\label{defhp1}
\end{equation}
Let us represent~\eqref{defhp1} in the following form (cf.~\eqref{padedenom}):
\begin{equation}
(Q_{n,1}f+Q_{n,2}f^2)(z)=-Q_{n,0}(z)+0\cdot z^{n+1}+\dots+
0\cdot z^{3n+1}+O(z^{3n+2}).
\label{def2hp1}
\end{equation}
From~\eqref{def2hp1} it follows that similarly to~\eqref{padedenom}
we obtain a system of $2n+1$ homogenous equations with respect to $2n+2$ unknowns, the coefficients of two polynomials $Q_{n,1}$ and $Q_{n,2}$. While these polynomials are found, the polynomial $Q_{n,0}$ can be found similarly to~\eqref{padenum}. Clearly, that to solve the problem~\eqref{defhp1} it should be given $3n+2$ coefficients, $c_0,\dots,c_{3n+1}$, of the expansion~\eqref{expan}.

Now let $P_{2n,0}\not\equiv0,P_{2n,1},P_{2n,2}$ be the type II HP-polynomials of degree $\leq 2n$ for the pair of functions $f,f^2$, i.e. as $z\to0$
\begin{equation}
\begin{aligned}
(P_{2n,0}f-P_{2n,1})(z)&=O(z^{3n+1}),\\
(P_{2n,0}f^2-P_{2n,2})(z)&=O(z^{3n+1}).
\end{aligned}
\label{defhp2}
\end{equation}
Similarly to the definition of PA and type I HP-polynomials, the relations~\eqref{defhp2} are equivalent to the following relations
\begin{equation}
\begin{aligned}
(P_{2n,0}f)(z)&=P_{2n,1}(z)+0\cdot z^{2n+1}+\dots+0\cdot z^{3n}+O(z^{3n+1}),\\
(P_{2n,0}f^2)(z)&=P_{2n,2}(z)
+0\cdot z^{2n+1}+\dots+0\cdot z^{3n}+
O(z^{3n+1}).
\end{aligned}
\label{def2hp2}
\end{equation}
From~\eqref{def2hp2} it follows that similarly to~\eqref{padedenom} and~\eqref{defhp1} we obtain a system of $2n$ homogenous equations with respect to $2n+1$ unknowns, the coefficients of the polynomial $P_{2n,0}$. When this polynomial is found, the polynomials $P_{2n,1}$ and $P_{2n,2}$ can be found similarly to~\eqref{padenum}. Thus just as in the PA case, the polynomial $P_{2n,0}$ is the main object of the construction and the polynomials $P_{2n,1}$ and $P_{2n,2}$ are secondary (cf.~\eqref{padenum}):
\begin{equation}
P_{2n,1}(z)=\sum_{k=0}^{2n}c_k(P_{2n,0}f)z^k, \quad
P_{2n,2}(z)=\sum_{k=0}^{2n}c_k(P_{2n,0}f^2)z^k.
\label{numhp2}
\end{equation}
Clearly, that to solve the problem~\eqref{defhp2} it should be given $3n+1$ coefficients, $c_0,\dots,c_{3n}$, of the expansion~\eqref{expan}.

Thus, we use $2n+1$ coefficients to obtain the core polynomial $Q_n$ of degree $\leq{n}$, i.e. the Pad\'e denominator. And we use $3n+1$ coefficients to obtain the core type II HP-polynomial $P_{2n,0}$ of degree $\leq{2n}$. In the first case the productivity of the procedure equals
$k_{\mathrm{prod, PA}}=n/(2n+1)\approx1/2$, while in the second case it equals $k_{\mathrm{prod, HP}}=2n/(3n+1)\approx2/3>1/2$.
Below we will give some other arguments in favor of using HP polynomials and the corresponding rational approximants instead of PA.

At first we discuss the difference between the limit zero distribution of Pad\'e polynomials $Q_n$ and type II HP-polynomials $P_{2n,0}$ as $n\to\infty$.

Stahl's Theory in fact consists of two parts: the geometric part and the analytical part. In the first, geometric part, it is established that there exists a unique Stahl's compact set $S$ of minimal capacity. In the second, analytic part, it is established that both polynomials $Q_n$ and $P_n$ have just the same limit zero distribution which coincides with Robin's measure of $S$. These facts imply finally the convergence of the diagonal PA in capacity inside (i.e. on the compact subsets) the domain $D=\myh{\CC}\setminus{S}$. Until the geometric part is established, it impossible to establish the analytic part. From this side the situation with HP-polynomials for $f,f^2$ is quite complicated, namely till now it is unclear in what terms the geometrical part should be solved. In particular, in general it is unclear what compact set should be considered instead of Stahl's compact set when describing the limit zero distribution of HP-polynomials $P_{2n,0}$. When $f$ is a cubic function the problem is completely solved in terms of Nuttall's global partition (with respect to expansion point $z=0$ or $\zeta=\infty$) of the corresponding three-sheeted Riemann surface $\RS_3$ into three open sheets; see~\cite{Nut84}, \cite{KoPaSuCh17}. In particular, the zeros of HP-polynomials $P_{2n,0}$ are attracted to the projection onto the complex plane $\myh{\CC}$ of the boundary between the zero sheet and the first sheet of $\RS_3$.
For a general case of two functions $f_1,f_2$ there was introduced in~\cite{RaSu13} a new approach to the problem based on the so-called {\it Nuttall's condenser} $(E,F)$ which replaces for the pair of functions Stahl's compact set $S$. The two plates of the condenser interact to each other. However in~\cite{Sue18b} and \cite{Sue19} it was stated that $F$ is the core plate and the plate $E$ is completely determined by $F$. For some cases it is possible to prove that $E$ coincides with the projection of the boundary between the zero sheet and the first sheet of the corresponding $\RS_3$ and $F$ coincides with the projection of the boundary between the first and the second sheets of $\RS_3$. Therefore it is quite believable that for a pair of functions $f,f^2$ the compact set $E$ is an analog of Stahl's compact set~$S$.

\begin{figure}
\begin{center}
\includegraphics[width=0.75\textwidth]{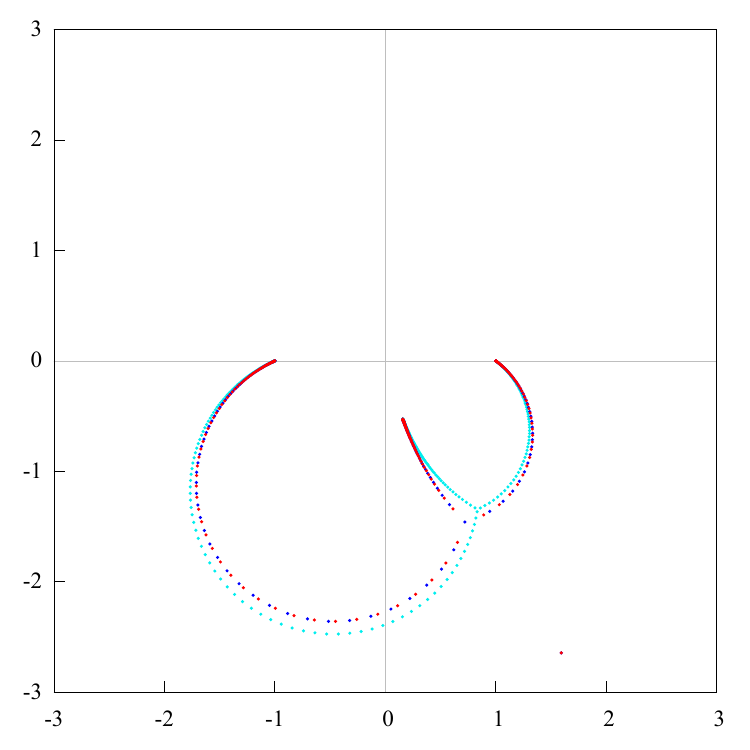}
\vskip-5mm
\caption{The zeros (blue points) and the poles (red points) of PA $[200/200]_f$ as well as the zeros (light blue points) of the type II HP-polynomial $P_{400,0}$ for the function $f$ given by~\eqref{cube2021} are plotted. Evidently the corresponding compact sets $S$ and $E$ are different from each other. Also there are a Chebotar\"ev's point of zero density on $S$ and a Chebotar\"ev's point on $E$ of positive density.}
\label{fig_cube2021_pade_hepas}
% \end{center}
% \end{figure}

% \begin{figure}
% \begin{center}
\includegraphics[width=0.75\textwidth]{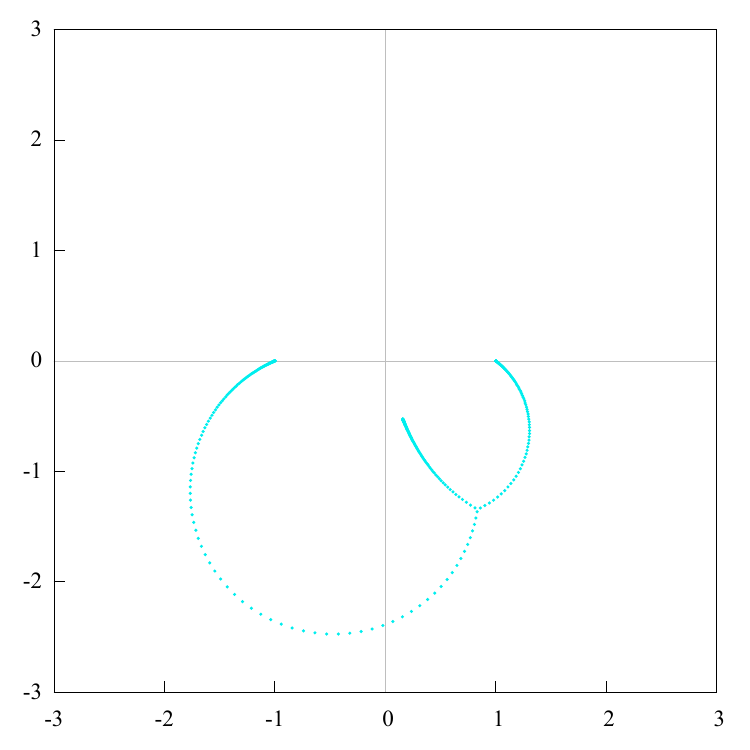}
\vskip-5mm
\caption{The zeros (light blue points) of the type II HP-polynomials $P_{400,0}$ for the function $f$ given by~\eqref{cube2021} are plotted. The corresponding compact set $E$ is a continuum with a unique Chebotar\"ev's point of positive density.}
\label{fig_cube2021_hepas}
\end{center}
\end{figure}

\begin{figure}
\begin{center}
\includegraphics[width=0.75\textwidth]{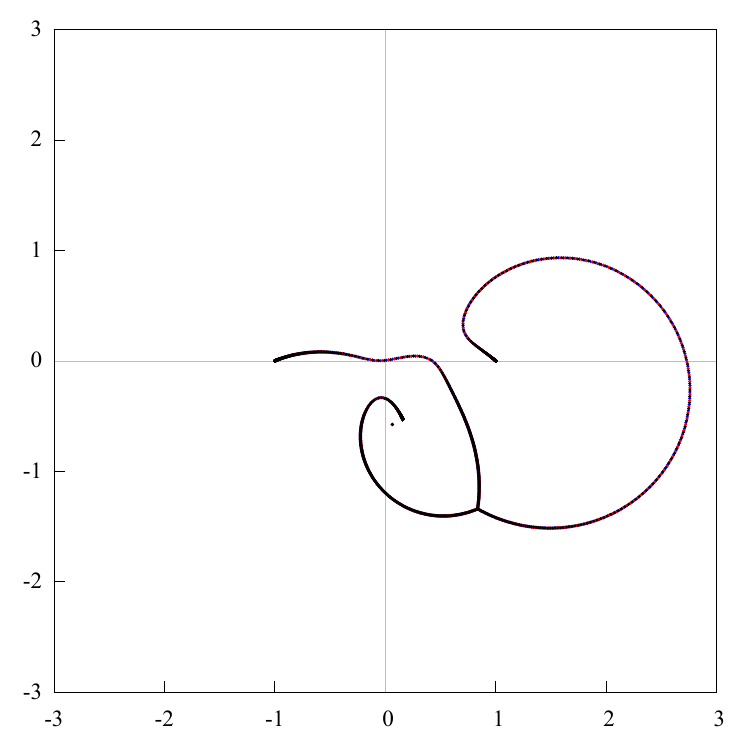}
\vskip-5mm
\caption{The zeros (blue, red and black points) of the type I HP-polynomials $Q_{200,0},Q_{200,1}$ and $Q_{200,2}$ for the function $f$ given by~\eqref{cube2021} are plotted. The corresponding compact set $F$ is a continuum with a unique Chebotar\"ev's point of positive density.}
\label{fig_cube2021_hepa}
% \end{center}
% \end{figure}

% \begin{figure}
% \begin{center}
\includegraphics[width=0.75\textwidth]{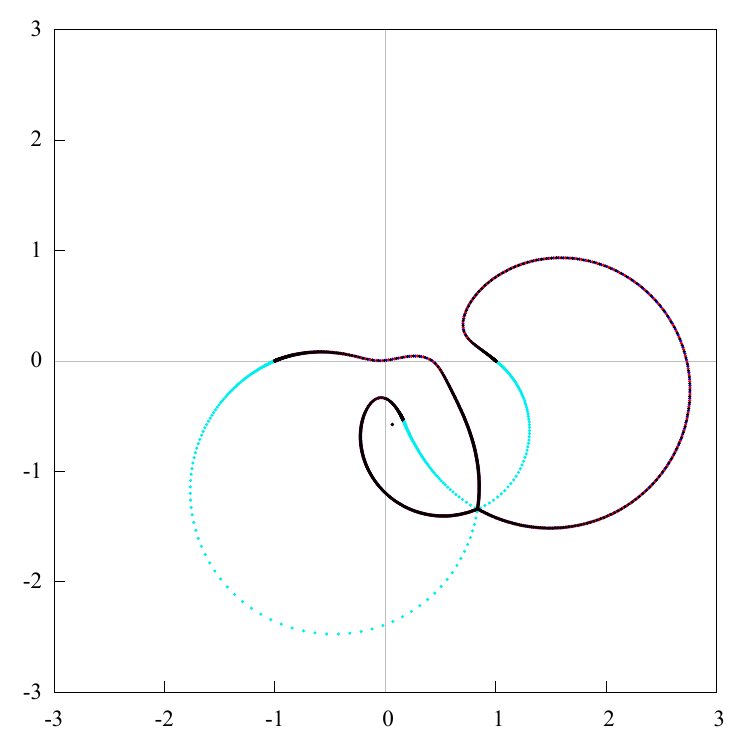}
\vskip-5mm
\caption{We combine the zeros (light blue points) of the type II HP-polynomials $P_{400,0}$ as well as the zeros (blue, red and black points) of the type I HP-polynomials $Q_{200,0},Q_{200,1}$ and $Q_{200,2}$ for the function $f$ given by~\eqref{cube2021}. The Chebotar\"ev's points of $E$ and $F$ coincide.}
\label{fig_cube2021_hepas_hepa}
\end{center}
\end{figure}

We illustrate these facts on a very simple function.
Let a function $f\in\HH(0)$ be given by the explicit representation
\begin{equation}
f(z)=(1-z^2)^{1/3}(1-za)^{-2/3},
\label{cube2021}
\end{equation}
where $a=(0.3+i)\sqrt{3}$, $f(0)=1$.
Thus $f$ is an algebraic function of order three. The corresponding Riemann surface (RS) of $f$, $\mathfrak R_3=\mathfrak R_3(f)$, is a three-sheeted branching covering of the Riemann sphere $\myh\CC$ with the branch points $\Sigma=\{\pm1,1/a\}$ and is of genus one. Due to Nuttall~\cite{Nut84} the RS $\mathfrak R_3$ is divided into three open (nonintersecting) sheets $\mathfrak R_3^{(0)},\mathfrak R_3^{(1)}$ and $\mathfrak R_3^{(2)}$ with
$\Gamma^{(0,1)}=\partial\mathfrak R_3^{(0)}$,
$\Gamma^{(1,2)}=\partial\mathfrak R_3^{(2)}$ and $\partial\mathfrak R_3^{(2)}=\Gamma^{(0,1)}\cup \Gamma^{(1,2)}$. It has been conjectured by Nuttall~\cite{Nut84} in 1984 and proven~\cite{KoPaSuCh17} in 2017 that $\mathfrak R_3\setminus\myo{\RS_3^{(2)}}$ is a domain on $\RS_3$. Let $\pi$ be the canonical projection of $\RS_3$ onto $\myh\CC$, $\pi\colon\RS_3\to\myh\CC$. Set $E=\pi(\Gamma^{(0,1)})$, $F=\pi(\Gamma^{(1,2)})$ (see Fig.~\ref{fig_cube2021_hepas} and Fig.~\ref{fig_cube2021_hepa}). From~\cite{Kom21} it follows that $E\cap F=\{\pm1,1/a,v\}$, where $v$ is a unique Chebotar\"ev's point of positive density for both $E$ and $F$ (see Fig.~\ref{fig_cube2021_hepas_hepa}).
The ordered pair $(E,F)$ forms {\it Nuttall's condenser} for the case under consideration.

%%% Fig 4-7

In Fig.~\ref{fig_cube2021_pade_hepas} the zeros (blue points) and
the poles (red points) of Pad\'e $[200/200]_f$ as well as the zeros (light blue points) of the type II HP-polynomial $P_{400,0}$ are plotted. Evidently the corresponding compact sets $S$ and $E$ are different from each other. Also there is a Chebotar\"ev's point of zero density on $S$, while Chebotar\"ev's point of $E$ is of positive density.

In Fig.~\ref{fig_cube2021_hepas}
the zeros of the type II HP polynomial $P_{2n,0}$ of order
$2n=400$ are plotted (light blue points). The corresponding compact set $E$ is a continuum with a unique Chebotar\"ev's point $v$ of positive density.

In Fig.~\ref{fig_cube2021_hepa}
the zeros of the type I HP polynomials $Q_{n,0},Q_{n,1}$ and $Q_{n,2}$ for $n=200$ are plotted (blue, red and black points). The corresponding compact set $F$ is a continuum with a unique Chebotar\"ev's point $v$ of positive density.

Jointly, the compact sets $E$ and $F$ form Nuttall's condenser, see Fig.~\ref{fig_cube2021_hepas_hepa}.

It is proved (see~\cite{Nut84}, \cite{KoPaSuCh17}) that for the function given by~\eqref{cube2021} as $n\to\infty$,
all but one zeros of $P_{2n,0}$ are attracted to $E$, all but one zeros of $Q_{n,j}$ are attracted to $F$ and for $z\in\myh{\CC}\setminus{E}$ and $z\in\myh{\CC}\setminus{F}$ we have respectively that
\begin{equation}
\frac{P_{2n,1}}{P_{2n,0}}(z)\overset{\mcap}\to f(z^{(0)}),\quad z\notin{E},\quad
-\frac{Q_{n,1}}{Q_{n,2}}(z)
\overset{\mcap}\to f(z^{(0)})+f(z^{(1)}),\quad z\notin{F}.
\label{conv_hp}
\end{equation}
Note that
\begin{equation*}
\frac{P_{2n,2}}{P_{2n,0}}(z)\overset{\mcap}\to f^2(z^{(0)}),\quad z\notin{E},
\end{equation*}
but that gives us no any additional information on the extension of $f\in\HH(0)$.

In conclusion, from~\eqref{conv_hp} it follows that {\it rational} functions based on type I and type II HP-polynomials, jointly give us a possibility to evaluate the function $f$ not only on the zero sheet but also on the first sheet.
In addition the plate $E$ can be recognized via the zeros of the type II HP-polynomials $P_{2n.0}$ and similarly the plate $F$ can be recognized via the zeros of the type I HP-polynomials $Q_{n,j}$.

Recall that the well-known Shafer's quadratic approximants $S_n$ (see~\cite{Sha74}) are determined from the equation (cf.~\eqref{defhp1})
\begin{equation}
Q_{n,0}(z)+Q_{n,1}(z)S_n+Q_{n,2}S_n^2\equiv0.
\end{equation}
Thus
\begin{equation}
S_{n;\pm}(z)=\frac{-Q_{n,1}\pm\sqrt{Q^2_{n,1}(z)-4Q_{n,0}(z)Q_{n,2}(z)}}{2Q_{n,2}(z)}.
\label{defsha}
\end{equation}

According to Henrici's definition~\cite[Sec.~2]{Hen66} of constructive approximation, Shafer's quadratic approximants are not constructive. Nevertheless, the discriminant of \eqref{defsha}
\begin{equation}
D_n(z):=Q^2_{n,1}(z)-4Q_{n,0}(z)Q_{n,2}(z)
\label{defdis}
\end{equation}
is constructive and keeps very useful information about the analytic structure of $f$ (see~\cite{Sue22b}). Note that in generic case $\mdeg{D_n}=2n$.

It is well-known that in a generic case the type I and the type II HP-polynomials are connected with each other~\cite{Nut84}. For example, let us consider a tuple $[1,f,f^2]$. For $n\in\NN$ set $\nn_1=(n,n,n-1)$ and $\nn_2=(n,n-1,n)$ be two multiindexes with $|\nn_1|=|\nn_2|=3n-1$. Let us define the type I HP-polynomials $Q_{\nn_k,j}$, $k=1,2$, $j=0,1,2$, by the relations:
\begin{align}
(Q_{\nn_1,0}+Q_{\nn_1,1}f+Q_{\nn_1,2}f^2)(z)&=O\(z^{3n+1} \),\label{defhp1m1}\\
(Q_{\nn_2,0}+Q_{\nn_2,1}f+Q_{\nn_2,2}f^2)(z)&=O\(z^{3n+1} \),\label{defhp1m2}
\end{align}
where $\mdeg{Q_{\nn_1,0}}=\mdeg{Q_{\nn_1,1}}=n$, $\mdeg{Q_{\nn_1,2}}=n-1$,
$\mdeg{Q_{\nn_2,0}}=\mdeg{Q_{\nn_2,2}}=n$, $\mdeg{Q_{\nn_2,1}}=n-1$.
It is easy to check that to solve~\eqref{defhp1m1} and~\eqref{defhp1m2} we need $3n+1$ coefficients $c_0,\dots,c_{3n}$ of a given power series~\eqref{expan}.
It is also easy to check (see~\cite{NiSo88}) that in a generic case
\begin{equation}
P_{2n,0}=
\begin{vmatrix}
Q_{\nn_1,1} & Q_{\nn_1,2}\\
Q_{\nn_2,1} & Q_{\nn_2,2}
\end{vmatrix}.
\label{hp2hp1}
\end{equation}
Thus using a generalization of Viskovatov's algorithm~\cite{IkSu21} and different multiindexes it is possible based on~\eqref{hp2hp1} to compute the type II HP-polynomials in terms of the type I HP-polynomials.

\subsection{}\label{s2s2}

Now let us consider a fourth-order algebraic function $f\in\HH(0)$. Let $\RS_4=\RS_4(f)$ be the four-sheeted Riemann surface of $f$ and $\RS_4^{(0)},\RS_4^{(1)},\RS_4^{(2)},\RS_4^{(4)}$ be the corresponding Nuttall's partition of $\RS_4$ with respect to the point $z=0$. Let $F^{(0)}$ be the projection of the boundary between the zero and the first sheets, $F^{(1)}$ be the projection of the boundary between the first and the second sheets and $F^{(2)}$ be the projection of the boundary between the second and the third sheets.

For the tuple $[1,f,f^2,f^3]$ and a fixed $n\in\NN$ let us define similarly to~\eqref{defhp1} the type I HP-polynomials of degree $\leq{n}$ by the relation
\begin{equation}
(Q_{n,0}+Q_{n,1}f+Q_{n,2}f^2+Q_{n,3}f^3)(z) = O(z^{4n+3}),\quad z\to0.
\label{defhp1m}
\end{equation}
Note that in fact we defined type I HP-polynomials for a multi-index $\nn=(n,n,n,n)$. Also similarly to~\eqref{defhp2} let us define the type II HP-polynomials of degree $\leq{3n}$ for the triple of the functions $f,f^2,f^3$ in the following way
\begin{equation}
\begin{aligned}
(P_{3n,0}f-P_{3n,1})(z)&=O(z^{4n+1}),\\
(P_{3n,0}f^2-P_{3n,2})(z)&=O(z^{4n+1}),\\
(P_{3n,0}f^3-P_{3n,3})(z)&=O(z^{4n+1}).
\end{aligned}
\label{defhp2m}
\end{equation}

It is proved (see~\cite{Nut84}, \cite{KoPaSuCh17}) that as $n\to\infty$,
all but a finite number of zeros of $P_{3n,0}$ are attracted to $F^{(0)}$, all but a finite number of zeros of $Q_{n,j}$ are attracted to $F^{(2)}$ and for $z\in\myh{\CC}\setminus{F^{(0)}}$ and $z\in\myh{\CC}\setminus{F^{(2)}}$ we have respectively that
\begin{equation}
\frac{P_{3n,1}}{P_{3n,0}}(z)\overset{\mcap}\to f(z^{(0)}),\, z\notin{F^{(0)}},\quad
-\frac{Q_{n,1}}{Q_{n,2}}(z)
\overset{\mcap}\to f(z^{(0)})+f(z^{(1)})+f(z^{(2)}),\, z\notin{F^{(2)}}.
\label{conv_hpm}
\end{equation}
Clearly, the relations~\eqref{conv_hpm} are not completed, since we can't find the values $f(z^{(0)}),f(z^{(1)})$ and $f(z^{(2)})$ from them. Fortunately, in~\cite{Kom21} a new construction of HP-polynomials was introduced to fill this gap in the general case of an $m$-order algebraic function. In the partial case under consideration, this construction leads to a system of polynomials
$H_{2n,j}$ of degree $\leq{2n}$ such that as $n\to\infty$
all but a finite number of zeros of $H_{2n,j}$ are attracted to $F^{(1)}$ and
\begin{equation}
\frac{H_{2n,1}}{H_{2n,0}}(z)
\overset{\mcap}\longrightarrow f(z^{(0)})+f(z^{(1)}),\, z\notin{F^{(1)}}.
\label{conv_hpm2}
\end{equation}
Thus the relation~\eqref{conv_hpm2} complements the relations~\eqref{conv_hpm}.

In conclusion, from~\eqref{conv_hpm} and~\eqref{conv_hpm2} it follows that {\it rational} functions based on HP-polynomials, jointly give us a possibility to evaluate the function $f$ on the zero, first and second sheets of the corresponding RS.
In addition the plates $F^{(0)},F^{(1)}$ and $F^{(2)}$ can be recognized via the zeros of HP-polynomials $P_{3n.0},H_{2n,j}$ and $Q_{n,j}$ respectively. This information will be used below while analyzing the analytic properties of the frequency function of the Van der Pol equation.

Note that to find the polynomials $Q_{n,j}$ satisfying~\eqref{defhp1m} we need $4n+3$ coefficients of the expansion~\eqref{expan} and to find the polynomial $P_{3n,0}$ satisfying~\eqref{defhp2m} we need $4n+1$ coefficients of~\eqref{expan}. Let us suppose that we are given $4n+1$ coefficients of~\eqref{expan}.
Let introduce three multiindexes $\nn_1=(n,n,n-1,n-1)$, $\nn_2=(n,n-1,n,n-1)$, $\nn_3=(n,n-1,n-1,n)$ and the corresponding type I HP polynomials:
\begin{equation}
\begin{aligned}
(Q_{\nn_1,0}+Q_{\nn_1,1}f+Q_{\nn_1,2}f^2+Q_{\nn_1,3}f^3)
(z)&=O\(z^{4n+1}\),\\
(Q_{\nn_2,0}+Q_{\nn_2,1}f+Q_{\nn_2,2}f^2+Q_{\nn_2,3}f^3)
(z)&=O\(z^{4n+1}\),\\
(Q_{\nn_3,0}+Q_{\nn_3,1}f+Q_{\nn_3,2}f^2+Q_{\nn_3,3}f^3)
(z)&=O\(z^{4n+1}\).
\end{aligned}
\label{defhp1mi}
\end{equation}
Then in a generic case (see~\cite{Nut84}, \cite{NiSo88} and cf.~\eqref{hp2hp1})
\begin{equation}
P_{3n,0}=
\begin{vmatrix}
Q_{\nn_1,1} & Q_{\nn_1,2} & Q_{\nn_1,3}\\
Q_{\nn_2,1} & Q_{\nn_2,2} & Q_{\nn_2,3}\\
Q_{\nn_3,1} & Q_{\nn_3,2} & Q_{\nn_3,3}
\end{vmatrix}
\label{hp2hp1m}
\end{equation}
and (see~\cite{Sue18d} and \cite{Sue20})
\begin{equation}
H_{2n,0}=
\begin{vmatrix}
Q_{\nn_1,2}&Q_{\nn_1,3}\\
Q_{\nn_2,2}&Q_{\nn_2,3}
\end{vmatrix},\quad
H_{2n,1}=
\begin{vmatrix}
Q_{\nn_1,1}&Q_{\nn_1,3}\\
Q_{\nn_2,1}&Q_{\nn_2,3}
\end{vmatrix}.
\label{hpnhp1m}
\end{equation}
Thus again based on a generalization of the Viskovatov algorithm~\cite{IkSu21b}, \cite{IkSuHePa21} and different multiindexes it is possible to compute HP-polynomials through type I HP-polynomials.

% \newpage\clearpage
\section{Two examples}\label{s3}

\subsection{}\label{s3s1}
Let a multivalued analytic function $f\in\HH(0)$ be given as a solution of the following
cubic equation
\begin{equation}
w^3+3(z^2+3z+5)w+2(z^3+2z^2+z+1)=0.
\label{cube_2016_1}
\end{equation}
To be more precise, suppose that $f$ is given by the Cardano formula
\begin{align}
f(z):=
&\root3\of{
-(z^3+2z^2+z+1)+\sqrt{(z^2+3z+5)^3+(z^3+2z^2+z+1)^2}}\nonumber\\
-&\frac{z^2+3z+5}
{\root3\of{-(z^3+2z^2+z+1)+\sqrt{(z^2+3z+5)^3+(z^3+2z^2+z+1)^2}}}.
\label{cube_2016_2}
\end{align}
Then $f\in\HH(0)$. The discriminant $D$,
\begin{equation}
D=(z^2+3z+5)^3+(z^3+2z^2+z+1)^2,
\label{cube_2016_3}
\end{equation}
has six simple roots at the points, say, $z_1,z_2,z_3$ and $\myo z_1,\myo
z_2,\myo z_3$, where $\Re z_j<0$, $\Im z_j>0$, $j=1,2,3$. All these six points are branch
points of the function~$f$ given by~\eqref{cube_2016_2} and the points $z=0$ an
$z=\infty$ are not such.

%%% Fig 8-13
\begin{figure}
\begin{center}
\includegraphics[width=0.75\textwidth]{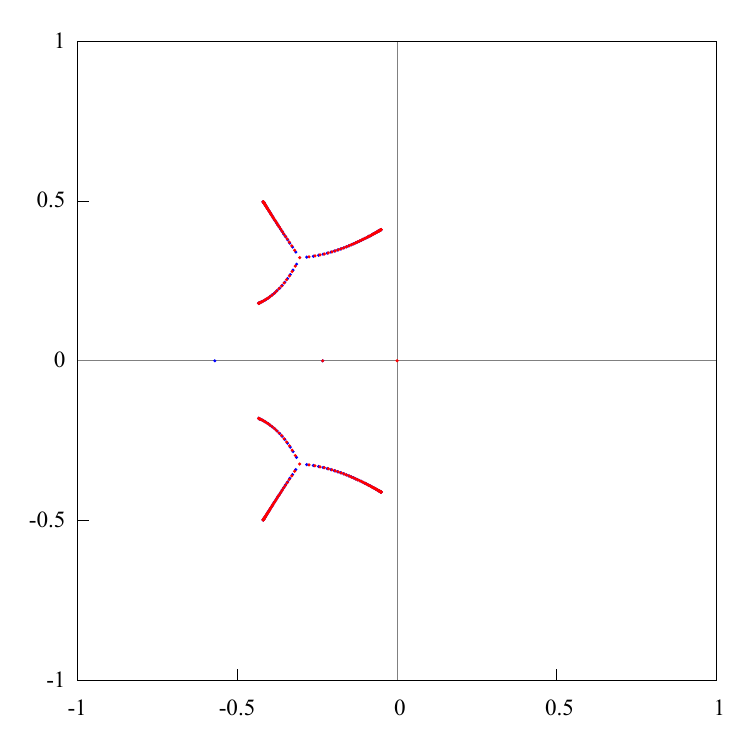}
\vskip-5mm
\caption{The zeros (blue points) and the poles (red points) of PA $[200/200]_f$ for $f$ given by~\eqref{cube_2016_2} are plotted after the transformation $z\mapsto\zeta=1/z$. The corresponding Stahl's compact set $S$ consists of two continuums each of which has a Chebotar\"ev's point of zero density.}
\label{fig_cube_2016_pade_200}
% \end{center}
% \end{figure}

% \begin{figure}
% \begin{center}
\includegraphics[width=0.75\textwidth]{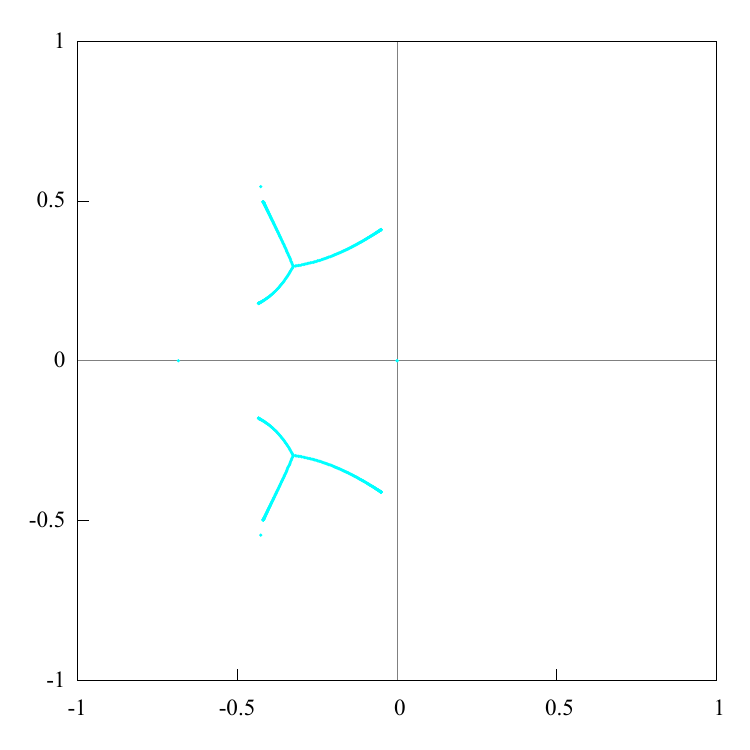}
\vskip-5mm
\caption{The zeros (light blue points) of the type II HP-polynomials $P_{600,0}$ for $f$ given by~\eqref{cube_2016_2} are plotted. The corresponding compact set $E$ consists of two continuums each of which has a Chebotar\"ev's point of positive density.}
\label{fig_cube_2016_hepas_300}
\end{center}
\end{figure}

\begin{figure}
\begin{center}
\includegraphics[width=0.75\textwidth]{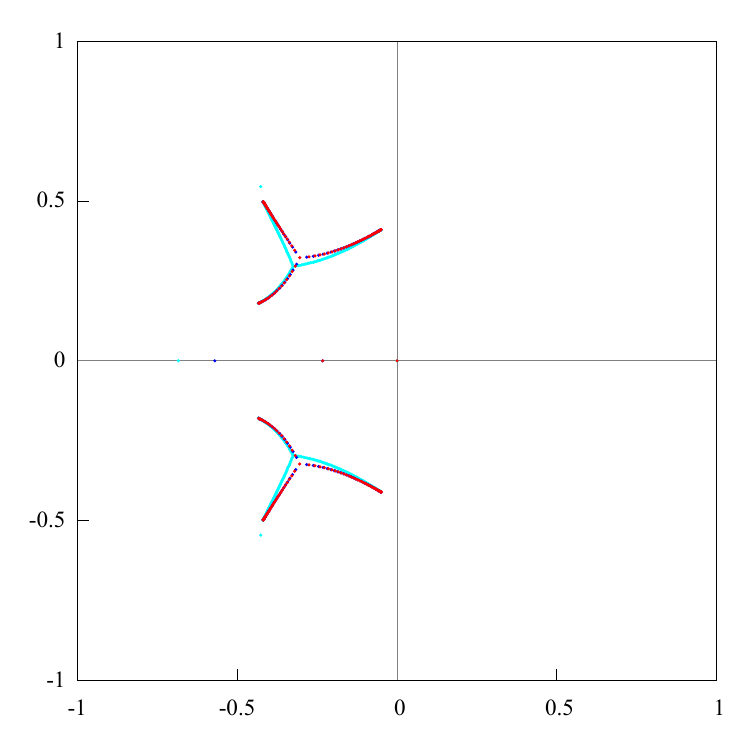}
\vskip-5mm
\caption{We combine the Fig.~\ref{fig_cube_2016_pade_200} and the Fig.~\ref{fig_cube_2016_hepas_300}. Evidently, Stahl's compact set $S$ and the compact set $E$ differ from each other.}
\label{fig_cube_2016_pade_hepas}
% \end{center}
% \end{figure}

% \begin{figure}
% \begin{center}
\includegraphics[width=0.75\textwidth]{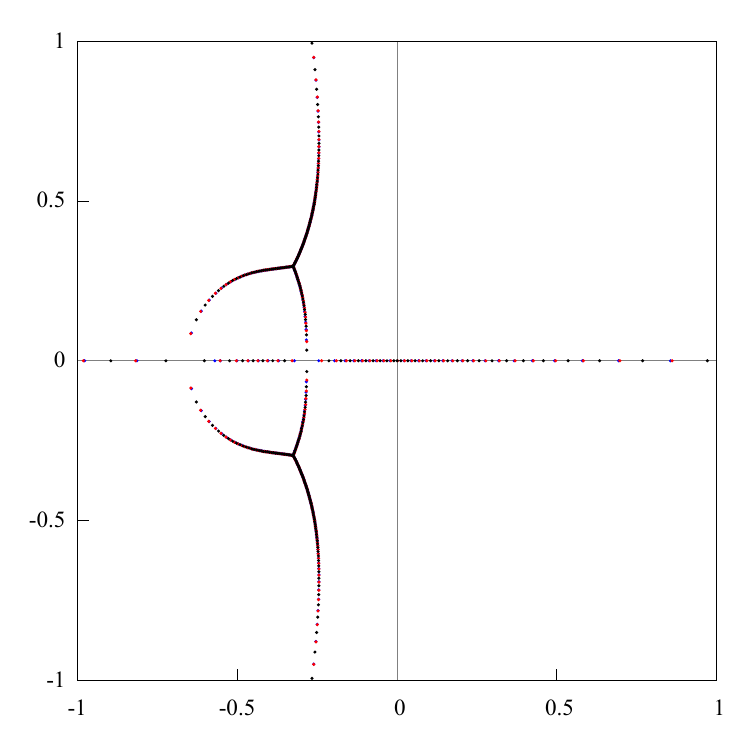}
\vskip-5mm
\caption{The zeros (blue, red and black points) of the type I HP-polynomials $Q_{300,0},Q_{300,1}$ and $Q_{300,2}$ for $f$ given by~\eqref{cube_2016_2} are plotted. The corresponding compact set $F$ is a continuum with two Chebotar\"ev's points of zero density and two Chebotar\"ev's points of positive density. The compact set $F$ shares the complex plane into six domains. In the case under consideration the compact set $F$ does not contain any of the branch points of~$f$.}
\label{fig_cube_2016_hepa_300}
\end{center}
\end{figure}

\begin{figure}
\begin{center}
\includegraphics[width=0.75\textwidth]{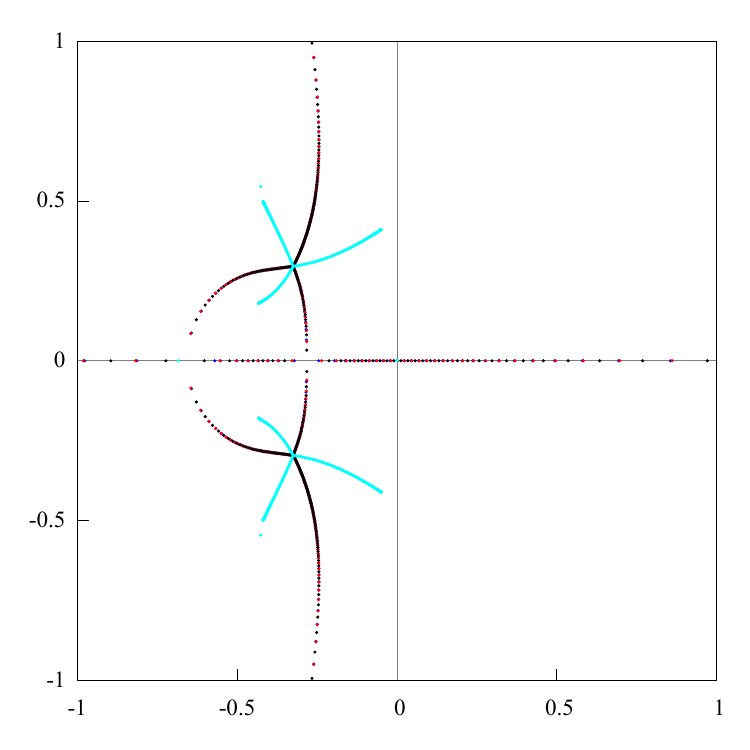}
\vskip-5mm
\caption{We combine the Fig.~\ref{fig_cube_2016_hepas_300} and the Fig.~\ref{fig_cube_2016_hepa_300}. Evidently, the Chebotar\"ev's points of positive density for the compact sets $E$ and $F$ coincide with each other. The ordered pair $(E,F)$ forms the corresponding Nuttall's condenser for~$f$.}
\label{fig_cube_2016_hepas_hepa}
% \end{center}
% \end{figure}

% \begin{figure}
% \begin{center}
\includegraphics[width=0.75\textwidth]{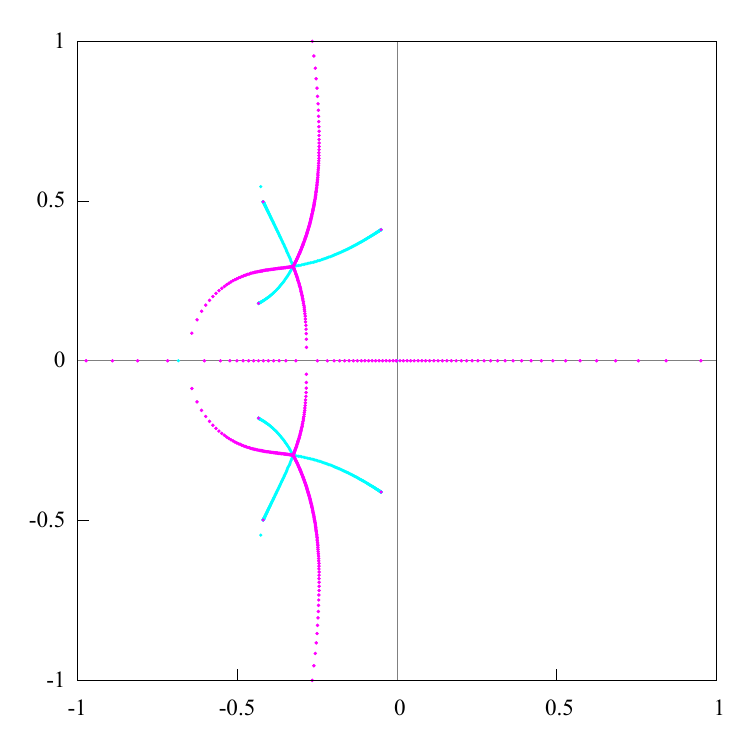}
\vskip-5mm
\caption{The zeros (light blue points) of the type II HP-polynomials $P_{600,0}$ as well as the zeros (violet points) of the discriminant $D_{300}$ of degree $600$ for $f$ given by~\eqref{cube_2016_2} are plotted. The zeros of $D_{300}$ simulate the compact set $F$ similar to the zeros of the type I HP-polynomials $Q_{300,0},Q_{300,1}$ and $Q_{300,2}$; see Fig.~\ref{fig_cube_2016_hepa_300}. In addition, exactly six zeros of $D_{300}$ mark six endpoints of the compact set $E$ all of which are the branch points of the square-root type.}
\label{fig_cube_2016_hepas_hepad}
\end{center}
\end{figure}

In Fig.~\ref{fig_cube_2016_pade_200} the zeros (blue points) and the poles (red points) of PA $[200/200]$ for $f$ are plotted after the transformation $z\mapsto\zeta=1/z$. The numerical distribution of these zeros is in a great accordance with Stahl's theorem, see~\cite{Sta97b}.

In Fig.~\ref{fig_cube_2016_hepas_300} the zeros (light blue points) of the type II HP-polynomial $P_{600,0}$ are plotted after the transformation $z\mapsto\zeta=1/z$. These zeros simulate the $E$ plate of Nuttall's condenser with Chebotar\"ev's point of positive density.

In Fig.~\ref{fig_cube_2016_pade_hepas} we combine Fig.~\ref{fig_cube_2016_pade_200} and Fig.~\ref{fig_cube_2016_hepas_300}. The difference between these zero distributions is caused by the fact that in general Stahl's compact set $S$ and the compact set $E$ give the solutions of two different extremal problems.

In Fig.~\ref{fig_cube_2016_hepa_300}
the zeros (blue, red and black points) of the type I HP-polynomials $Q_{n,0},Q_{n,1},Q_{n,2}$ of degree $n=300$ are plotted after the transformation $z\mapsto\zeta=1/z$. These zeros simulate the $F$ plate of Nuttall's condenser with two Chebotar\"ev's points of positive density and two Chebotar\"ev's points of zero density. It should be also noted that the compact set $F$ shares the complex plane in four domains.

In Fig.~\ref{fig_cube_2016_hepas_hepa} the zeros of the type I HP-polynomials $Q_{300,0},Q_{300,1},Q_{300,2}$ are plotted as well as the zeros of the type II HP-polynomial $P_{600,0}$ after the transformation $z\mapsto\zeta=1/z$. These zeros simulate Nuttall's condenser $(E,F)$.

In Fig.~\ref{fig_cube_2016_hepas_hepad} the zeros (light blue) of the type II HP-polynomial $P_{600,0}$ as well as the zeros (violet points) of the discriminant $D_n:=Q_{n,1}^2-4Q_{n,0}Q_{n,2}$ for $n=300$ are plotted. Clearly each of the six vertexes of the compact set $E$ is marked exactly by one of the zeros of $D_n$. According to the ideas of~\cite{Sue22b} that means that each of the vertexes is a square-root branch point of $f$. This is in a full accordance with the representation~\eqref{cube_2016_2} of $f$.

\subsection{}\label{s3s2}
Let a multivalued analytic function $f$ be given as a solution of the following cubic equation
\begin{equation}
w^3-3(z-1)^2w+2(z-3)^3=0.
\label{cube2_2016_1}
\end{equation}
To be more precise, let $f$ be given by the Cardano formula
\begin{align}
f(z):=\root3\of{-(z-3)^3+\sqrt{(z-3)^6-(z-1)^6}}\nonumber\\
+\frac{(z-1)^2}{\root3\of{-(z-3)^3+\sqrt{(z-3)^6-(z-1)^6}}}.
\label{cube2_2016_2}
\end{align}
Then $f\in\HH(0)$. The discriminant $D$,
\begin{equation}
D=(z-3)^6-(z-1)^6,
\label{cube2_2016_3}
\end{equation}
has five simple roots at the points $z_1=2$, $z_2,z_3$ and $\myo z_2,\myo
z_3$, where $\Re{z}_j,\Im z_j>0$, $j=2,3$. All these five points and the infinity
point $z=\infty$ as well are the branch
points of the function~$f$ given by~\eqref{cube2_2016_2}.

%%% Fig 14-19
\begin{figure}
\begin{center}
\includegraphics[width=0.75\textwidth]{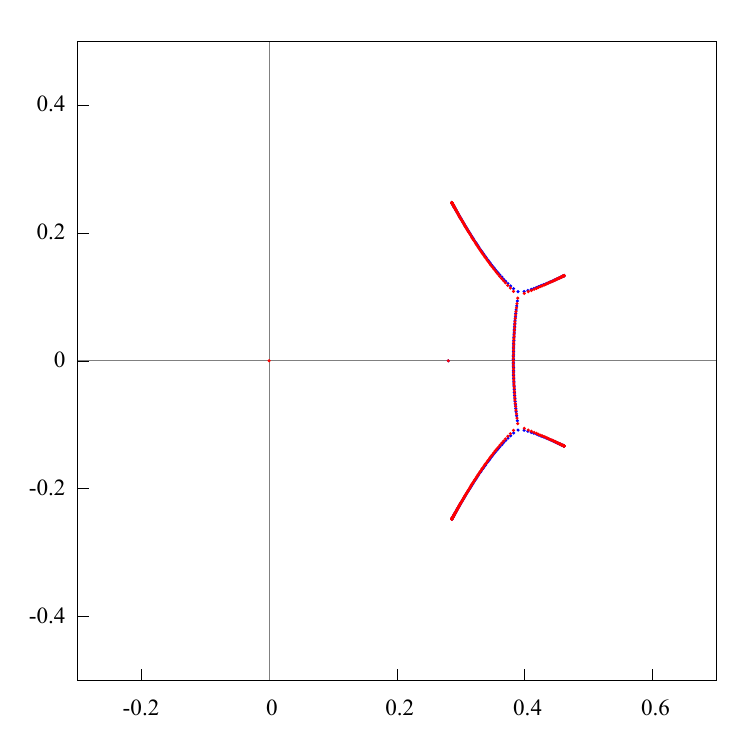}
\vskip-5mm
\caption{The zeros (blue points) and the
poles (red points) of PA $[300/300]_f$ for $f$ given by~\eqref{cube2_2016_2} are plotted after the transformation $z\mapsto\zeta=1/z$. The corresponding Stahl's compact set $S$ is a continuum with two Chebotar\"ev's points of zero density.}
\label{fig_cube2_2016_pade_300}
% \end{center}
% \end{figure}

% \begin{figure}
% \begin{center}
\includegraphics[width=0.75\textwidth]{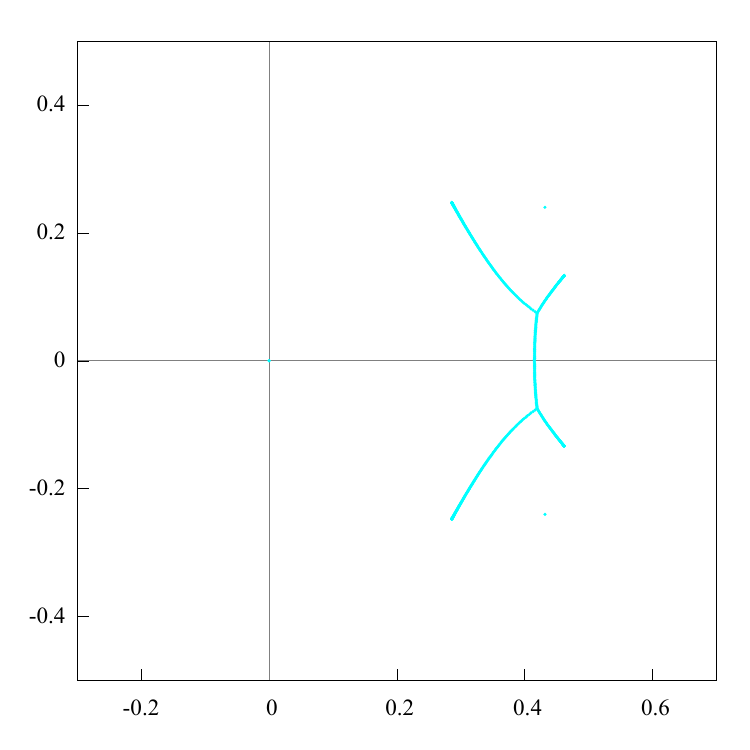}
\vskip-5mm
\caption{The zeros (light blue points) of the type II HP-polynomial $P_{600,0}$ for $f$ given by~\eqref{cube2_2016_2} are plotted. The corresponding compact set $E$ is a continuum with two Chebotar\"ev's points of positive density.}
\label{fig_cube2_2016_hepas_300}
\end{center}
\end{figure}

\begin{figure}
\begin{center}
\includegraphics[width=0.75\textwidth]{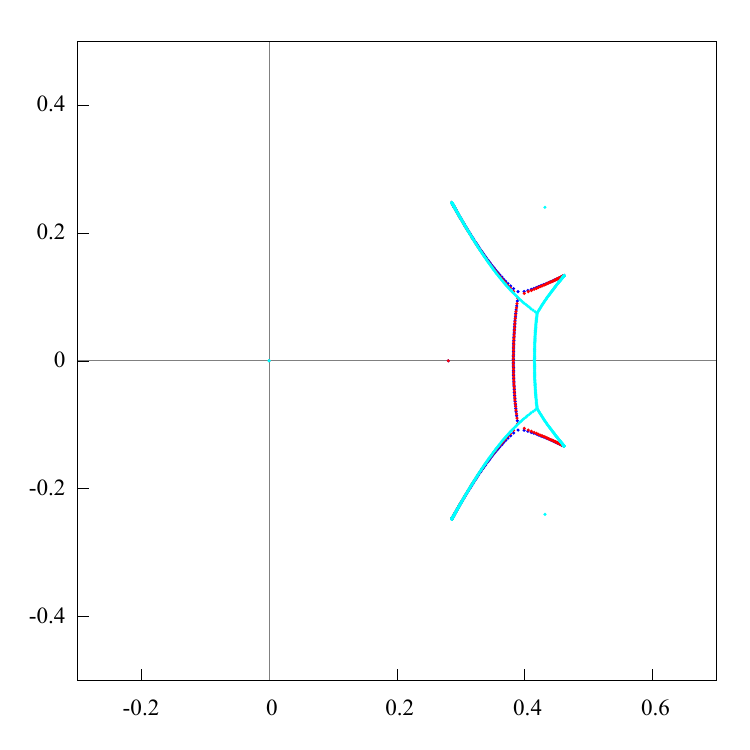}
\vskip-5mm
\caption{We combine the Fig.~\ref{fig_cube2_2016_pade_300} and the Fig.~\ref{fig_cube2_2016_hepas_300}. Evidently, Stahl's compact set $S$ and the compact set $E$ differ from each other.}
\label{fig_cube2_2016_pade_hepas}
% \end{center}
% \end{figure}

% \begin{figure}
% \begin{center}
\includegraphics[width=0.75\textwidth]{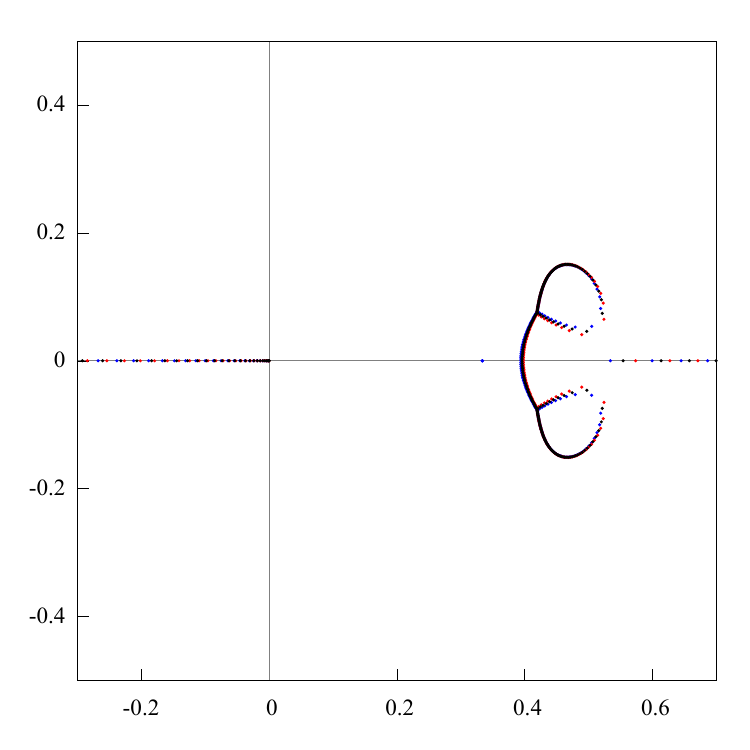}
\vskip-5mm
\caption{The zeros (blue, red and black points) of the type I HP-polynomials $Q_{300,0},Q_{300,1}$ and $Q_{300,2}$ for $f$ given by~\eqref{cube2_2016_2} are plotted. The corresponding compact set $F$ consists of two continuums with two Chebotar\"ev's points of positive density. The compact set $F$ shares the complex plane into three domains. In the case under consideration the compact set $F$ does not contain any of the branch points of~$f$.}
\label{fig_cube2_2016_hepa_300}
\end{center}
\end{figure}

\begin{figure}
\begin{center}
\includegraphics[width=0.75\textwidth]{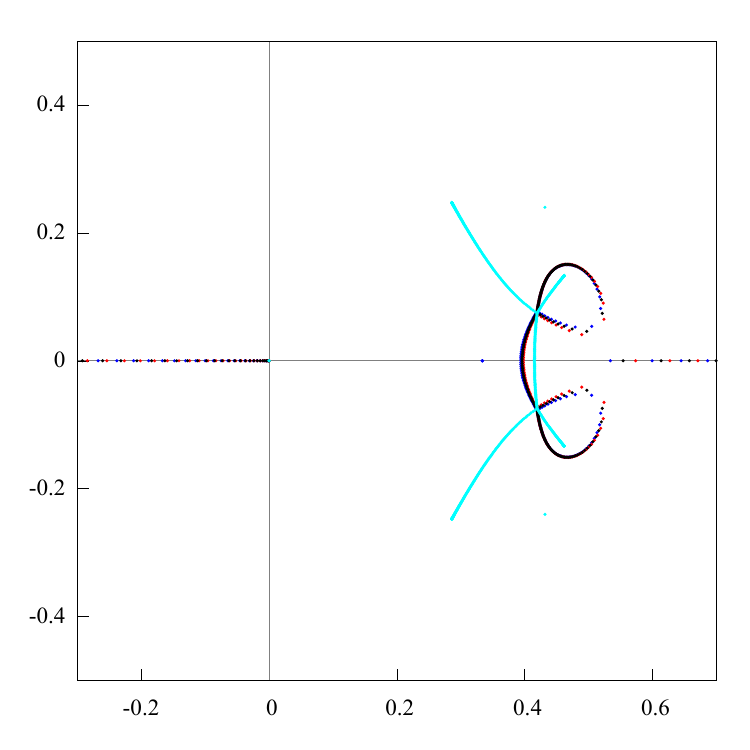}
\vskip-5mm
\caption{We combine the Fig.~\ref{fig_cube2_2016_hepas_300} and the Fig.~\ref{fig_cube2_2016_hepa_300}. Evidently, the Chebotar\"ev's points of positive density for the compact sets $E$ and $F$ coincide with each other. The ordered pair $(E,F)$ forms the corresponding Nuttall's condenser for~$f$.}
\label{fig_cube2_2016_hepas_hepa_300}
% \end{center}
% \end{figure}

% \begin{figure}
% \begin{center}
\includegraphics[width=0.75\textwidth]{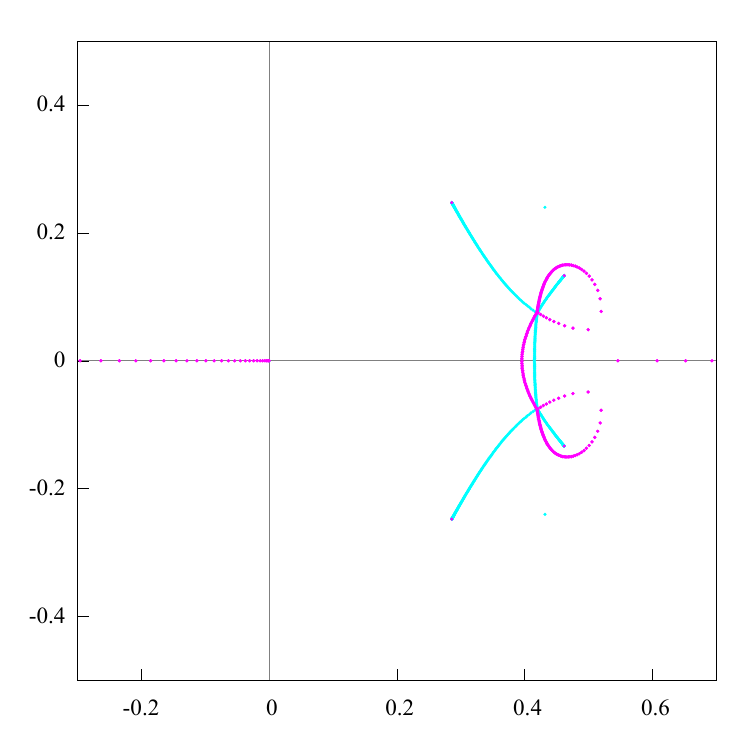}
\vskip-5mm
\caption{The zeros (light blue points) of the type II HP-polynomial $P_{600,0}$ as well as the zeros (violet points) of the discriminant $D_{300}$ for $f$ given by~\eqref{cube2_2016_2} are plotted. The zeros of $D_{300}$ simulate the compact set $F$ similar to the zeros of the type I HP-polynomials $Q_{300,0},Q_{300,1}$ and $Q_{300,2}$. In addition, exactly four zeros of $D_{300}$ mark four endpoints of the compact set $E$ all of which are branch points of square-root type.}
\label{fig_cube2_2016_hepas_hepad}
\end{center}
\end{figure}

%%% Fig 20-28
\begin{figure}
\begin{center}
\includegraphics[width=0.75\textwidth]{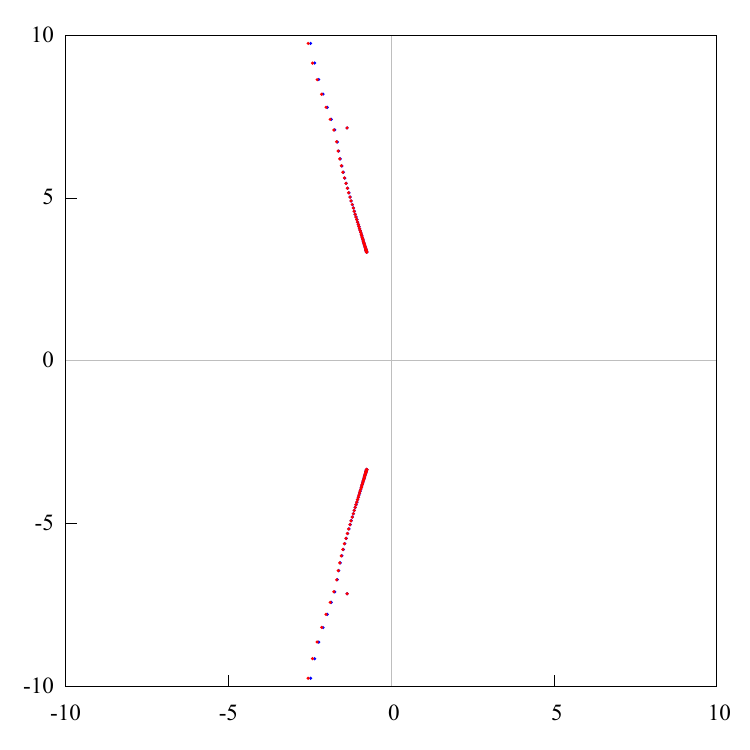}
\vskip-5mm
\caption{The zeros (blue points) and the poles (red points) of PA $[143/143]$ for the frequency function $\nu$ are plotted. In accordance to Stahl's theorem they accumulate to the ``nearest'' (with respect to the origin) pair of branch points of~$\nu$. It is impossible to recognize the type of these ``blue-red'' singularities based on the PA-analysis only.}
\label{fig_pade_143(10)}
% \end{center}
% \end{figure}

% \begin{figure}
% \begin{center}
\includegraphics[width=0.75\textwidth]{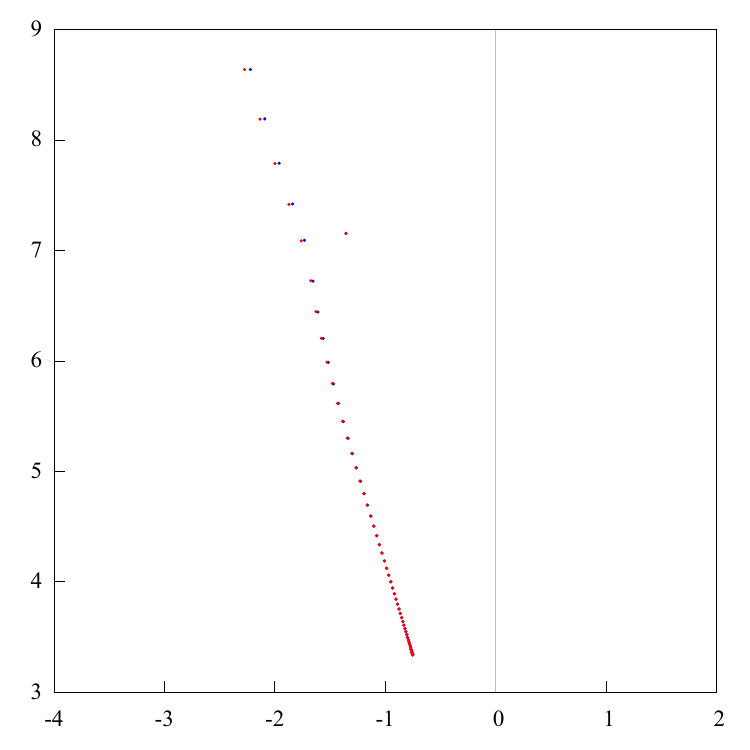}
\vskip-5mm
\caption{A part of the Fig.~\ref{fig_pade_143(10)} is presented in a different scale.}
\label{fig_pade_143(39)}
\end{center}
\end{figure}

\begin{figure}
\begin{center}
\includegraphics[width=0.75\textwidth]{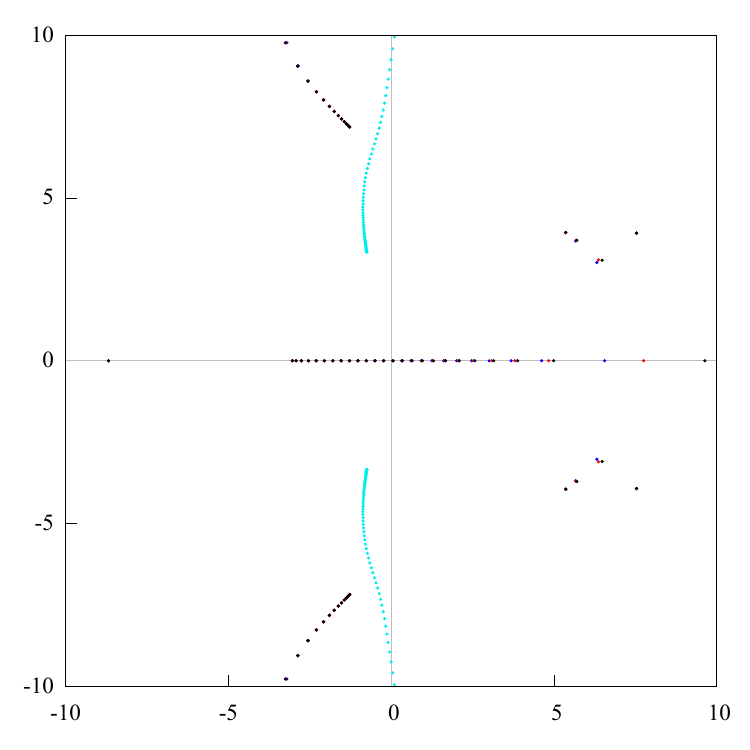}
\vskip-5mm
\caption{The zeros (light blue points) of the type II HP-polynomial $P_{190,0}$ and the zeros (blue, red and black points) of the type I HP-polynomial $Q_{95,0},Q_{95,1}$ and $Q_{95,2}$ are plotted. These zeros simulate the $E$ plate and the $F$ plate of Nuttall's condenser respectively. The endpoints of the $E$ plate mark just the same singularities of the function $\nu$ as the endpoints of Stahl's compact set. Since these ``light-blue'' singularities don't belong to the $F$ plate, they are exactly square-root type singularities.}
\label{fig_hepa_hepas_95(10)}
\end{center}
\end{figure}

\begin{figure}
\begin{center}
\includegraphics[width=0.75\textwidth]{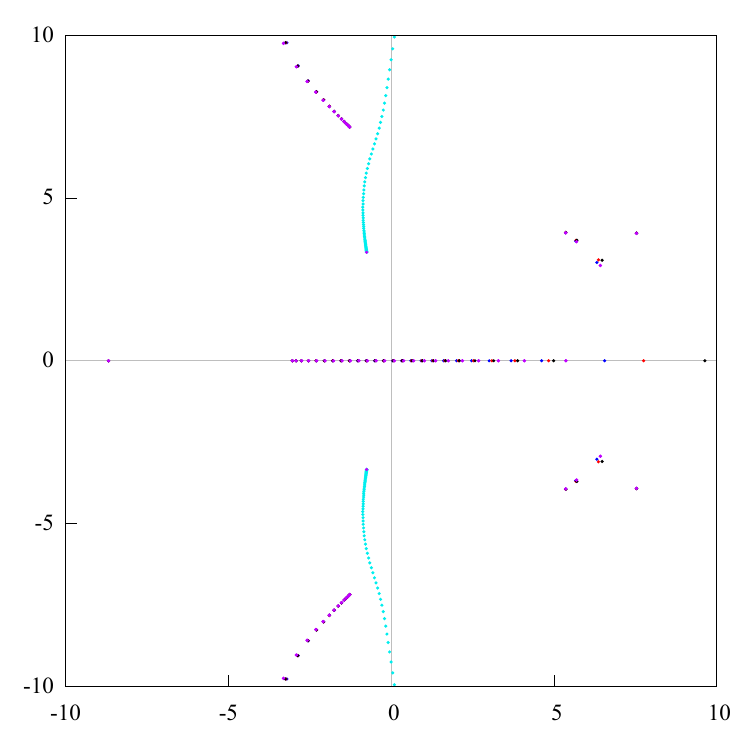}
\vskip-5mm
\caption{The zeros (light blue points) of the type II HP-polynomial $P_{190,0}$, the zeros (blue, red and black points) of the type I HP-polynomial $Q_{95,0},Q_{95,1}$ and $Q_{95,2}$ and the zeros (violet points) of the discriminant polynomial $D_{95}$ are plotted. Exactly two zeros of $D_{95}$ mark the endpoints of the $E$ plate. This fact gives an additional evidence to the square-root type singularities of $\nu$.}
\label{fig_hepa_hepas_hepad_95(10)}
% \end{center}
% \end{figure}

% \begin{figure}
% \begin{center}
\includegraphics[width=0.75\textwidth]{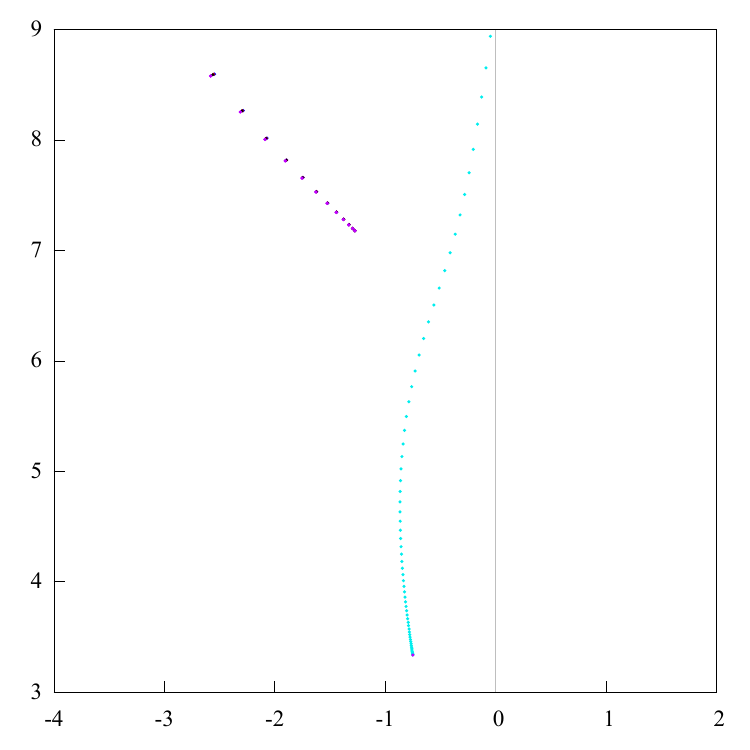}
\vskip-5mm
\caption{A part of the Fig.~\ref{fig_hepa_hepas_hepad_95(10)} is presented in a different scale.}
\label{fig_hepa_hepas_hepad_95(39)}
\end{center}
\end{figure}

\begin{figure}
\begin{center}
\includegraphics[width=0.75\textwidth]{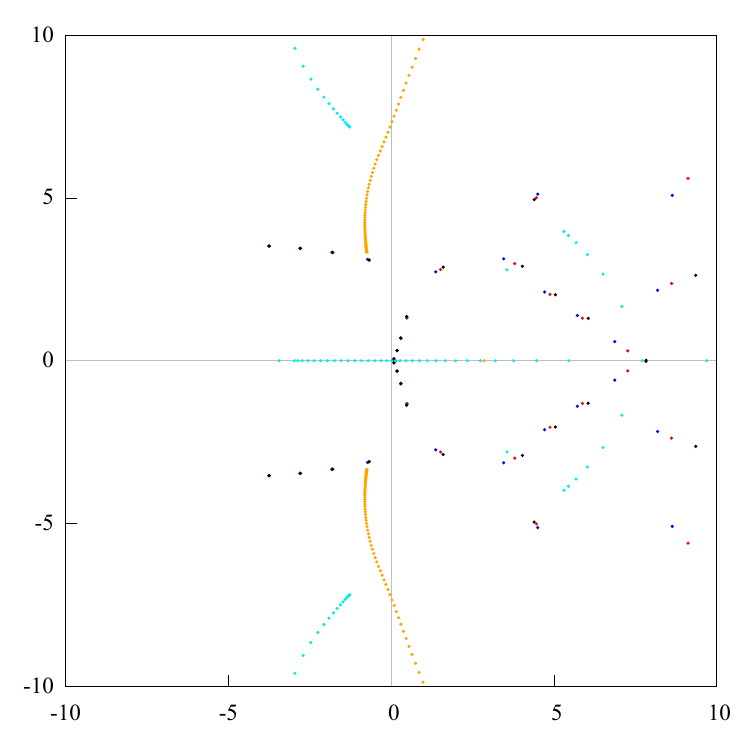}
\vskip-5mm
\caption{The zeros (blue, red and black points) of the type I HP-polynomials $Q_{71,j}$ for the tuple $[1,\nu,\nu^2,\nu^3]$ are plotted. Also the zeros (yellow points) of the type II HP-polynomial $P_{213,0}$ as well as the zeros (light blue points) of the generalized HP-polynomial of degree $141$ are plotted. Due to some general results, the distribution of zeros of these three types of HP-polynomials corresponds to Nuttall's partition of the four-sheeted Riemann surface associated with~$\nu$.}
\label{fig_mhepas_mhep_71(10)}
% \end{center}
% \end{figure}

% \begin{figure}
% \begin{center}
\includegraphics[width=0.75\textwidth]{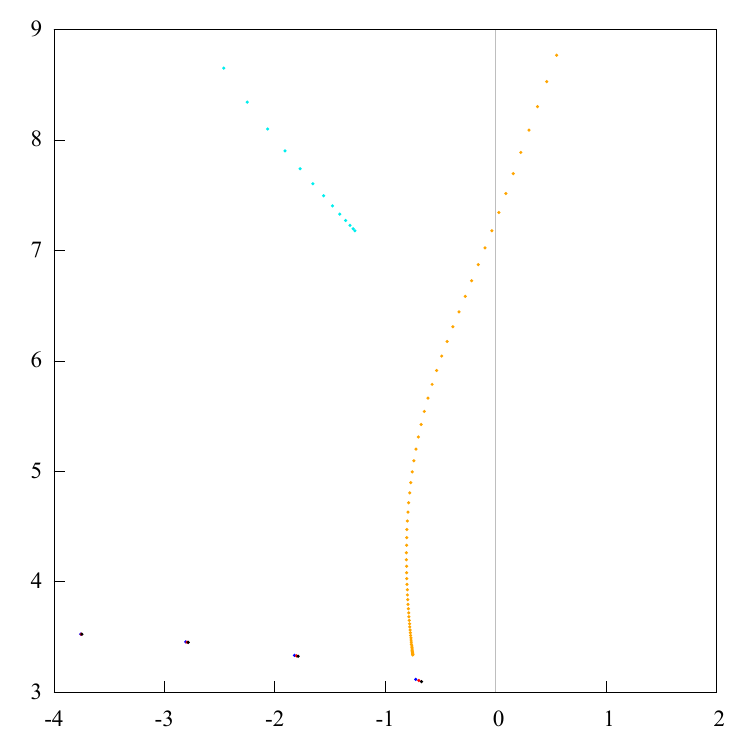}
\vskip-5mm
\caption{A part of the Fig.~\ref{fig_mhepas_mhep_71(10)} is presented in a different scale.}
\label{fig_mhepas_mhep_71(39)}
\end{center}
\end{figure}

\begin{figure}
\begin{center}
\includegraphics[width=0.75\textwidth]{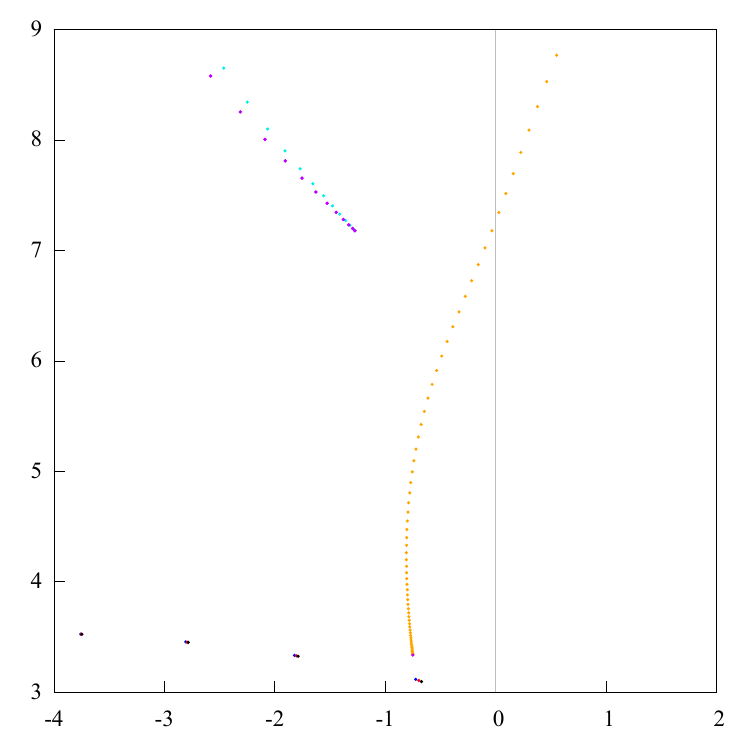}
\vskip-5mm
\caption{We combine the zeros of the three types of the HP-polynomials with the zeros of the discriminant polynomial $D_{95}$.}
\label{fig_hepa_hepad_95_mhepas_71(39)}
% \end{center}
% \end{figure}

% \begin{figure}
% \begin{center}
\includegraphics[width=0.75\textwidth]{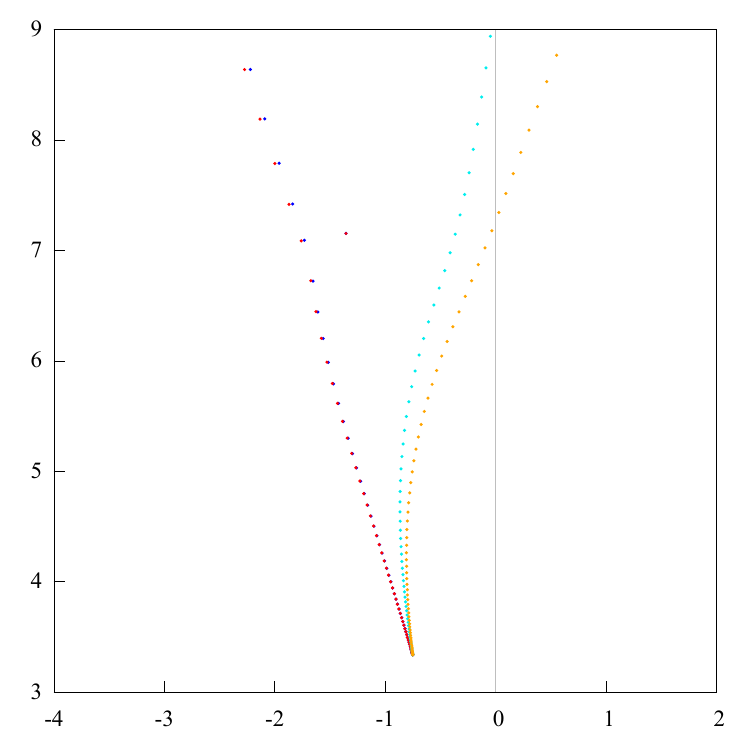}
\vskip-5mm
\caption{We combine the zeros (blue points) and the poles (red points) of PA $[143/143]$ for the frequency function $\nu$, the zeros (light blue points) of the type II HP-polynomial $P_{180,0}$ for the system $\nu,\nu^2$ and the zeros (yellow points) of the type II HP-polynomial $P_{213,0}$ for the system $\nu,\nu^2,\nu^3$.}
\label{fig_pade_hepas_mhepas(39)}
\end{center}
\end{figure}

In Fig.~\ref{fig_cube2_2016_pade_300} the zeros (blue points) and the poles (red points) of the PA $[300/300]$ for $f$ are plotted after the transformation $z\mapsto\zeta=1/z$. The numerical distribution of these zeros is in a great accordance with Stahl's theorem, see~\cite{Sta97b}.

In Fig.~\ref{fig_cube2_2016_hepas_300} the zeros (light blue points) of the type II HP-polynomial $P_{600,0}$ are plotted after the transformation $z\mapsto\zeta=1/z$. These zeros simulate the $E$ plate of Nuttall's condenser with Chebotar\"ev's point of positive density.

In Fig.~\ref{fig_cube2_2016_pade_hepas} we combine Fig.~\ref{fig_cube2_2016_pade_300} and Fig.~\ref{fig_cube2_2016_hepas_300}. The difference between these zero distributions is caused by the fact that in general Stahl's compact set $S$ and the compact set $E$ give the solutions of two different extremal problems.

In Fig.~\ref{fig_cube2_2016_hepa_300}
the zeros (blue, red and black points) of the type I HP-polynomials $Q_{n,0},Q_{n,1},Q_{n,2}$ of degree $n=300$ are plotted after the transformation $z\mapsto\zeta=1/z$. These zeros simulate the $F$ plate of Nuttall's condenser with two Chebotar\"ev's points of positive density and two Chebotar\"ev's points of zero density. It should also be noted that the compact set $F$ shares the complex plane in three domains.

In Fig.~\ref{fig_cube2_2016_hepas_hepa_300} the zeros of
the type I HP-polynomials $Q_{300,0},Q_{300,1},Q_{300,2}$ are plotted as well as the zeros of the type II HP-polynomial $P_{600,0}$ after the transformation $z\mapsto\zeta=1/z$. These zeros simulate Nuttall's condenser $(E,F)$.

In Fig.~\ref{fig_cube2_2016_hepas_hepad} the zeros (light blue) of the type II HP-polynomial $P_{600,0}$ as well as the zeros (violet points) of the discriminant $D_n:=Q_{n,1}^2-4Q_{n,0}Q_{n,2}$ for $n=300$ are plotted. Clearly each of the four vertexes of the compact set $E$ is marked exactly by one of the zeros of $D_n$. According to the ideas of~\cite{Sue22b} that means that each of the vertexes is a square-root type singular point of $f$. This is in a full accordance with the representation~\eqref{cube2_2016_2} of~$f$.

\section{Free Van der Pol equation}\label{s4}

\subsection{}\label{s4s1}

The well-known free Van der Pol equation is a mathematical model of a physical
device which generates undamped periodic oscillations of an electric current. After a certain normalization, the equation takes the form
\begin{equation}
\frac{d^2U}{dt^2}+\eps(U^2-1)\frac{dU}{dt}+U=0,
\label{vaneqini}
\end{equation}
where $U=U(t;\eps)$ is a quantity related to the current strength in the electric circuit, $t$ is time, and the
physical characteristics of the device itself are described by a small parameter $\eps$.
It is also well known that the Van der Pol equation~\eqref{vaneqini} for any $\eps>0$ has a unique limit cycle in the
phase plane $(U,dU/dt)$ and this limit cycle is stable.

The properties of the limit cycle of the free Van der Pol equation~\eqref{vaneqini} were studied numerically in several papers by using
Pad\'e analysis or Hermite--Pad\'e analysis; see~\cite{AnGe82}, \cite{AnGeDa84}, \cite{Sue10}, \cite{AmBoFe18}.

Let us describe in short the standard procedure for constructing a formal solution to the Van der Pol
equation~\eqref{vaneqini} as a series in powers of the small parameter $\eps$ whose coefficients are trigonometric
polynomials depending on the new variable $x=\nu(\eps)t$.
We make the change of variables $x=\nu(\eps)t$ in~\eqref{vaneqini}, where $\nu=\nu(\eps)=2\pi/T(\eps)$ is the frequency of the limit cycle and $T=T(\eps)$ is its period. Set now
$u(x;\eps)=U(t;\eps)$. The new function $u$ is periodic in $x$ variable with period $2\pi$ and thus the equation~\eqref{vaneqini} takes the form
\begin{equation}
\nu^2\ddot{u}+\eps\nu(u^2-1)\dot{u}+u=0,
\label{vaneq}
\end{equation}
where $\dot{u}=du/dt$, $\ddot{u}=d^2u/dt^2$, $u=u(x;\eps)$, $\dot{u}(0;\eps)\equiv0$, $u(0;\eps)=A(\eps)>0$. 

Let us represent the functions $u(x;\eps)$ and $\nu(\eps)$ as formal power series in~$\eps$:
\begin{equation}
u(x;\eps)=\sum_{k=0}^\infty u_k(x)\eps^k,
\quad \nu(\eps)=\sum_{k=0}^\infty\myt{\nu}_k\eps^k.
\label{vanexp}
\end{equation}
Then it turns out that the function $\nu$ depends on $\eps^2$ only and thus we have that
\begin{equation}
\nu(\eps)=1+\sum_{j=1}^\infty\nu_j\eps^{2j}
\label{nuexp}
\end{equation}
and the coefficients $\nu_j$ can be computed step by step based on initial data and via some recurrence relations; see~\cite{AnGe82}, \cite{AnGeDa84}, \cite{Sue10}, \cite{AmBoFe18}.

The first several terms of the series~\eqref{nuexp} are
\begin{equation*}
\nu(\eps)=1-\frac1{16}\eps^2+\frac{17}{3072}\eps^4+
\frac{35}{884736}\eps^6-\frac{678899}{5096079360}
\eps^8+\frac{28160413}{2293235712000}\eps^{10}+ O(\eps^{12}).
\end{equation*}
For more details see~\cite{AnGe82}, \cite{AnGeDa84}, \cite{Sue10}, \cite{AmBoFe18}, \cite{OrCh23}.

The main purpose of this paper is to demonstrate the advantage of HP-analysis over PA-analysis. In particular, for this purpose we computed $N=287$ ($N=2n+1=3m+2=4\ell+3$, $n=143$, $m=95$, $\ell=71$) Taylor coefficients $\nu_0,\dots,\nu_{286}$ of the expansion~\eqref{nuexp} of the frequency function with the precision of $400$ decimal digits.
In~\cite{AmBoFe18} it was mentioned
that the application of the generating function method developed by Fern\'andez et al.~\cite{FeArCa87} shows that the pair of two conjugate and nearest to the origin singularities of $\nu$ are of square-root type. We confirm here this result basing on the analysis of numerical distribution of the zeros of the type I HP-polynomials $Q_{95,j}$ for the tuple $[1,\nu,\nu^2]$, numerical distribution of the zeros of the type II HP-polynomials $P_{180,j}$ for the pair of functions $\nu,\nu^2$ as well as numerical distribution of the zeros of the discriminant $D_{95}$ given~\eqref{defdis}. Moreover we treat the numerical distribution of zeros of type I HP-polynomials for the tuple $[1,\nu,\nu^2,\nu^3]$, numerical distribution of zeros of type II HP-polynomials for the system of the functions $\nu,\nu^2,\nu^3$ as well as numerical distribution of zeros of a generalization of the classical construction of Hermite--Pad\'e polynomials introduced in~\cite{Kom21}.

For the general definition of type I HP-polynomials for a tuple $[1,f,f^2,\dots,f^m]$ as well as the definition of type II HP-polynomials for a system $f,f^2,\dots,f^m$ see~\cite{Nut84}, \cite{NiSo88}.

\subsection{}\label{s4s2}

In Fig.~\ref{fig_pade_143(10)}--Fig.~\ref{fig_pade_hepas_mhepas(39)} the results of computation of the zeros of Pad\'e and Hermite--Pad\'e polynomials are presented. The results are based on the computed $N=287$ Taylor coefficients $\nu_0,\dots,\nu_{286}$ of the power series~\eqref{nuexp}. The main purpose of the analysis of the obtained numerical results is to establish that {\it two pairs} of the ``nearest'' branch type singularities of the frequency function $\nu(\eps)$ are of square-root type. Recall that in~\cite{AmBoFe18} a similar conclusion was made for the first pair only and was based on a different analysis.
The justification of our analysis is based on the theoretical results which state that the limit zero distribution of the HP-polynomials simulates Nuttall's partition of the associated three-sheeted Riemann surface; see~\cite{KoPaSuCh17}, \cite{Sue18d}, \cite{Kom21}.

In Fig.~\ref{fig_pade_143(10)} and Fig.~\ref{fig_pade_143(39)} the zeros (blue points) and the poles (red points) of PA $[90/90]$ for the function $\nu$ are plotted in a different scale. In accordance with Stahl's theorem they accumulate to the ``nearest'' (with respect to the origin) branch type singularities of~$\nu$. It is impossible to recognize the type of these ``blue-red'' singularities based on the PA-analysis only.

In Fig.~\ref{fig_hepa_hepas_95(10)} the zeros (light blue points) of type II HP-polynomial $P_{180,0}$ and the zeros (blue, red and black points) of the type I HP-polynomial $Q_{95,0},Q_{95,1}$ and $Q_{95,2}$ are plotted. These zeros simulate the $E$ plate and $F$ plate of Nuttall's condenser respectively. An end point of the $E$ plate marks just the same singularity of the function $\nu$ as the end point of Stahl's compact set (see Fig.~\ref{fig_pade_143(10)} and Fig.~\ref{fig_pade_143(39)}). Thus the ``light-blue'' singularity coincides with the ``blue-red'' singularity (see Fig.~\ref{fig_pade_hepas_mhepas(39)}). Since this ``light-blue'' singularity doesn't belong to the $F$ plate, it is exactly a square-root type singularity. Indeed, since the condenser $(E,F)$ corresponds to Nuttall's partition of the associated three-sheeted Riemann surface of $\nu$, then it follows from the local structure of $(E,F)$ that when we go twice around this ``light-blue'' singular point, we return to the same value of the multivalued analytic function~$\nu(\eps)$.

As an additional evidence of our conclusion on the square-root type singularity of the point under consideration, let us analyze the numerical distribution of the zeros of the discriminant polynomial $D_{95}$ of degree $180$ (see~\eqref{defdis}). According to the results and conjecture of~\cite{Sue22b}, the zeros of $D_m$ mark the square-type singularities of the initial function as {\it isolated singularities}. This is similar to the well known property of PA to mark polar singularities of the initial function as isolated singularities.
In Fig.~\ref{fig_hepa_hepas_hepad_95(10)} and Fig.~\ref{fig_hepa_hepas_hepad_95(39)} the zeros (light blue points) of the type I HP-polynomials $P_{180,0}$, the zeros (blue, red and black points) of the type II HP-polynomials $Q_{95,0},Q_{95,1},Q_{95,2}$ and the zeros (violet points) of the discriminant polynomial $D_{95}$ (of degree $180$) are plotted. Evidently a single zero of $D_{95}$ marks the end point of the $E$ plate.

There is another pair of conjugate singular points of $\nu$ in these figures which coincide with some end points of the $F$ plate of Nuttall's condenser. However we can say nothing on the type of this pair of conjugate singularities. To analyze the situation we should consider the next level of HP-polynomials for the tuple $[1,\nu,\nu^2,\nu^3]$ and for the system $\nu,\nu^2,\nu^3$.

To create Fig.~\ref{fig_mhepas_mhep_71(10)}--Fig.~\ref{fig_mhepas_mhep_71(39)} we considered the next level of HP-polynomials for the tuple $[1,\nu,\nu^2,\nu^3]$ and for the system $\nu,\nu^2,\nu^3$ as wells as Komlov's generalization of HP-polynomials~\cite{Kom21}, \cite{Kom21b}.
In Fig.~\ref{fig_mhepas_mhep_71(10)} the zeros (blue, red and black points) of the type I HP polynomials $Q_{71,j}$ for the tuple $[1,\nu,\nu^2,\nu^3]$ are plotted. Also the zeros (yellow points) of the type II HP-polynomial $P_{213,0}$ as well as the zeros (light blue points) of the generalized HP-polynomial of degree $141$ are plotted. Due to the results of~\cite{Nut84}, \cite{KoPaSuCh17}, \cite{Kom21} the distribution of the zeros of these three types of the HP-polynomials correspond to Nuttall's partition of the four-sheeted Riemann surface associated with the function~$\nu$. Analyzing from this point of view the mutual location of the zeros of these three types of the HP-polynomials, we conclude that the second pair of the conjugated singular points of $\nu$ is also of the square-root type.
In Fig.~\ref{fig_hepa_hepad_95_mhepas_71(39)} we combine the zeros of three types of the HP-polynomials with the zeros of the discriminant polynomial $D_{180}$.

\newpage
\clearpage

\section{Hermite--Pad\'e Polynomials and Katz's Points}\label{s5}

In some actual problems of molecular physics there are computations which are based on the evaluation of the power series of perturbation theory (SPT). In the case of the problem of model operator perturbation, perhaps the most common formulation of perturbation theory is the Rayleigh--Schr\"{o}dinger perturbation theory (see~\cite{Sch26}). In it, the full Hamiltonian is divided into two parts -- the starting approximation $H_0$ with a known solution and the perturbation $V(z)$, which is frequently expressed in a linearly dependent form $V(z) = zV$.
In accordance with the Katz's results~\cite{Katz62} it is supposed that the initial multivalued analytic function of energy $E(z)$ that corresponds to SPT, is a real-valued function on the real axis which has only quadratic branching at the pairs of complex-conjugated points. The main problem is to evaluate the multivalued energy function $E(z)$ at the point $z=1$ via the SPT given at the point $z=0$.
If the function $E(z)$ has singular points on the closed unit disk, $\myo{\mathbb D}:=\{z\in\mathbb C:|z|\leq1\}$, then the desired evaluation is impossible directly via given power series and it is necessary to involve methods of summation of divergent power series.

Thus as the first step, it is necessary to recognize are there any singular points of the given branch of the function $E(z)$ in the disk $\myo{\mathbb D}$ or not. When such points exist, they are called {\it Katz's points}. In fact Katz's results~\cite{Katz62} establish relations between different branches of the multivalued analytic function $E(z)$. Namely according to Katz's results, each perturbed state $E_{\myj}(z)$ (i.e., a branch of the energy function $E(z)$) can be analytically continued to another perturbed state $E_{\myk}(z)$ through their common branch points.
When two different branches $E_{\myj}(z)$ and $E_{\myk}(z)$ have common Katz's points (on the closed unit disk $\myo{\mathbb D}$), these points are called {\it resonance points} corresponding to these ${\myj}$\,th and ${\myk}$\,th states; see~\cite{ByDu16},~\cite{KrDoSyPa20},~\cite{Duc17},~\cite{ChDoKr22}.

Here we consider three cases of states of the water molecule $\mathrm{H}_2\mathrm{O}$: $(100)$--$(020)$, $(300)$--$(102)$ and $(201)$--$(121)$--$(041)$ (see~\cite{Duc17} and \cite{ChDoKr22}). 

Based on these three examples we describe here in short an algorithm of how to find the common Katz's points for two states. For given $N=2n+1=3m+2$ Taylor coefficients of a power series at the point $z=0$, the algorithm is based on the diagonal PA $[n/n]$, type II HP polynomials $P_{2m,0}$ and the discriminant polynomials $D_m$ of degree $2m$ (see Section~\ref{s1s3} above).

\subsection{\texorpdfstring{Case of the states $\myk=(100)$ and $\myj=(020)$}{Case of the states k=(100) and j=(020)}}\label{s5s1}

The energies $E_{(100)}$ and $E_{(020)}$ of the states $(100)$ and $(020)$ are connected to each other by Fermi resonance points (see~\cite[Sec 3.3.5]{Duc17},
\cite{KrDoSyPa20} and \cite{ChDoKr22}). Since for both cases, $(100)$ and $(020)$, there are $401=2n+1=3m+2$ (with $n=200$ and $m=133$) Taylor coefficients of power expansions of $E_{(100)}$ and $E_{(020)}$ at the point $z=0$, then according to Section~\ref{s1s3} we can apply for our usage the PA $[200/200]$, the type II HP polynomial $P_{266,0}$ and the discriminant polynomial $D_{133}$ of degree $266$.

Recall that according to Stahl's Theory~\cite{Sta97b}, all but a finite 
number of zeros and poles of the diagonal PA $[n/n]_f$ of $f$ are attracted to Stahl's compact set $S=S(f)$ of~$f$. Moreover, each pole of $f$ attracts as many poles of $[n/n]_f$ as its multiplicity. Just the same is valid for the zeros of $f$. In both cases the rate of attraction is exponential (see also~\cite{GoSu04}). Therefore the first step of the algorithm is to apply this property of zeros and poles of PA to find out if there are common Katz's points for the states $(100)$ and $(020)$ or not.

Notice, that actually in Stahl's Theory it is stated that the number of those exceptional zeros and poles of $[n/n]_f$ is $o(n)$ as $n\to\infty$. 
However from the paper~\cite{ApYa15} (see also~\cite{Sue22})
it follows that in the case of analytic functions which is under consideration here, the number of the exceptional zeros and poles of $[n/n]_f$ is finite and depends on the function $f$ only.

\begin{figure}
\begin{center}
\includegraphics[width=0.75\textwidth]{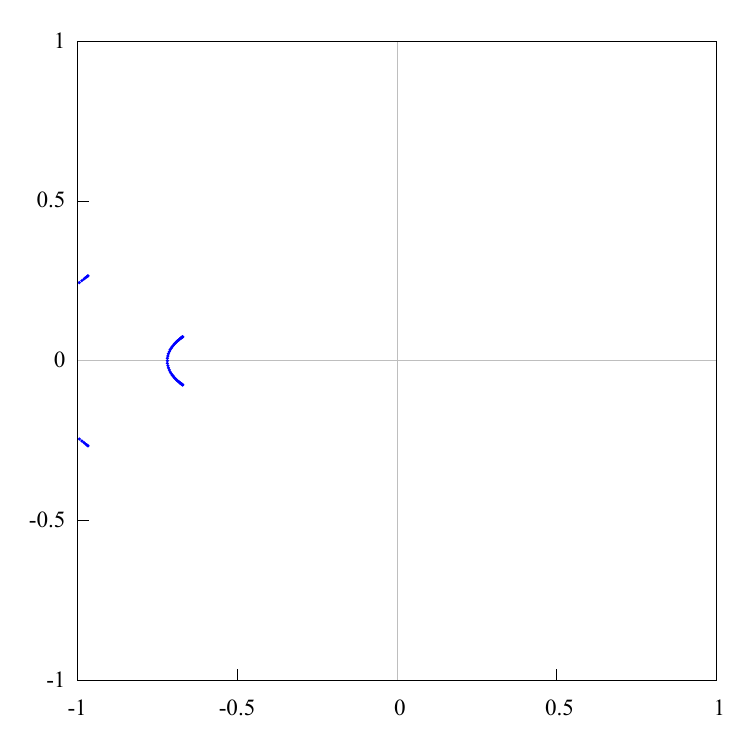}
\vskip-5mm
\caption{The part of the zeros (blue points) of the PA $[200/200]_f$ for the energy function $f=E_{(100)}(z)$ that belongs to the $xy$-square $[-1,1]\times[-1,1]$. According to Stahl's Theory~\cite{Sta97b} the zeros of $[200/200]_f$ simulate Stahl's compact set $S=S(f)$ of $f$.}
\label{fig_PA_100(200)}
% \end{center}
% \end{figure}

% \begin{figure}
% \begin{center}
\includegraphics[width=0.75\textwidth]{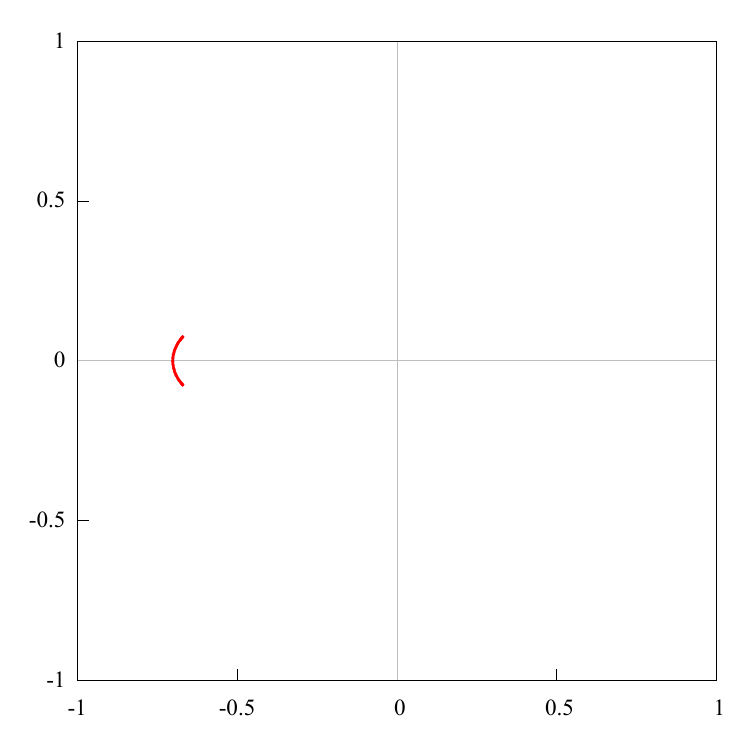}
\vskip-5mm
\caption{The part of the poles (red points) of the PA $[200/200]_f$ for the energy function $f=E_{(020)}(z)$ that belongs to the $xy$-square $[-1,1]\times[-1,1]$. According to Stahl's Theory~\cite{Sta97b} the poles of $[200/200]_f$ simulate Stahl's compact set $S=S(f)$ of $f$.}
\label{fig_PA_020(200)}
\end{center}
\end{figure}

\begin{figure}
\begin{center}
\includegraphics[width=0.75\textwidth]{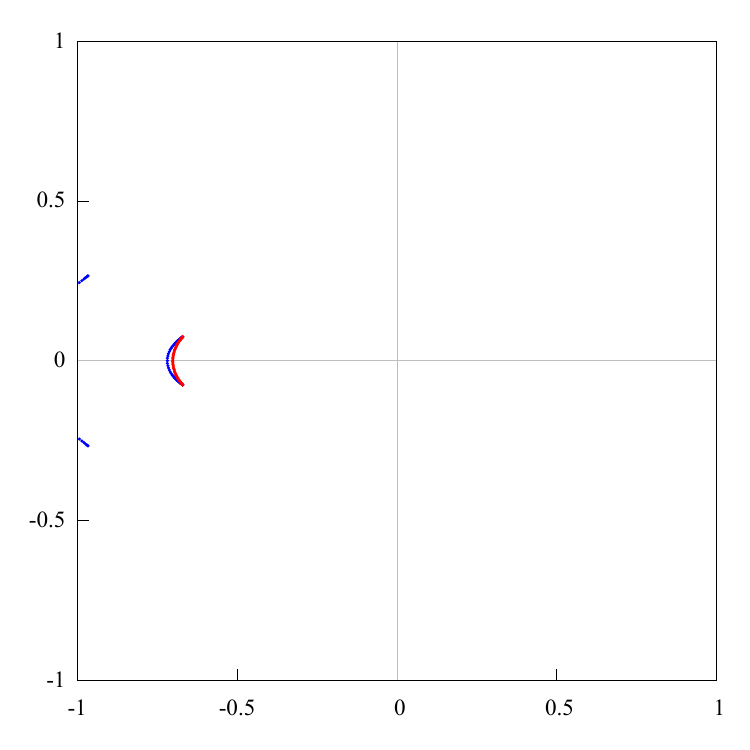}
\vskip-5mm
\caption{We combine Fig.~\ref{fig_PA_100(200)} of the zeros (blue points) of the PA $[200/200]_{E_{(100)}}$ and Fig.~\ref{fig_PA_020(200)} of the poles (red points) of the PA $[200/200]_{E_{(020)}}$ to find out that there are exactly two resonance points, $a$ and $\myo{a}$, for the states $(100)$ and $(020)$.}
\label{fig_PA_100_020(200)}
% \end{center}
% \end{figure}

% \begin{figure}
% \begin{center}
\includegraphics[width=0.75\textwidth]{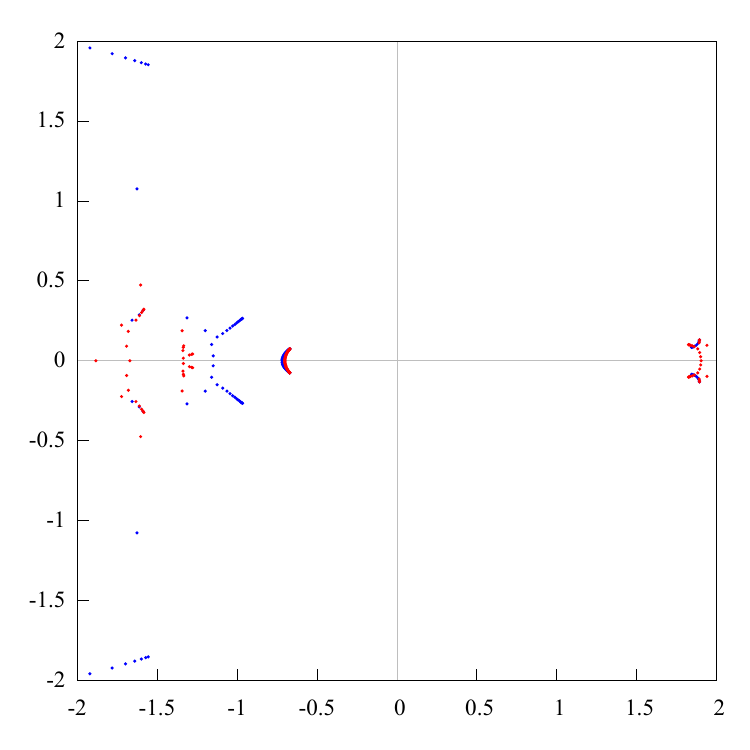}
\vskip-5mm
\caption{The parts of the zeros (blue points) of the PA $[200/200]_{E_{(100)}}$ and the poles (red poles) of the PA $[200/200]_{E_{(020)}}$ are plotted in twice less scale than in Fig.~\ref{fig_PA_100_020(200)}. These parts belong to the $xy$-square $[-2,2]\times[-2,2]$.}
\label{fig_PA_100_020(200)2}
\end{center}
\end{figure}

In Fig.~\ref{fig_PA_100(200)} the part of the zeros (blue points) of the PA $[200/200]_f$ for the function $f=E_{(100)}(z)$ is plotted. Similarly in Fig.~\ref{fig_PA_020(200)} the part of the poles (red points) of the PA $[200/200]_f$ for $f=E_{(020)}(z)$ is plotted. These are the parts that
belong to the $xy$-square $[-1,1]\times[-1,1]$. Finally in Fig.~\ref{fig_PA_100_020(200)} we combine the figures \ref{fig_PA_100(200)} and \ref{fig_PA_020(200)} to find out that there are exactly two resonance points, $a$ and $\myo{a}$, for the states $(100)$ and $(020)$. Namely these $a$ and $\myo{a}$ are the common end points for the arcs of Stahl's compact sets for the functions and $E_{(100)}(z)$ and $E_{(020)}(z)$ that
belong to the $xy$-square $[-1,1]\times[-1,1]$, i.e. they are the common end points of blue and red arcs.
Notice that outside the unit disk there are additional branch points common for the states $(100)$ and $(020)$; see Fig.~\ref{fig_PA_100_020(200)2}. However, these points are not resonance ones.

%%% Fig 29-38

\begin{figure}
\begin{center}
\includegraphics[width=0.75\textwidth]{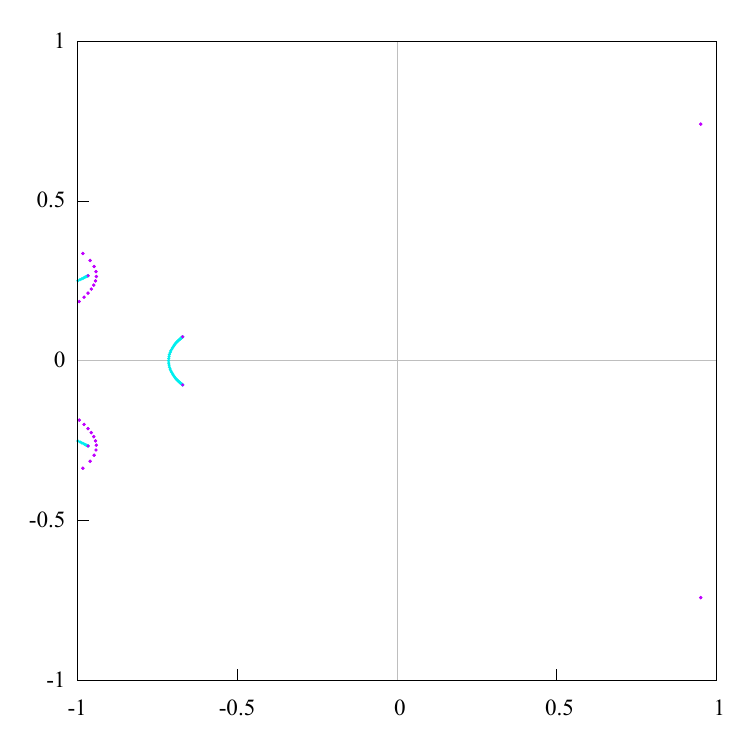}
\vskip-5mm
\caption{The part of the zeros (light blue points) of the type II HP polynomial $P_{266,0}$ as well as the part of the zeros (violet points) of the discriminant polynomial $D_{133}$ for the energy function $E_{(100)}(z)$.}
\label{fig_HP2_100(133)}
% \end{center}
% \end{figure}

% \begin{figure}
% \begin{center}
\includegraphics[width=0.75\textwidth]{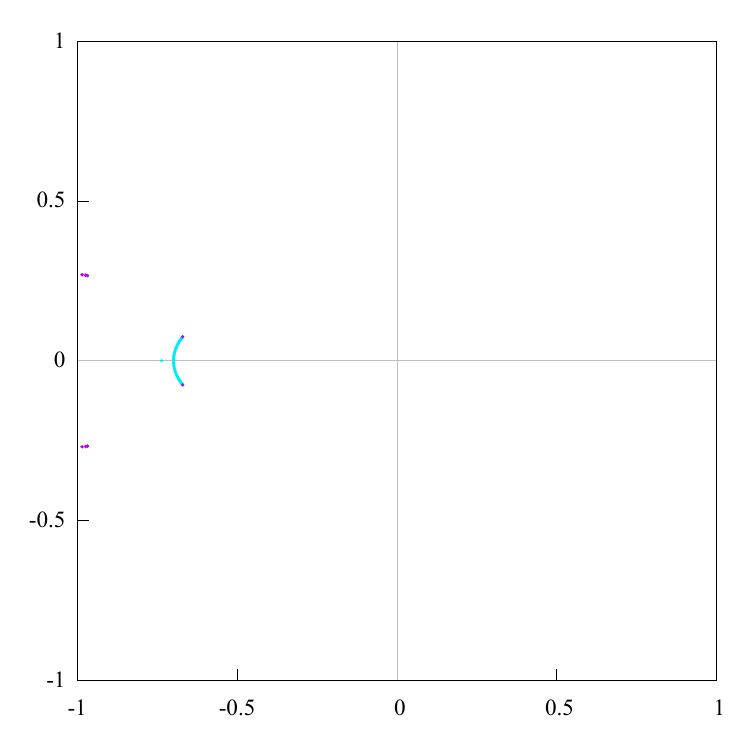}
\vskip-5mm
\caption{The part of the zeros (light blue points) of the type II HP polynomial $P_{266,0}$ as well as the part of the zeros (violet points) of the discriminant polynomial $D_{133}$ for the energy function $E_{(020)}(z)$.}
\label{fig_HP2_020(133)}
\end{center}
\end{figure}

\begin{figure}
\begin{center}
\includegraphics[width=0.75\textwidth]{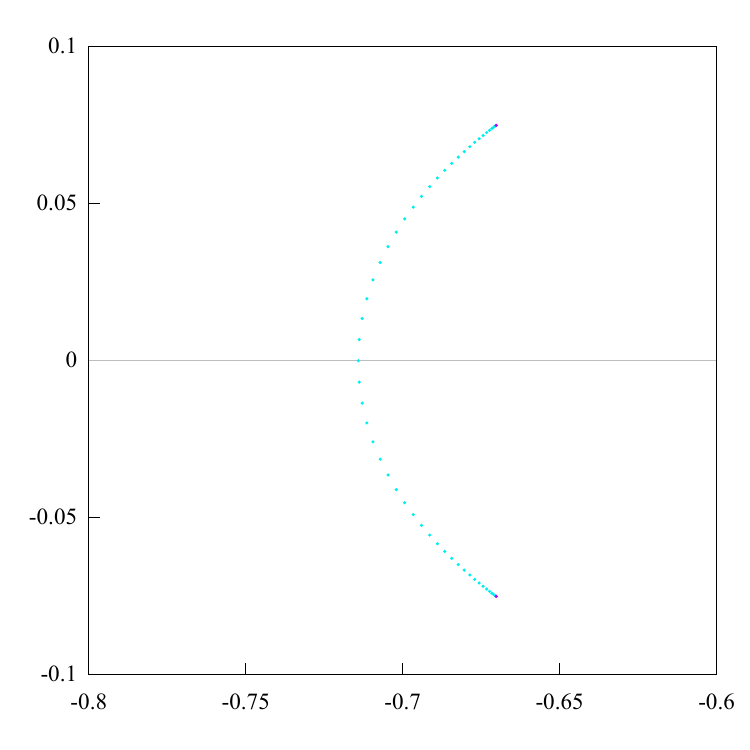}
\vskip-5mm
\caption{The part of the Fig.~\ref{fig_HP2_100(133)} in a different scale.}
\label{fig_HP2_100(133)2}
% \end{center}
% \end{figure}

% \begin{figure}
% \begin{center}
\includegraphics[width=0.75\textwidth]{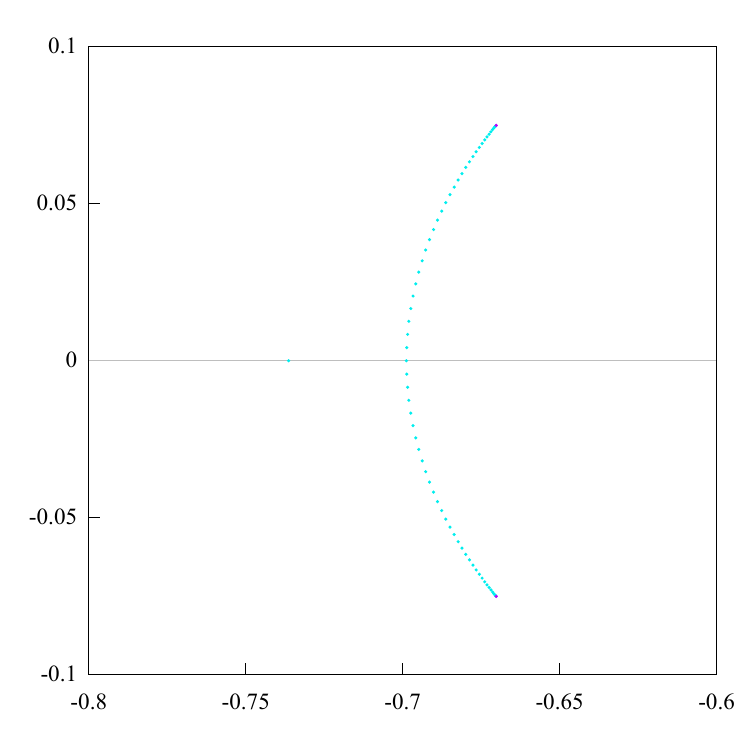}
\vskip-5mm
\caption{The part of the Fig.~\ref{fig_HP2_020(133)} in a different scale.}
\label{fig_HP2_020(133)2}
\end{center}
\end{figure}

\begin{figure}
\begin{center}
\includegraphics[width=0.75\textwidth]{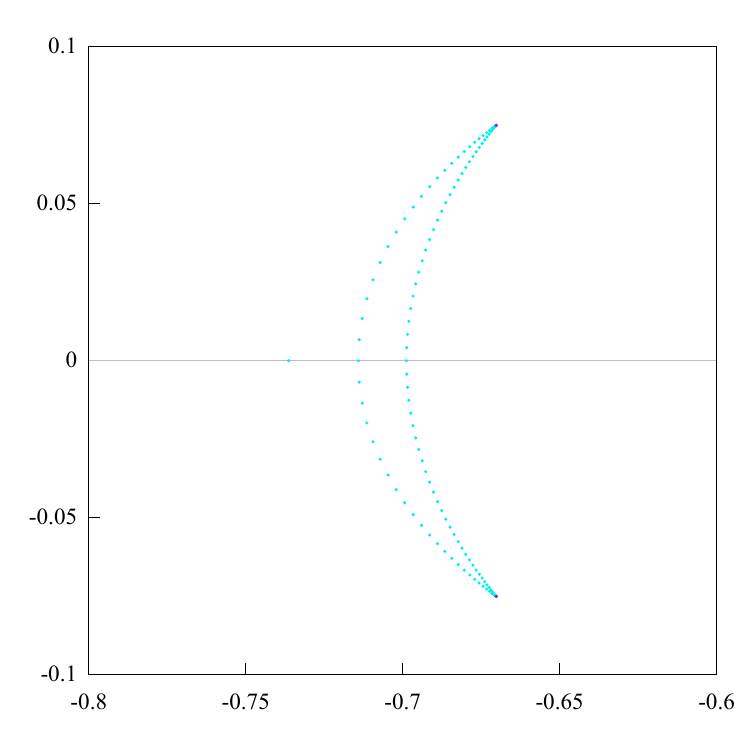}
\vskip-5mm
\caption{We combine Fig.~\ref{fig_HP2_100(133)2} and Fig.~\ref{fig_HP2_020(133)2}.}
\label{fig_HP2_100_020(133)}
% \end{center}
% \end{figure}

% \begin{figure}
% \begin{center}
\includegraphics[width=0.75\textwidth]{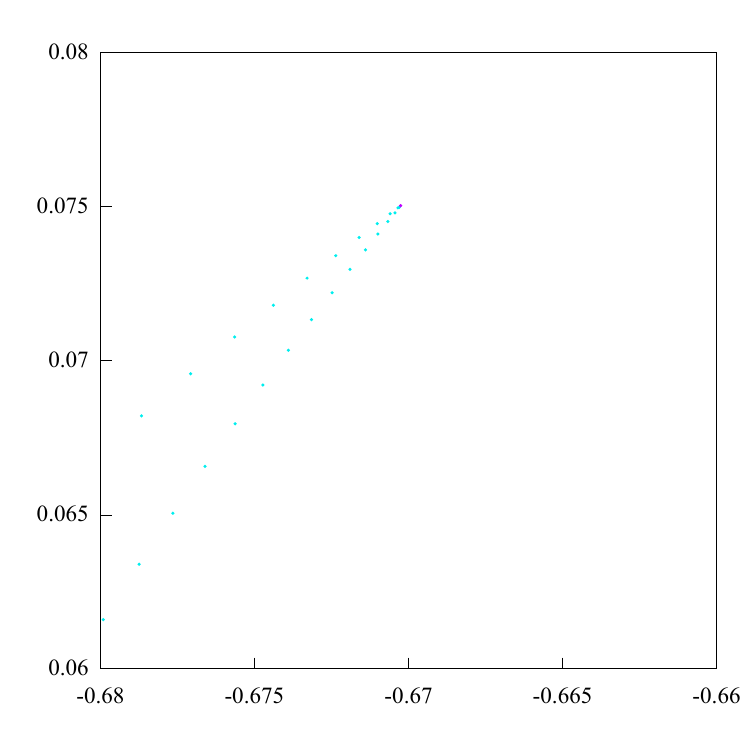}
\vskip-5mm
\caption{The part of Fig.~\ref{fig_HP2_100_020(133)} in a different scale.}
\label{fig_HP2_100_020(133)2}
\end{center}
\end{figure}

%%% Table 1-5
\begin{table}
\begin{center}
\begin{tabular}{|p{0.05\linewidth}|p{0.4\linewidth}|p{0.4\linewidth}|}\hline
 $m$& State $(100)$, $\Re a$&State $(100)$, $\Im a$\\\hline
132 & \textbf{-0.67026119461016981952140}&
\textbf{0.07503319816291785587446}\\\hline
133 & \textbf{-0.67026119461016981952140}&
\textbf{0.07503319816291785587446}\\\hline
\end{tabular}
\caption{State (100), the zeros of the discriminant polynomials $D_{132}$ and $D_{133}$ (of degrees $264$ and $266$ respectively) which correspond to Katz's points $a$ and $\myo{a}$ under consideration. Here and in what follows we mark in bold matching significant digits in computing results.}
\label{tab_100_dis}
% \end{center}
% \end{table}

\vspace{7mm}

% \begin{table}
% \begin{center}
\begin{tabular}{|p{0.05\linewidth}|p{0.4\linewidth}|p{0.4\linewidth}|}\hline
 $m$&State $(020)$, $\Re b$&State $(020)$, $\Im b$\\\hline
132 & \textbf{-0.67026119461016981952140}&
\textbf{0.07503319816291785587446}\\\hline
133 & \textbf{-0.67026119461016981952140}&
\textbf{0.07503319816291785587446}\\\hline
\end{tabular}
\caption{State (020), the zeros of the discriminant polynomials $D_{132}$ and $D_{133}$ (of degrees $264$ and $266$ respectively) which correspond to Katz's points $b$ and $\myo{b}$ under consideration. In view of this table and Table~\ref{tab_100_dis} it is very likely that $a=b$.}
\label{tab_020_dis}
% \end{center}
% \end{table}

\vspace{7mm}

% \begin{table}
% \begin{center}
\begin{tabular}{|p{0.05\linewidth}|p{0.4\linewidth}|p{0.4\linewidth}|}\hline
 $m$&State $(100)$, $\Re a$&State $(100)$, $\Im a$\\\hline
132 & \textbf{-0.670}60837331030142464500 &
\textbf{0.0747}6845430134830985743 \\\hline
133 & \textbf{-0.670}59809849194255074000&
\textbf{0.0747}7777945085218412696\\\hline
\end{tabular}
\caption{State (100), the zeros of type II HP-polynomials $P_{264,0}$ and $P_{266,0}$ (of degrees $266$ and $264$ respectively) which corresponds to the end points of the first plate $E$ of Nuttall's condenser, i.e. to Katz's points.}
\label{tab_100_hp2}
% \end{center}
% \end{table}

\vspace{7mm}

% \begin{table}
% \begin{center}
\begin{tabular}{|p{0.05\linewidth}|p{0.4\linewidth}|p{0.4\linewidth}|}\hline
 $m$&State $(020)$, $\Re b$&State $(020)$, $\Im a$\\\hline
132 & \textbf{-0.670}99012883492243843996 &
\textbf{0.074}10300601602209254537 \\\hline
133 & \textbf{-0.670}30716286939747295259&
\textbf{0.074}97625147484154139891\\\hline
\end{tabular}
\caption{State (020), the zeros of the II HP-polynomials $P_{264,0}$ and $P_{266,0}$ (of degrees $264$ and $266$ respectively) which correspond to the end point of the first plate $E$ of Nuttall's condenser, i.e. to Katz's point.}
\label{tab_020_hp2}
% \end{center}
% \end{table}

\vspace{7mm}

% \begin{table}
% \begin{center}
\begin{tabular}{|p{0.05\linewidth}|p{0.4\linewidth}|p{0.4\linewidth}|}\hline
 $m$& State $(300)$, $\Re a_1$ &State $(300)$, $\Im a_1$\\\hline
132 & \textbf{-0.533443}51574580692125444&
\textbf{0.762229}90907521253553984\\\hline
133 & \textbf{-0.533443}34795658276857734&
\textbf{0.762229}70415991163293731\\\hline
 $m$& State $(102)$, $\Re b_1$ &State $(102)$, $\Im b_1$\\\hline
132 & \textbf{-0.53344}299528745937443528&
\textbf{0.762229}83076337916506589\\\hline
133 & \textbf{-0.53344}319272109942749749&
\textbf{0.762229}96892599080028438\\\hline
\end{tabular}
\caption{State $(300)$, the zeros of the discriminant polynomials $D_{132}$ and $D_{133}$ (of degrees $264$ and $266$ respectively) which correspond to Katz's points $a_1$ and $\myo{a}_1$ of $E_{(300)}$. State $(102)$, the zeros of the discriminant polynomials $D_{132}$ and $D_{133}$ (of degrees $264$ and $266$ respectively) which correspond to Katz's points $b_1$ and $\myo{b}_1$ of $E_{(102)}$.}
\label{tab_300_102_dis1}
\end{center}
\end{table}

%%% Table 6-9
\begin{table}
\begin{center}
\begin{tabular}{|p{0.05\linewidth}|p{0.4\linewidth}|p{0.4\linewidth}|}\hline
 $m$& State $(300)$, $\Re a_2$ &State $(300)$, $\Im a_2$\\\hline
132 & \textbf{-0.26829795679962550141549}&
\textbf{0.19161346883586603091241}\\\hline
133 & \textbf{-0.26829795679962550141549}&
\textbf{0.19161346883586603091241}\\\hline
 $m$& State $(102)$, $\Re b_2$ &State $(102)$, $\Im b_2$\\\hline
132 & \textbf{-0.26829795679962550141549}&
\textbf{0.19161346883586603091241}\\\hline
133 & \textbf{-0.26829795679962550141549}&
\textbf{0.19161346883586603091241}\\\hline
\end{tabular}
\caption{State $(300)$, the zeros of the discriminant polynomials $D_{132}$ and $D_{133}$ (of degrees $264$ and $266$ respectively) which correspond to Katz's points $a_2$ and $\myo{a}_2$ of $E_{(300)}$. State $(102)$, the zeros of the discriminant polynomials $D_{132}$ and $D_{133}$ (of degrees $264$ and $266$ respectively) which correspond to Katz's points $b_2$ and $\myo{b}_2$ of $E_{(102)}$.}
\label{tab_300_102_dis2}
% \end{center}
% \end{table}

\vspace{5mm}

% \begin{table}
% \begin{center}
\begin{tabular}{|p{0.05\linewidth}|p{0.4\linewidth}|p{0.4\linewidth}|}\hline
 $m$& State $(300)$, $\Re a_3$ &State $(300)$, $\Im a_3$\\\hline
132 & \textbf{-0.346134977}37296942030131&
\textbf{0.038093272}86190358525117\\\hline
133 & \textbf{-0.346134977}23274944106753&
\textbf{0.038093272}94232679467132\\\hline
 $m$& State $(102)$, $\Re b_3$ &State $(102)$, $\Im b_3$\\\hline
132 & \textbf{-0.34353580034967537}678176&
\textbf{0.03830510430931258}595242\\\hline
133 & \textbf{-0.34353580034967537}787002&
\textbf{0.03830510430931258}829996\\\hline
\end{tabular}
\caption{State $(300)$, the zeros of the discriminant polynomials $D_{132}$ and $D_{133}$ (of degrees $264$ and $266$ respectively) which correspond to Katz's points $a_3$ and $\myo{a}_3$ of $E_{(300)}$. State $(102)$, the zeros of the discriminant polynomials $D_{132}$ and $D_{133}$ (of degrees $264$ and $266$ respectively) which correspond to Katz's points $b_3$ and $\myo{b}_3$ of $E_{(102)}$. Thus, it is very probably that $a_3\neq b_3$ (see also Fig.~\ref{fig_HP_2_300_102_d_300_102_z2}).}
\label{tab_300_102_dis3}
% \end{center}
% \end{table}

\vspace{5mm}

% \begin{table}
% \begin{center}
\begin{tabular}{|p{0.05\linewidth}|p{0.4\linewidth}|p{0.4\linewidth}|}\hline
 $m$& State $(201)$, $\Re a_1$ &State $(201)$, $\Im a_1$\\\hline
132 & \textbf{-0.32632295243225535632149}&
\textbf{0.03517114877422451767753}\\\hline
133 & \textbf{-0.32632295243225535632149}&
\textbf{0.03517114877422451767753}\\\hline
 $m$& State $(121)$, $\Re a$ &State $(121)$, $\Im a$\\\hline
132 & \textbf{-0.32632295243225535632149}&
\textbf{0.03517114877422451767753}\\\hline
133 & \textbf{-0.32632295243225535632149}&
\textbf{0.03517114877422451767753}\\\hline
\end{tabular}
\caption{State $(201)$, the zeros of the discriminant polynomials $D_{132}$ and $D_{133}$ (of degrees $264$ and $266$ respectively) which correspond to Katz's points $a_1$ and $\myo{a}_1$ of $E_{(201)}$. State $(121)$, the zeros of the discriminant polynomials $D_{132}$ and $D_{133}$ (of degrees $264$ and $266$ respectively) which correspond to Katz's points $a$ and $\myo{a}$ of $E_{(121)}$.}
\label{tab_201_121_dis}
% \end{center}
% \end{table}

\vspace{5mm}

% \begin{table}
% \begin{center}
\begin{tabular}{|p{0.05\linewidth}|p{0.4\linewidth}|p{0.4\linewidth}|}\hline
 $m$& State $(201)$, $\Re a_2$ &State $(201)$, $\Im a_2$\\\hline
132 & \textbf{-0.4054}5657692987582982683&
\textbf{0.014}73569561762927184216\\\hline
133 & \textbf{-0.4054}1511097862812977399&
\textbf{0.014}83582036989545697530\\\hline
 $m$& State $(041)$, $\Re b$ &State $(041)$, $\Im b$\\\hline
132 & \textbf{-0.40541518613840703096783}&
\textbf{0.01479134185042952294418}\\\hline
133 & \textbf{-0.40541518613840703096783}&
\textbf{0.01479134185042952294418}\\\hline
\end{tabular}
\caption{State $(201)$, the zeros of the discriminant polynomials $D_{132}$ and $D_{133}$ (of degrees $264$ and $266$ respectively) which correspond to Katz's points $a_2$ and $\myo{a}_2$ of $E_{(201)}$. State $(041)$, the zeros of the discriminant polynomials $D_{132}$ and $D_{133}$ (of degrees $264$ and $266$ respectively) which correspond to Katz's points $b$ and $\myo{b}$ of $E_{(041)}$.}
\label{tab_201_041_dis}
\end{center}
\end{table}

%%% Fig 39-43
\begin{figure}
\begin{center}
\includegraphics[width=0.75\textwidth]{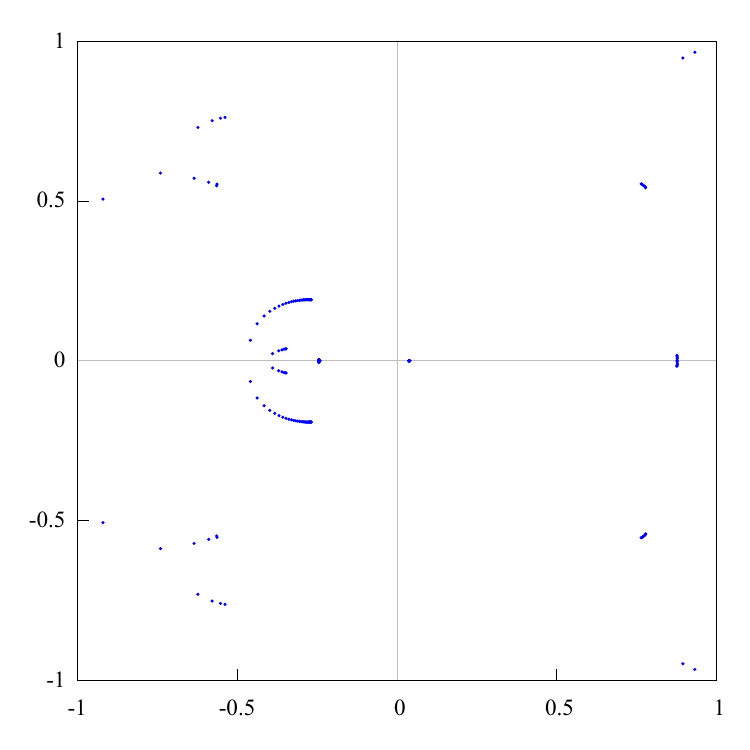}
\vskip-5mm
\caption{The part of the zeros (blue points) of the PA $[200/200]_f$ for the energy function $f=E_{(300)}(z)$ that belongs to the $xy$-square $[-1,1]\times[-1,1]$. According to Stahl's Theory~\cite{Sta97b} the poles of $[200/200]_f$ simulate Stahl's compact set $S=S(f)$ of $f$.}
\label{fig_PA_300}
% \end{center}
% \end{figure}

% \begin{figure}
% \begin{center}
\includegraphics[width=0.75\textwidth]{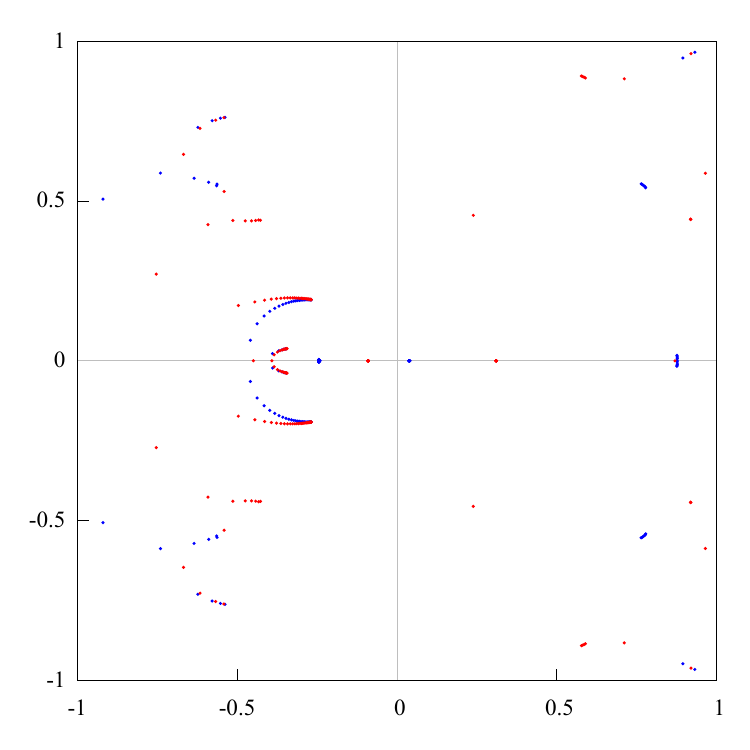}
\vskip-5mm
\caption{The part of the zeros (blue points) of the PA $[200/200]_{E_{(300)}}$ and the part of the poles (red points) of the PA $[200/200]_{E_{(102)}}$ that belongs to the $xy$-square $[-1,1]\times[-1,1]$. From this figure we find out that there are probably three pairs of the resonance points for the states $(300)$ and $(102)$.}
\label{fig_PA_300_102}
\end{center}
\end{figure}

\begin{figure}
\begin{center}
\includegraphics[width=0.75\textwidth]{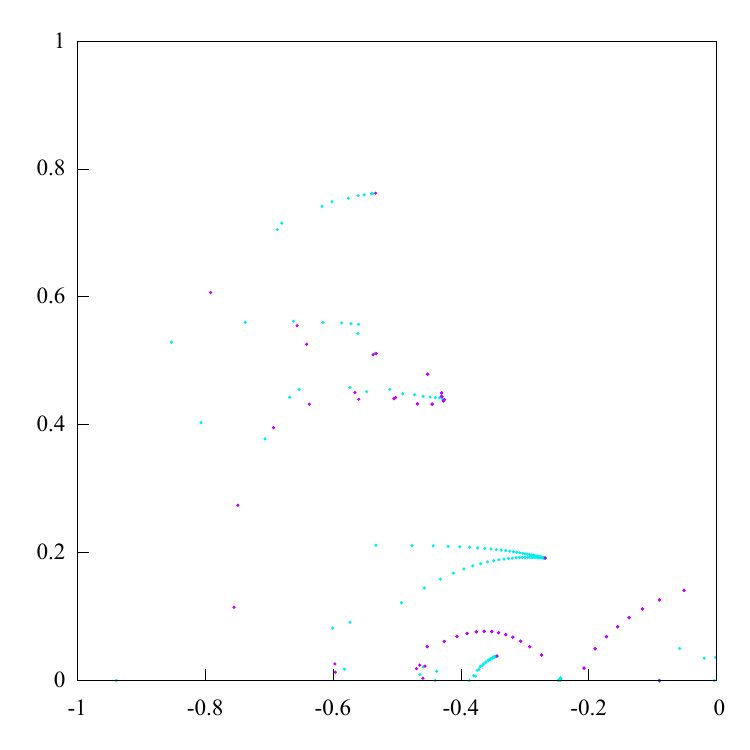}
\vskip-5mm
\caption{The part of the zeros (light blue points) of the type II HP polynomials $P_{266,0}$ for the energy functions $E_{(300)}(z)$ and $E_{(102)}$ as well as the part of the zeros (violet points) of the discriminant polynomial $D_{133}$ for the energy function $E_{(102)}(z)$.}
\label{fig_HP_2_300_102_d_102}
% \end{center}
% \end{figure}

% \begin{figure}
% \begin{center}
\includegraphics[width=0.75\textwidth]{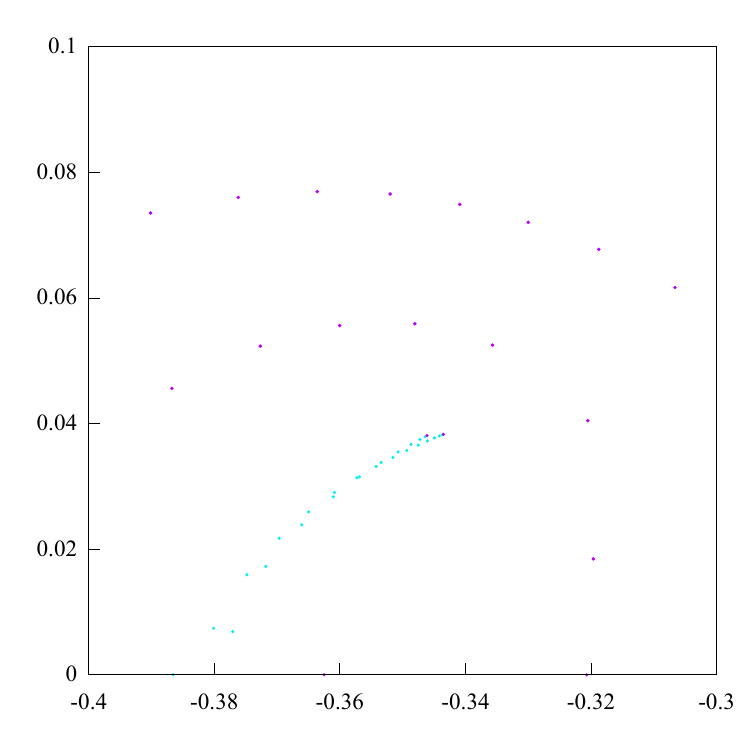}
\vskip-5mm
\caption{The part of the zeros (light blue points) of the type II HP polynomials $P_{266,0}$ for the energy functions $E_{(300)}(z)$ and $E_{(102)}$ as well as the part of the zeros (violet points) of the discriminant polynomials $D_{133}$ for these energy functions.}
\label{fig_HP_2_300_102_d_300_102_z}
\end{center}
\end{figure}

\begin{figure}[t]
\begin{center}
\includegraphics[width=0.7\textwidth]{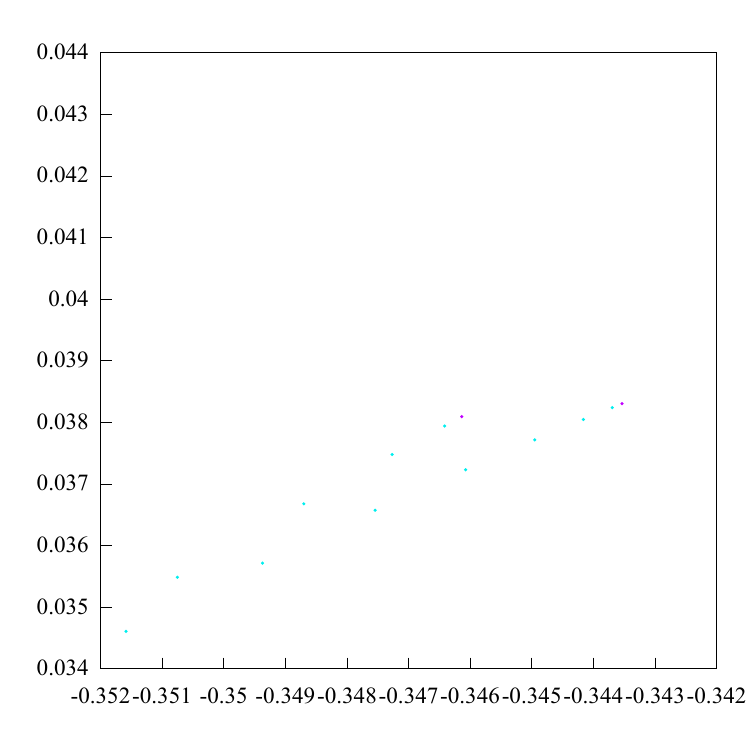}
\vskip-5mm
\caption{The part of Fig.~\ref{fig_HP_2_300_102_d_300_102_z} is represented in another scale. In view of the figure it is very likely that $a_3\neq b_3$ (see also Table~\ref{tab_300_102_dis3}).}
\label{fig_HP_2_300_102_d_300_102_z2}
\end{center}
\vspace{-5mm}
\end{figure}

The second step of the algorithm is intended to locate these two Katz's points and after that to compute them as precisely as possible. For this purpose we apply type II HP polynomials $P_{266,0}$ and the discriminant polynomials $D_{133}$ for both states $(100)$ and $(020)$.
In accordance with the results of~\cite{KoKrPaSu16},~\cite{KoPaSuCh17} and~\cite{Sue22b}, the limit zero distribution of
type II HP polynomials $P_{2m,0}$ corresponds to the plate $E$ of Nuttall's condenser $(E,F)$ which consists of a finite number of analytic arcs. In particular, zeros of $P_{2m,0}$ simulate these arcs of $E$ and mark their end points. However the corresponding zeros of $P_{2m,0}$ converge to the end points very slowly, namely, at a power-law rate.
In their part, the zeros of the discriminant polynomials $D_{m}$ simulate the second plate $F$ of Nuttall's condenser. However there are a finite number of zeros of $D_{m}$ which mark as $m\to\infty$ the end points of $E$. These zeros of $D_m$ converge to the corresponding end points of $E$ very fast, namely at exponential rate~\cite{Sue22b}.

In Fig.~\ref{fig_HP2_100(133)} the zeros of the type II HP polynomial $P_{266,0}$ (light blue points) and the zeros of the discriminant polynomial $D_{133}$ (violet points) for the function $E_{(100)}(z)$ are plotted. Similarly in Fig.~\ref{fig_HP2_020(133)} the zeros of type II HP polynomial $P_{266,0}$ (light blue points) and the zeros of discriminant polynomial $D_{133}$ (violet points) for the function $E_{(020)}(z)$ are plotted.

In Fig.~\ref{fig_HP2_100(133)2} we present a part of Fig.~\ref{fig_HP2_100(133)} in a different scale. In Fig.~\ref{fig_HP2_020(133)2} we present a part of the Fig.~\ref{fig_HP2_020(133)} but also in a different scale. In
Fig.~\ref{fig_HP2_100_020(133)} we combine the figures~\ref{fig_HP2_100(133)2} and~\eqref{fig_HP2_020(133)2} and in Fig.~\ref{fig_HP2_100_020(133)2} we represent a part of Fig.~\ref{fig_HP2_100_020(133)} but again in a different scale. The last figure gives us an approximate location of one of the two Katz's points -- the end of the $E$-plate of Nuttall's condenser. Namely, from Fig.~\ref{fig_HP2_100_020(133)2} it follows that within the boundaries of the $xy$-square specified in this figure, there is only one zero of the discriminant polynomial $D_{133}$.

Given the polynomial $D_{133}$, this zero can be computed very easy as well as the corresponding zero of the polynomial $P_{266,0}$. To understand the accuracy 
of the computations, we combine into Table~\ref{tab_100_dis} the results of computations of the desired zero of the discriminant polynomials $D_{132}$ and $D_{133}$ corresponding to the state $(100)$. From the table it directly follows that the rate of convergence is exponential. Just the same conclusion follows from Table~\ref{tab_020_dis} about the desired zero of the discriminant polynomials $D_{132}$ and $D_{133}$ corresponding to the state $(020)$. In contrast to those facts, the convergence of the zeros of the type II HP polynomials is at power-law rate; see
Tables~\ref{tab_100_hp2} and~\ref{tab_020_hp2}.

\subsection{\texorpdfstring{Case of the states $\myk=(300)$ and $\myj=(102)$}{Case of the states k=(300) and j=(102)}}\label{s5s2}

To analyze the existence of resonance points for the states $\myk=(300)$ and $\myj=(102)$ we present in this section the arguments similarly to Sec.~\ref{s5s1}.

To find out if there are resonance points for the energies $E_{(300)}$ and $E_{(102)}$ of the states $(300)$ and $(102)$, we apply the PA $[200/200]$ as a first step to analyze the problem; see Fig.~\ref{fig_PA_300} and Fig.~\ref{fig_PA_300_102}. To be more precise, in the figures~\ref{fig_PA_300} and
~\ref{fig_PA_300_102} a part of the zeros (blue points) of the PA $[200/200]_{E_{(300)}}$ is depicted as well as a part of the poles (red points) of the PA $[200/200]_{E_{(102)}}$. From these figures it follows that there are probably three pairs of the resonance points to the states $(300)$ and $(102)$. All these probable resonance points locate to the left of the origin of coordinates. Thus following the algorithm proposed above in Sec.~\ref{s5s1}, we restrict our attention to the square $[-1,0]\times[0,1]$ in the $xy$-plane instead of the original square $[-1,1]\times[-1,1]$; see Fig.~\ref{fig_HP_2_300_102_d_102}. From this figure it follows that there are probably three pairs of resonance points. We denote these candidates for resonance points by $a_1,a_2,a_3$ for the state $E_{(300)}$ and by $b_1,b_2,b_3$ for the state $E_{(102)}$, $\Im a_j,\Im b_j>0$. So the problem is to check if $a_j=b_j$ or not.

This validation is doing
 by the second step of the algorithm which is based on the
type II HP polynomial $P_{266,0}$ and the discriminant polynomial $D_{133}$ of degree $266$ for the functions $E_{(300)}$ and $E_{(102)}$.Thus we localize all the probable resonance points and after that compute them as precisely as possible. The results are presented in
Table~\ref{tab_300_102_dis1}, Table~\ref{tab_300_102_dis2} and Table~\ref{tab_300_102_dis3} and in
Fig.~\ref{fig_HP_2_300_102_d_300_102_z} and Fig.~\ref{fig_HP_2_300_102_d_300_102_z2}. So it definitely follows from figures~\ref{fig_HP_2_300_102_d_300_102_z} and~\ref{fig_HP_2_300_102_d_300_102_z2} and Table~\ref{tab_300_102_dis3} that $a_3\neq b_3$. Therefore there are only two resonance pairs for the states $(300)$ and $(102)$ presented in the tables~\ref{tab_300_102_dis1} and~\ref{tab_300_102_dis2}.

\subsection{\texorpdfstring{Case of the states $\myk_1=(201)$, $\myk_2=(121)$ and $\myk_3=(041)$}{Case of the states k1=(201), k2=(121) and k3=(041)}}\label{s5s3}

Again, to analyze the existence of resonance points for the states $\myk_1=(201)$, $\myk_2=(121)$ and $\myk_3=(041)$ we present in this section the arguments similarly to Sec.~\ref{s5s1}.

To find out if there are resonance points for the energies $E_{(201)}$, $E_{(121)}$ and $E_{(041)}$ of the states $(201)$, $(121)$ and $(041)$, we apply the PA $[200/200]$ as a first step for analysis of the problem. Following the algorithm proposed above in Sec.~\ref{s5s1}, we come to the square $[-0.5,-0.3]\times[-0.1,0.1]$ in the $xy$-plane instead of the original square $[-1,1]\times[-1,1]$.

\begin{figure}
\begin{center}
\includegraphics[width=0.75\textwidth]{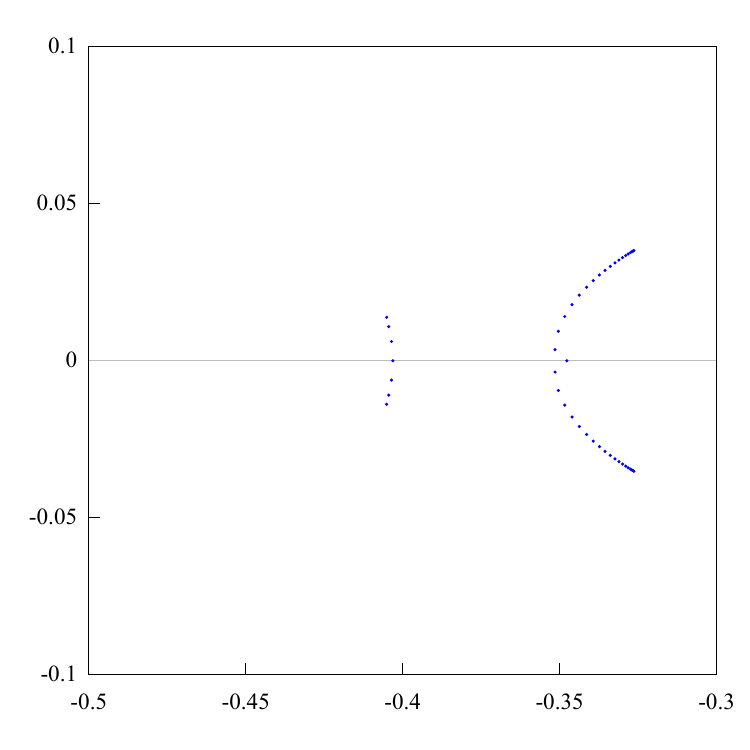}
\vskip-5mm
\caption{The part of the zeros (blue points) of the PA $[200/200]_f$ for the energy function $f=E_{(201)}(z)$. According to Stahl's Theory~\cite{Sta97b} the poles of $[200/200]_f$ simulate Stahl's compact set $S=S(f)$ of $f$.}
\label{fig_PA_201}
% \end{center}
% \end{figure}

% \begin{figure}
% \begin{center}
\includegraphics[width=0.75\textwidth]{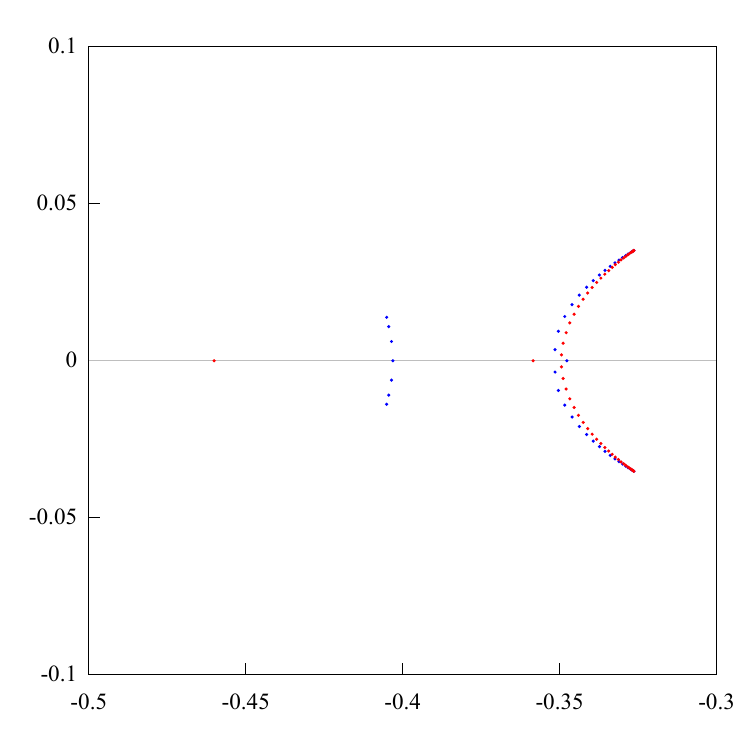}
\vskip-5mm
\caption{The part of the zeros (blue points) of the PA $[200/200]_{E_{(201)}}$ and the part of the poles (red points) of the PA $[200/200]_{E_{(121)}}$. From this figure we find out that there are exactly two resonance points, $a_1=a$ and $\myo{a}_1=a$, for the states $(201)$ and $(121)$.}
\label{fig_PA_201_121}
\end{center}
\end{figure}

\begin{figure}
\begin{center}
\includegraphics[width=0.75\textwidth]{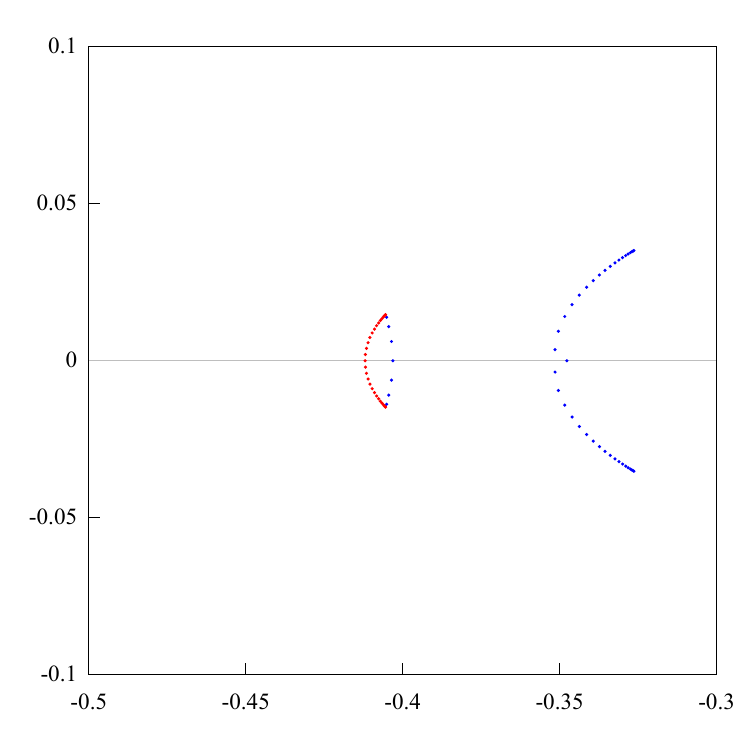}
\vskip-5mm
\caption{The part of the zeros (blue points) of the PA $[200/200]_{E_{(201)}}$ and a part of the poles (red points) of the PA $[200/200]_{E_{(041)}}$. From this figure we find out that there are exactly two resonance points, $a_2=b$ and $\myo{a}_2=\myo{b}$, for the states $(201)$ and $(041)$.}
\label{fig_PA_201_041}
% \end{center}
% \end{figure}

% \begin{figure}
% \begin{center}
\includegraphics[width=0.75\textwidth]{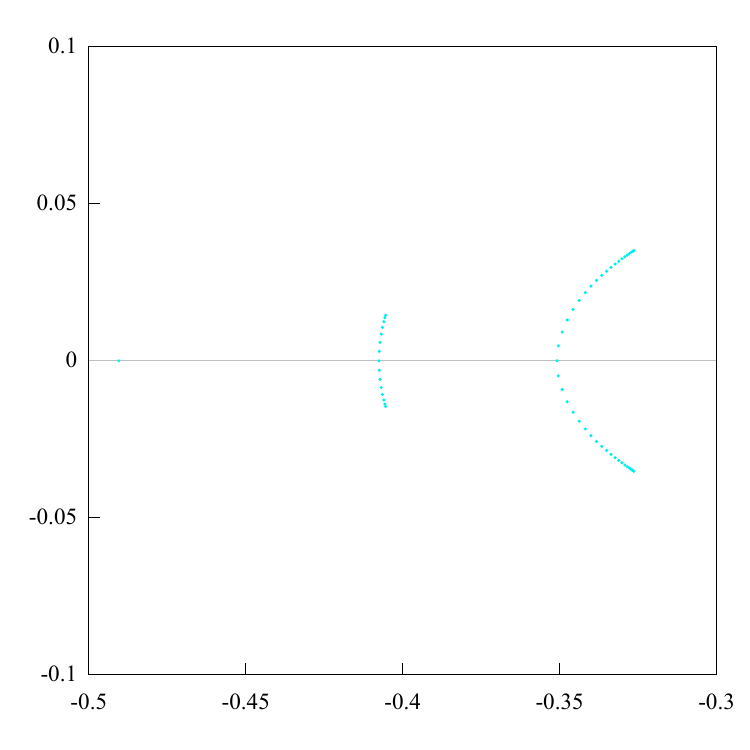}
\vskip-5mm
\caption{The part of the zeros (light blue points) of the type II HP polynomial $P_{266,0}$ for the energy function $E_{(201)}(z)$.}
\label{fig_HP_2_201}
\end{center}
\end{figure}

In the figures~\ref{fig_PA_201},~\ref{fig_PA_201_121} and
~\ref{fig_PA_201_041} the part of the zeros (blue points) of the PA $[200/200]_{E_{(201)}}$ is depicted as well as the part of the poles (red points) of the PA $[200/200]_{E_{(121)}}$ and the PA $[200/200]_{E_{(041)}}$. From these figures it follows that there are probably two pairs of the resonance points to the state $(201)$: one pair that is common to the state $(121)$ and the second pair that is common to the state $(041)$.

\begin{figure}
\begin{center}
\includegraphics[width=0.75\textwidth]{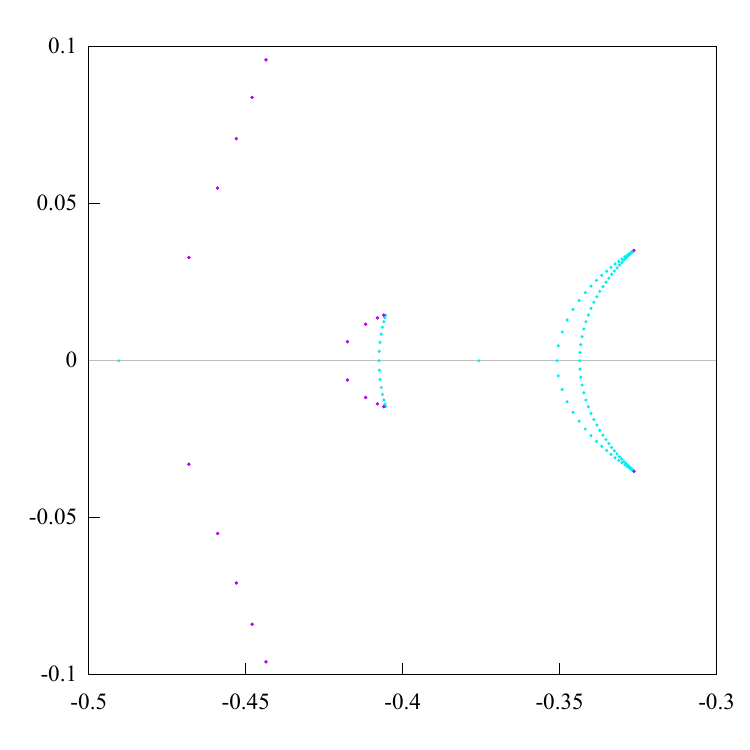}
\vskip-5mm
\caption{The part of the zeros (light blue points) of the type II HP polynomial $P_{266,0}$ for the energy function $E_{(201)}(z)$ as well as the part of the zeros (light blue points) of the type II HP polynomial $P_{266,0}$ and the part of the zeros (violet points) of the discriminant polynomial $D_{133}$ for the energy function $E_{(121)}(z)$ (cf. Fig.~\ref{fig_HP_2_201}).}
\label{fig_HP_2_201_121_D_121}
% \end{center}
% \end{figure}

% \begin{figure}
% \begin{center}
\includegraphics[width=0.75\textwidth]{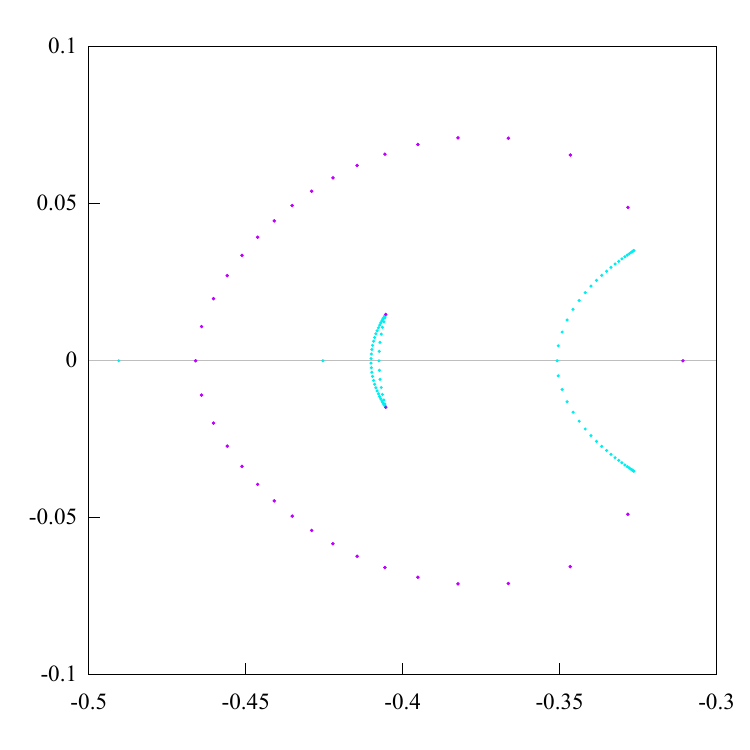}
\vskip-5mm
\caption{The part of the zeros (light blue points) of the type II HP polynomial $P_{266,0}$ for the energy function $E_{(201)}(z)$ as well as the part of the zeros (light blue points) of the type II HP polynomial $P_{266,0}$ and the part of the zeros (violet points) of the discriminant polynomial $D_{133}$ for the energy function $E_{(041)}(z)$ (cf. Fig.~\ref{fig_HP_2_201}).}
\label{fig_HP_2_201_041_D_041}
\end{center}
\end{figure}

These preliminary results can be improved and detailed by the second step of the algorithm which is based on the
type II HP polynomial $P_{266,0}$ and the discriminant polynomial $D_{133}$ of degree $266$ for the functions $E_{(201)},E_{(121)}$ and $E_{(041)}$; see figures~\ref{fig_HP_2_201},~\ref{fig_HP_2_201_121_D_121}
and~\ref{fig_HP_2_201_041_D_041}. Thus we localize these four resonance points and after that compute them as precisely as possible. The results are presented in Tables~\ref{tab_201_121_dis} and~\ref{tab_201_041_dis}.

\begin{figure}
\begin{center}
\includegraphics[width=0.75\textwidth]{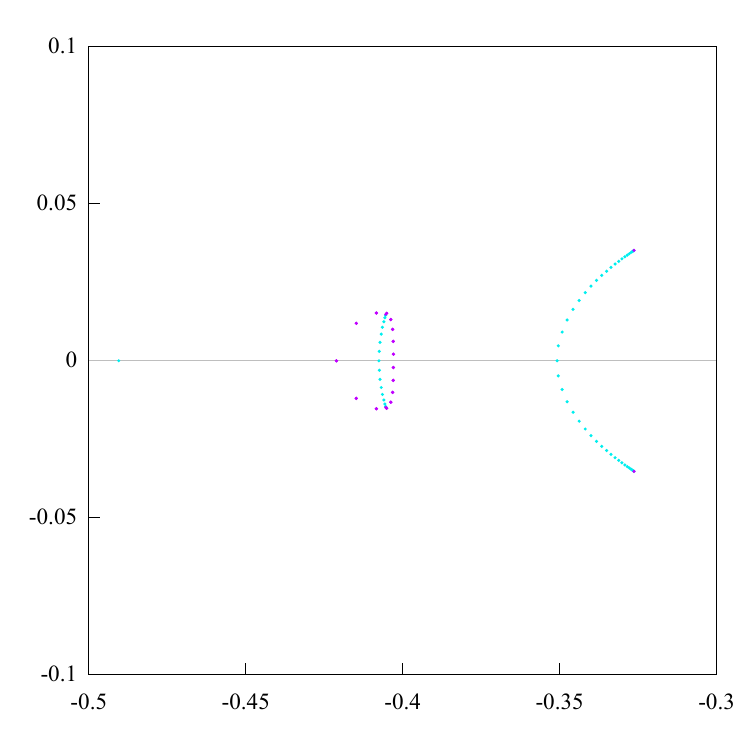}
\vskip-5mm
\caption{The part of the zeros (light blue points) of the type II HP polynomial $P_{266,0}$ as well as the part of the zeros (violet points) of the discriminant polynomial $D_{133}$ for the energy function $E_{(201)}(z)$.}
\label{fig_HP_2_201_z}
% \end{center}
% \end{figure}

% \begin{figure}
% \begin{center}
\includegraphics[width=0.75\textwidth]{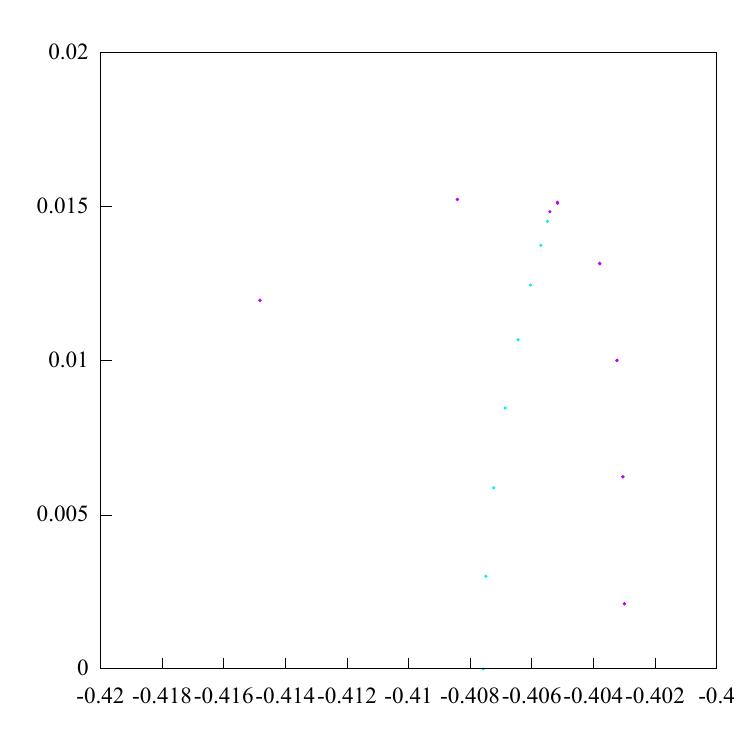}
\vskip-5mm
\caption{The part of the zeros (light blue points) of the type II HP polynomial $P_{266,0}$ as well as the part of the zeros (violet points) of the discriminant polynomial $D_{133}$ for the energy function $E_{(201)}(z)$ in a different scale than in Fig.~\ref{fig_HP_2_201_z}.}
\label{fig_HP_2_201_z2}
\end{center}
\end{figure}

Let us give some comments to these tables. In view of Table~\ref{tab_201_121_dis} it is very likely that the resonance point $a_1=a$ is recognized via the algorithm very reliably. The significant digits of all the appropriate zeros of the discriminant polynomials are stabilizing via the computations very fast, i.e. up to 23 significant digits at least (in fact, even more significant digits are stabilized). In contrast to that, the situation in Table~\ref{tab_201_041_dis} is very different from the above. Namely, the significant digits of the appropriate zeros of the two discriminant polynomials $D_{132}$ and $D_{133}$ for the state $E_{(041)}$ are stabilized also very fast. However it is not the case for the state $E_{(201)}$. The reason for that is the following. As it was explained above, zeros of the type II HP polynomials accumulate to the first plate~$E$ of Nuttall's condenser $(E,F)$. At the same time, zeros of
type I HP polynomials as well as zeros of discriminant polynomials accumulate to the second plate~$F$.
In view of Fig.~\ref{fig_HP_2_201_z} and Fig.~\ref{fig_HP_2_201_z2} it is very likely that in the case of the state $E_{(201)}$ these two plates $E$ and $F$ are very close to each other. This is the reason for a very slow convergence of the appropriate zeros of the discriminant polynomials in the case under consideration.

%%% Fig 44-51

% \newpage
% \clearpage

\section{Some theoretical results}\label{s6}

In this section, we present some theoretical results to support our point of view about the advantage of the HP-analysis over the PA-analysis.

Let a function $f_{*}(\zeta)\in\HH(\infty)$ be given by the following explicit representation
\begin{equation}
f_{*}(\zeta):=\[\(A-\frac1{\varphi(\zeta)}\)\(B-\frac1{\varphi(\zeta)}\)\]^{-1/2},
\quad \zeta\notin[-1,1],
\label{funzhu}
\end{equation}
where $1<A<B$ and such branch of the function $(\cdot)^{1/2}$ is chosen that $\varphi(\zeta)=\zeta+(\zeta^2-1)^{1/2}\sim 2\zeta$ and $f_*(\zeta)\sim1/\sqrt{AB}$
as $\zeta\to\infty$. The function $f_{*}$ is an algebraic function of fourth degree with four branch points $\{\pm1,a,b\}$, where $a=(A+1/A)/2$, $b=(B+1/B)/2$, $1<a<b$. The segment $E=[-1,1]$ is Stahl's compact set $S(f_*)$ for $f_{*}$ given by~\eqref{funzhu} and the domain $D:=\myh{\CC}\setminus{E}$ is the corresponding Stahl's domain. In~\cite{Sue18c} it is proved that $f_{*}$ is a Markov type function and the pair $f_*,f_*^2$ forms a Nikishin system~\cite{NiSo88}:
\begin{align*}
f_*(\zeta)&=\frac1{\sqrt{AB}}+\myh{\sigma}(\zeta),\quad \supp{\sigma}=E, \\
f_*^2(\zeta)&=\frac1{AB}+\frac1{\sqrt{AB}}\myh{\sigma}(\zeta) +\myh{s}(\zeta),\quad \supp{s}=E,
\end{align*}
where $ds(\zeta)=\myh{\sigma}_2(\zeta)\,d\sigma(\zeta)$, $\supp{\sigma_2}=[a,b]$,
and the measures $\sigma$ and~$\sigma_2$
have the following explicit representations
\begin{align}
d\sigma(x):&=
\frac{\sqrt{1-x^2_1}}{4\pi\sqrt{AB}\sqrt{(a-x)(b-x)}}
\[\frac{h_2(x)}{h_1(x)}+\frac{h_1(x)}{h_2(x)}\]\,dx,
\quad x\in[-1,1],
\label{zhumes1}\\
d\sigma_2(t):&=
\frac1{\pi}\frac{dt}{\sqrt{(\pfi(t)-A)(B-\pfi(t))}},\quad
t\in(a,b),
\label{zhumes2}
\end{align}
with
\begin{align*}
h_1(x):&=\(A-(x+i\sqrt{1-x^2})\)^{1/2}+\(A-(x-i\sqrt{1-x^2})\)^{1/2}>0, \ x\in[-1,1], \\
h_2(x):&=\(B-(x+i\sqrt{1-x^2})\)^{1/2}+\(B-(x-i\sqrt{1-x^2})\)^{1/2}>0, \ x\in[-1,1]. 
\end{align*}
Recall that
\begin{equation*}
\myh{\varkappa}(\zeta):=\int\frac{d\varkappa(x)}{\zeta-x},\quad \zeta\in\myh{\CC}\setminus\supp{\varkappa},\quad
\supp{\varkappa}\Subset\RR.
\end{equation*}

For a probability measure $\mu$, $\supp{\mu}\Subset\RR$, let
\begin{equation}
V^\mu(\zeta):=\int\log\frac1{|\zeta-t|}\,d\mu(t)
\label{logpot}
\end{equation}
be the logarithmic potential of $\mu$. Let $M_1(E)$ be the set of all probability measures supported on $E=[-1,1]$.
There exists~\cite{Lan66} a unique measure $\tau^{\vphantom{p}}_E\in M_1(E)$ with the following property
\begin{equation}
V^{\tau^{\vphantom{p}}_E}(x)\equiv\gamma^{\vphantom{p}}_E=\const,\quad x\in E.
\label{equrob}
\end{equation}
The measure $\tau^{\vphantom{p}}_E$ is called Robin's or the equilibrium measure for the compact set $E$, $\gamma^{\vphantom{p}}_E$ is Robin's constant for $E$. Let
\begin{equation*}
g_E(\zeta,\infty)=\log|\zeta|+\gamma^{\vphantom{p}}_E+o(1),\quad \zeta\to\infty,
\end{equation*}
be Green's function for the domain $D:=\myh{\CC}\setminus{E}$ with the logarithmic singularity at the infinity point $\zeta=\infty$. Then
\begin{equation}
g_E(\zeta,\infty)=\gamma^{\vphantom{p}}_E-V^{\tau^{\vphantom{p}}_E}(\zeta).
\label{green}
\end{equation}

Let $g_E(t,\zeta)$ be Green's function for the domain $D=\myh{\CC}\setminus{E}$ with the logarithmic singularity at a point $t=\zeta$,
\begin{equation}
G^\varkappa_E(\zeta)=\int_F g_E(t,\zeta)\,d\varkappa(t),\quad \zeta\in D,
\label{greenpote}
\end{equation}
be the corresponding Green's potential of a probability measure $\varkappa$ supported on $F$, $\varkappa\in M_1(F)$, $F=[a,b]$, $1<a<b$.

Similarly, let $g_F(x,\zeta)$ be Green's function for the domain $G:=\myh{\CC}\setminus{F}$ with the logarithmic singularity at the point $x=\zeta$,
\begin{equation}
G^\mu_F(\zeta)=\int_E g_F(x,\zeta)\,d\mu(x),\quad \zeta\in G,
\label{greenpotf}
\end{equation}
be the corresponding Green's potential of a probability measure $\mu$ supported on $E$, $\mu\in M_1(E)$.

It is proved~\cite{RaSu13} that there exists a unique measure $\lambda_E\in M_1(E)$ with the following property (cf.~\eqref{equrob}):
\begin{equation}
3 V^{\lambda_E}(x)+G_F^{\lambda_E}(x)\equiv\const,
\quad x\in \supp{\lambda_E}=E.
\label{mixequ}
\end{equation}
Let $\lambda_F=\beta_F(\lambda_E)\in M_1(F)$ be the balayage (see~\cite{Lan66}) of $\lambda_E$ from the domain $G$ onto its boundary $\partial G=F$.

For an arbitrary polynomial $Q\not\equiv0$ set
\begin{equation*}
\chi(Q):=\sum_{\eta:Q(\eta)=0}\delta_\eta,
\end{equation*}
where each zero $\eta$ of $Q$ is counted with accordance of its multiplicity.

Set $f(z):=f_*(1/z)\in\HH(0)$, $N=2n+1=3m+2$ and let $[n/n]_f(z)$ be the PA to $f(z)$ and $P_{2m,1}(z)/P_{2m,0}(z)$ be the rational type II HP approximation to $f(z)$ of order~$2m$.

The following statements hold true.
\begin{theorem}\label{the1}
Set $p_n(\zeta):=\zeta^nP_n(1/\zeta)$, $q_n(\zeta):=\zeta^nQ_n(1/\zeta)$.
We have that as $N\to\infty$
\begin{equation}
\frac1n\chi(p_n),\frac1n\chi(q_n)\overset{*}\longrightarrow
\tau^{\vphantom{p}}_E
\label{apweak}
\end{equation}
in the sense of the weak-$*$ convergence of measures,
and uniformly on the compact subsets of $D$
\begin{equation}
\bigl|f_*(\zeta)-[n/n]_f(1/\zeta)\bigr|^{1/N}\to
e^{-g_D(\zeta,\infty)}<1,\quad \zeta\in D.
\label{zhupade}
\end{equation}
\end{theorem}

\begin{theorem}\label{the2}
Set $p_{2m,0}(\zeta):=\zeta^{2m}P_{2m,0}(1/\zeta)$,
$p_{2m,1}(\zeta):=\zeta^{2m}P_{2m,1}(1/\zeta)$.
We have that as $N\to\infty$
\begin{equation}
\frac1{2m}\chi(p_{2m,0}),\frac1{2m}\chi(p_{2m,1})\overset{*}\longrightarrow
\lambda^{\vphantom{p}}_E
\label{hp2weak}
\end{equation}
in the sense of the weak-$*$ convergence of measures,
and uniformly on the compact subsets of $D$
\begin{equation}
\biggl|f_*(\zeta)-\frac{P_{2m,1}}{P_{2m,0}}(1/\zeta)\biggr|^{1/N}
\to e^{-\frac13G_E^{\lambda_F}(\zeta)-g_D(\zeta,\infty)}<
e^{-g_D(\zeta,\infty)},\quad \zeta\in D.
\label{zhuhp2conv}
\end{equation}
\end{theorem}

Theorem~\ref{the1} is in fact a consequence of Stahl's theorem on the convergence in capacity. The proof of Theorem~\ref{the2} will be given in the next paper.

Note that it follows from the Theorems~\ref{the1} and~\ref{the2} that when we are given $N = 2n + 1 = 3m +2$ coefficients of the power series \eqref{expan}, it is more efficient to use the rational type II HP-approximant corresponding to $m$ than to use the PA corresponding to~$n$.

\end{document}